\documentclass{amsproc}
\usepackage{amssymb}

\newtheorem{theorem}{Theorem}[section]
\newtheorem{lemma}[theorem]{Lemma}
\newtheorem{proposition}[theorem]{Proposition}
\newtheorem{corollary}[theorem]{Corollary}
\newtheorem{conjecture}[theorem]{Conjecture}

\theoremstyle{definition}
\newtheorem{definition}[theorem]{Definition}

\newtheorem{problem}[theorem]{Problem}
\newtheorem{question}[theorem]{Question}

\theoremstyle{remark}

\numberwithin{equation}{section}

\newcommand{\ds}{\displaystyle}

\newcommand{\cat}{\mbox{\rm cat}}
\newcommand{\sop}{{\bf Proof:} }
\newcommand{\eop}{\;\;\;  \Box}
\newcommand{\GF}{\Gamma_{\F}}

\newcommand{\cBF}{\cB(\F)}
\newcommand{\cGF}{\cG_{\F}}
\newcommand{\cRF}{\cR_{\F}}

\newcommand{\F}{{\mathcal F}}
\newcommand{\wF}{{\widehat{\F}}}

\newcommand{\G}{\Gamma}

\newcommand{\K}{{\mathbf K}}

\newcommand{\mE}{{\mathbb E}}
\newcommand{\mN}{{\mathbb N}}
\newcommand{\mR}{{\mathbb R}}
\newcommand{\mS}{{\mathbb S}}
\newcommand{\mT}{{\mathbb T}}
\newcommand{\mZ}{{\mathbb Z}}

\newcommand{\WF}{{\mathbf{W}(\F)}}
\newcommand{\NWF}{{\mathbf{NW}(\F)}}

\newcommand{\AF}{{\mathbf{A}_{\F}}}
\newcommand{\KF}{{\mathbf{K}_{\F}}}
\newcommand{\AR}{{\mathbf{A}_{\cR}}}
\newcommand{\KR}{{\mathbf{K}_{\cR}}}

\newcommand{\EF}{{\mathbf{E}_{\F}}}
\newcommand{\PF}{{\mathbf{P}_{\F}}}
\newcommand{\HF}{{\mathbf{H}_{\F}}}
\newcommand{\ER}{{\mathbf{E}_{\cR}}}
\newcommand{\PR}{{\mathbf{P}_{\cR}}}
\newcommand{\HR}{{\mathbf{H}_{\cR}}}

\newcommand{\BF}{{\mathbf{B}_{\F}}}
\newcommand{\SF}{{\mathbf{S}_{\F}}}
\newcommand{\FF}{{\mathbf{F}_{\F}}}
\newcommand{\BR}{{\mathbf{B}_{\cR}}}
\newcommand{\SR}{{\mathbf{S}_{\cR}}}
\newcommand{\FR}{{\mathbf{F}_{\cR}}}

\newcommand{\ZF}{{\mathbf{Z}_{\F}}}
\newcommand{\CF}{{\mathbf{C}_{\F}}}

\newcommand{\cA}{{\mathcal A}}
\newcommand{\cB}{{\mathcal B}}

\newcommand{\cG}{{\mathcal G}}

\newcommand{\cI}{{\mathcal I}}

\newcommand{\cL}{{\mathcal L}}
\newcommand{\cM}{{\mathcal M}}

\newcommand{\cO}{{\mathcal O}}
\newcommand{\cP}{{\mathcal P}}

\newcommand{\cR}{{\mathcal R}}
\newcommand{\cS}{{\mathcal S}}
\newcommand{\cT}{{\mathcal T}}
\newcommand{\cU}{{\mathcal U}}

\newcommand{\cW}{{\mathcal W}}

\newcommand{\bC}{{\mathbf C}}

\newcommand{\bE}{{\mathbf E}}

\newcommand{\bH}{{\mathbf H}}

\newcommand{\bP}{{\mathbf P}}

\newcommand{\bZ}{{\mathbf Z}}

\newcommand{\wtL}{{\widetilde L}}

\newcommand{\wtU}{{\widetilde U}}

\newcommand{\wtx}{{\widetilde{x}}}

\newcommand{\wty}{{\widetilde{y}}}

\newcommand{\whD}{{\widehat D}}

\newcommand{\whM}{{\widehat M}}

\begin{document}

\title{Classifying foliations}

\author{Steven Hurder}
\address{Department of Mathematics, University of Illinois at Chicago, 322 SEO (m/c 249), 851 S. Morgan Street, Chicago, IL 60607-7045, USA}
\email{hurder@uic.edu}
\thanks{The author was supported in part by NSF Grant \#0406254.}
\date{March 15, 2008 and, in revised form, October  20, 2008.}



\begin{abstract}
We give a survey of the   approaches to classifying foliations, starting with the Haefliger classifying spaces and the various results and examples about the secondary classes of foliations. Various dynamical properties of foliations are introduced and discussed, including expansion rate, local entropy, and orbit growth rates. This leads to a decomposition of the foliated space into Borel or measurable components with these various dynamical types. The dynamical structure is compared with the classification via secondary classes.

\end{abstract}

\maketitle
 
\vfill
\eject

\section{Introduction} \label{sec-intro}

A basic problem of foliation theory is    how to ``classify all the foliations'' of fixed codimension-$q$ on a given closed manifold $M$, assuming that at least one such foliation exists on $M$. 
This survey concerns this classification problem for foliations.

Kaplan \cite{Kaplan1941} proved the first complete classification result in the subject  in 1941. For a foliation $\F$ of the plane by lines (no closed orbits) the leaf space $\cT = \mR^2/\F$ is a (possibly non-Hausdorff) 1-manifold, and   $\F$ is characterized up to homeomorphism by the leaf space $\cT$. (See also    \cite{Haefliger2002b,HaefligerReeb1957,Wang1988}.) 
 Palmeira \cite{Palmeira1978} proved an analogous result for the  case of foliations of simply connected manifolds by hyperplanes.

Research on classification advanced dramatically in 1970, with three seminal works: Bott's Vanishing Theorem \cite{Bott1970},   Haefliger's construction of a ``classifying space'' for foliations  \cite{Haefliger1970,Haefliger1971},  and Thurston's profound  results on existence and classification of foliations \cite{Thurston1972,Thurston1974a,Thurston1974b,Thurston1976},  in terms of the homotopy theory of Haefliger's classifying spaces $B\G_q$.
The rapid progress during this period can be seen in the two survey works by H. Blaine Lawson: first was his  article   ``Foliations'' \cite{Lawson1974}, which gave a survey of the field up to approximately 1972; second was the   CBMS Lecture Notes \cite{Lawson1975} which included developments up to 1975, including the Haefliger-Thurston Classification results.  The work of many researchers in the 1970's filled in more details of this classification scheme, as we discuss below. 

The philosophy of the construction of $B\G_q$ is simple, as described by Haefliger  \cite{Haefliger1984}: for a codimension-$q$  foliation $\F$ of a manifold $M$,  one associates  a natural map $h_{\F} \colon M \to B\GF$ to a space $B\GF$ which is ``foliated'', with  all leaves in $B\GF$  contractible. The space  $B\GF$ represents a homotopy-theoretic version of the leaf space $M/\F$, much as one constructs the Borel quotient space $M_{G} = EG \times_G M$ for a Lie group action $G \times M \to M$ on a manifold.  Two foliations $\F_1$ and $\F_2$ are equivalent in this sense if their ``leaf spaces'' $B\G_{\F_1}$ and $B\G_{\F_2}$ are functorially homotopic. The universal space $B\G_q$ is   obtained by performing this operation on the universal groupoid $\G_q$. Imagine $B\G_q$  as the  direct limit of performing this classifying  construction on all foliations of codimension-$q$. This suggests just how large is this universal space.  This scheme of classification via the a canonical model of the  leaf space also underlies the classification of the $C^*$-algebras associated to   foliations in Connes' work \cite{Connes1994}.

In the approximately 40 years since Haefliger introduced the classifying spaces for foliations    in 1970, our knowledge of the homotopy theory of $B\G_q$ remains marginal. There is one exception, which happens when we ask about the classification of foliations which are transversely $C^1$. That is, their transverse holonomy maps are assumed to be $C^1$, and so they  are classified  by  a space $B\G_q^{1}$.  In 
this case, Tsuboi  \cite{Tsuboi1989a,Tsuboi1989b} proved in 1989 an absolutely remarkable result, that the natural map $\nu \colon B\G_q^{1} \to BO(q)$,  classifying the universal normal bundle of $C^1$-foliations, is a homotopy equivalence! For foliations whose transverse differentiability is $C^r$ with $r > 1$,  the study of the  homotopy type of the classifying spaces $B\G_q^{(r)}$ of $C^r$-foliations awaits a similar breakthrough.

During Spring semester 1982, there  was an emphasis  on foliation theory at the  Institute for Advanced Study, Princeton. In attendance were Paul Schweitzer, along with Larry Conlon, Andr\'{e} Haefliger, James Heitsch and the author among others. Lawrence Conlon had just arrived with the hand-written manuscript by 
G\'{e}rard Duminy, whose main result  and its proof were presented in seminar:

\begin{theorem}  [Duminy \cite{Duminy1982a,CantwellConlon1984,Hurder2002a}] \label{thm-duminy}
Let $\F$ be a $C^2$-foliation of codimension-one on a compact manifold $M$. If the Godbillon-Vey class $GV(\F) \in H^3(M)$ is non-trivial, then $\F$ has a resilient leaf, and hence $\F$ has an uncountable set of leaves with exponential growth.
 \end{theorem}
 
  Duminy's result  solved a conjecture   posed in 1974 by Moussu \& Pelletier \cite{MoussuPelletier1974} and Sullivan \cite{Schweitzer1978}: must   a $C^2$-foliation with non-zero Godbillon-Vey class   have a leaf with exponential growth type? More broadly, this conjecture  can be interpreted (and was)  as asking for connections between the values of the secondary classes and geometric or dynamical properties of the foliation. That some connection exists, between the dynamics of a foliation and its secondary classes, was suggested by the known examples, and by the philosophy that the secondary classes are evaluated on compact cycles in $M$, and if $\F$ does not have ``sufficiently strong recurrence'' or even ``chaotic behavior'' near such a cycle, then the secondary classes vanish  on it.

The method of proof of Duminy's Theorem, which was the culmination of several years of investigations by various researchers \cite{CantwellConlon1981c, DuminySergiescu1981, Herman1977, MMT1981, MoritaTsuboi1980, Nishimori1980, Tsuchiya1982,  Wallet1976},   suggested avenues of further research, whose pursuit during the past 26 years has led to a new understanding of foliations using ideas of dynamical systems and ergodic theory. The purpose of this paper is to survey some of these developments, with a highlight on some of the open questions.   Here is the  primary  question:

  \begin{question} To what extent do the dynamical and ergodic properties of a $C^r$-foliation on a closed manifold provide an  
  effective classification?  What aspects of the dynamical properties of a foliation $\F$   are determined by the homotopy class of the Haefliger  classifying map $h_{\F} \colon M \to B\G^r_q$?
 \end{question}
 
 The possibility of giving some solution to this problem depends, of course, on what we mean by ``classification''. 
 The literal answer is that ``classification''  is impossible to achieve, as non-singular vector fields on manifolds define foliations, and the   dynamical systems obtained from  vector fields are not ``classifiable'' in any reasonable sense. Thus, the much more complicated dynamics of foliations whose leaves may have dimension greater than one, are equally not classifiable. 
 
 On the other hand, it is possible to give broad descriptions of classes of foliations in terms of their dynamical and ergodic properties. For example, a foliation can be decomposed into its wandering and non-wandering components; or   into the union of leaves  with exponential and sub-exponential growth. This survey introduces  six   such decompositions of a foliation based on its dynamical properties (see \S\ref{sec-secclasses2})  including a new scheme, which has its roots in study of hyperbolic dynamical systems:

\begin{theorem} [Hurder \cite{Hurder2008b}] \label{thm-main1} 
Let $\F$ be a $C^1$-foliation on a closed manifold $M$. Then $M$ has a disjoint decomposition 
into  $\F$--saturated, Borel subsets:
\begin{equation}\label{decomp}
M = \EF \cup \PF  \cup \HF 
\end{equation}
\begin{itemize}
\item $\EF $ consists of ``elliptic leaves'' with ``bounded transverse expansion''
\item $\PF $ consists of ``parabolic  leaves'' with   ``slow-growth transverse expansion''
\item $\HF $ consists of ``(partially) hyperbolic  leaves'' with 
  ``exponential-growth transverse  expansion''.
\end{itemize}
\end{theorem}

The point of such a decomposition is to study the dynamical properties of the foliation $\F$ restricted to  each component, which then  suggests more focused  problems and approaches for further research. 
For example, a Riemannian foliation satisfies $M = \EF$; it is not known under what hypotheses  the converse is true. 
A  foliation is said to be {\it distal} if the orbits of pairs of distinct points remain a bounded distance apart under the action of the holonomy pseudogroup (see Definition~\ref{def-distal} below). A distal foliation 
satisfies $M = \EF  \cup \PF$. The class of foliations with $M = \EF  \cup \PF $ is called \emph{parabolic} in \S\ref{sec-parabolic}. Many classes of examples of parabolic foliations are known, but their full extent is not. 
 Finally, the set $\HF $ is the union of leaves which have some degree of ``non-uniformly partial hyperbolicity''. It is unknown in general what hypotheses are necessary   in order to conclude that $\F$ behaves chaotically on $\HF$, or  that the geometric entropy of $\F$ is positive on $\HF$.

The importance of the hyperbolic part of the decomposition $\HF$ is illustrated by the following generalization of 
Duminy's   Theorem~\ref{thm-duminy}:

\begin{theorem} [Hurder \cite{Hurder2008b}]  \label{thm-main2}
Let $\F$ be a $C^2$-foliation on a closed manifold $M$. Suppose that some residual secondary class $\Delta_{\F}^*(h_I   c_J) \in H^*(M;\mR)$ is non-zero. Then $\HF$ must have positive Lebesgue measure.
\end{theorem}
 If the codimension is one, then there is just one secondary class, the Godbillon-Vey class $GV(\F) = \Delta_{\F}^*(h_1c_1) \in H^3(M;\mR)$, and $GV(\F) \ne 0$   implies   the existence of resilient leaves, hence the existence of uncountably many leaves of exponential growth rate. Section~\ref{sec-q1} describes other results for   codimension-one foliations, which should admit extensions of some form to foliations with  codimension $q > 1$.
 
Note that these notes discuss only briefly the important topic of amenability for foliations,  in \S\ref{sec-amenable}. This is an important theme in the study of the dynamics and ergodic theory of foliations  \cite{Brooks1983,CFW1981,Heitsch1983,HM1990,HM1991,Kaimanovich2001,Plante1975b, PlanteThurston1976}.    We also omit all discussions of the topic of random walks on the leaves of foliations, and the properties of harmonic measures for foliations, which have proven to be a powerful tool for the study of foliation dynamics 
\cite{Candel2003,CanCon2003,DKN2007,Ghys1995,Kaimanovich1988,Kaimanovich1997,Kaimanovich2001}.
Other omissions and details of proofs from this survey are   developed more fully in \cite{Hurder2008b}.

 These notes are an expansion of a talk of the same title, given at the conference 
 ``Foliations, Topology and Geometry in Rio'', August 7, 2007, 
on the occasion of the 70th birthday of Paul Schweitzer. The author would like to thank the organizers for making this special event possible, and their   efforts at making this excellent meeting a success. Happy Birthday, Paul!

\section{Foliation groupoids}\label{sec-groupoid}

The defining property of a codimension-$q$ foliation $\F$ of a closed manifold $M$ is that locally, $\F$ is defined by a submersion  onto a manifold of dimension $q$. If the leaves of $\F$ form a fibration of $M$, then this local fibration property is global: there is a global submersion $\pi \colon M \to B$ onto a compact manifold $B$ whose fibers are the leaves of $\F$. In general, one only has that for each $x \in M$ there is some open $U_x \subset M$ and fibration $\pi_x \colon U_x \to B_x  \subset \mR^q$ such that the fibers of $\pi_x$ are connected submanifolds of dimension $p$, equal to some connected component of a leaf of $\F | B_x$.

This local submersion data defines a topological groupoid over $M$, denoted by $\GF$, whose object space is the disjoint union $B = \cup B_x$, and the morphisms are generated by local transformations $h_{\{x,y\}}$ defined whenever $U_x \cap U_y \ne \emptyset$. 
One can assume that $B_x = (-1,1)^q$ for all $x \in M$, and that one needs only work with a finite collection of open sets of $M$, enough to form a covering. The resulting groupoid $\GF$ is a compactly-generated subgroupoid of the groupoid $\G_q$ of local diffeomorphisms of $\mR^q$. The Haefliger classifying map   $h_{\F} \colon M \to B\G_q$ is defined   from this data.

Moerdijk observed in \cite{Moerdijk1991,Moerdijk2001c} that   the above data naturally defines an \'{e}tale groupoid, so one can form an associated  category of sheaves $\G(\cdot)$ on $M$ which admits a classifying topos $\cB\G_q$ for $\F$. More recently,  a third point of view has developed, that of a foliation as an example of a Lie groupoid   over $M$ which leads to much simplified formal constructions of the de~Rham and cyclic cohomology invariants of $\F$ \cite{Crainic2003,CrainicMoerdijk2000,CrainicMoerdijk2001, CrainicMoerdijk2004}. All approaches yield the same homotopy classification theory for foliations \cite{Moerdijk1991};  it is just a matter of taste how one defines these invariants. Our discussion here follows the original approach of Haefliger \cite{Haefliger1970,Haefliger1971,Haefliger1984}.

In this section, we make precise the objects being considered. Given the foliation $\F$, let $\cU = \{ \varphi_i \colon U_i \to (-1,1)^n \mid 1 \leq i \leq k\}$, be a covering of $M$ by foliation charts. That is, if we compose $\varphi_i$ with the projection $(-1,1)^n \to (-1,1)^q$ onto the last $q$-coordinates, where $n = p+q$, we obtain local submersions $\pi_i \colon U_i \to (-1,1)^q$ so that the fibers of $\pi_i$ are connected components of the leaves of $\F \mid U_i$. We assume that each chart $\varphi_i$ admits an extension to a foliation chart $\widetilde{\varphi}_i \colon \wtU_i \to (-2,2)^n$ where $\wtU_i$  contains the closure $\overline{U_i}$. Fix a Riemannian metric on $M$. Then we can also assume that   $\wtU_i$ is a convex subset of $M$ for the Riemannian distance function on $M$.

Define $\cT_i = \varphi_i^{-1}(\{0\} \times (-1,1)^q)$ and let $\cT = \cup \cT_i$ be the complete transversal for $\F$ associated to the covering.  Identify each $\cT_i$ with the the subset $(3i-1,3i+1)^q$ via a translation of the range of $\pi_i$.  Then  $\cT$ is identified with the disjoint union of open subsets of  $\mR^q$.  

A pair of indices $(i,j)$ is said to be \emph{admissible} if $U_i \cap U_j \ne \emptyset$. 

For $(i,j)$ admissible,   there is a local diffeomorphism $h_{j,i} \colon \cT_{i,j} \to \cT_{j,i}$ where 
$D(h_{j,i}) = \cT_{i,j} \subset \cT_i$ is the domain of $h_{j,i}$ and $R(h_{j,i}) = \cT_{j,i} \subset \cT_j$ is the range. The maps $\{h_{j,i} \mid (i,j) ~{\rm admissible}\}$ are the transverse change of coordinates defined by the foliation charts, and the assumptions imply that each map $h_{j,i}$ admits an extension to a compact subset of $\mR^q$. Hence, even though defined on open subsets of $\mR^q$ we have uniform estimates of these maps and their derivatives. They define      a \emph{compactly generated pseudogroup} on $\mR^q$:
 
\begin{definition}[Haefliger \cite{Haefliger2002a}]
A     pseudogroup of transformations $\cG$ of $\cT$ is \emph{compactly generated}     if there is 
\begin{itemize}
\item   a relatively compact open subset $\cT_0 \subset \cT$ meeting all orbits of $\cG$
\item  a finite set $\cG^{(1)} = \{g_1, \ldots , g_k\} \subset    \cG$  such that $\langle \cG^{(1)} \rangle =  \cG | \cT_0$; 
\item   $g_i \colon D({g_i}) \to R({g_i})$ is the restriction of  $\widetilde{g}_i \in   \cG$ with   
$\overline{D(g)} \subset D({\widetilde{g}_i})$.
\end{itemize}
\end{definition}

The foliation $\F$ is said to be $C^r$ if the maps $\cGF^{(1)} \equiv \{h_{j,i} \mid (i,j) ~{\rm admissible}\}$ are $C^r$, where $r = \ell + \alpha$, $\ell \geq 1$ is an integer and $0 \leq \alpha < 1$ is the H\"{o}lder modulus of continuity for the $C^{\ell}$ derivatives of the maps $h_{j,i}$.  The collection $\cGF^{(1)}$ of maps define  a   compactly generated pseudogroup acting on $\cT$,  denoted by  $\cG^r_{\F}$ when we need to emphasize that the degree of transverse regularity is $C^r$, and otherwise simply denoted by $\cGF$.

  The groupoid $\GF$ of $\F$ is the  space of germs associated to the elements of $\cGF$
$$\G({\cGF}) = \{[g]_x \mid g \in \cGF ~ \& ~ x \in D(g)\} ~ , ~ \GF = \G(\cGF)  $$
with source map $s[g]_x = x$ and range map $r[g]_x = g(x)$. Again, when we need to emphasize the degree of regularity, we write $\GF^r$ for the germs of maps in $\cG^r_{\F}$. Let 
\begin{equation}\label{eq-sr}
\GF^x = \{\gamma \in \GF \mid s(\gamma) = x\} ~ , ~
 \GF^{x,y} = \{\gamma \in \GF \mid s(\gamma) = x ~ , ~ r(\gamma) = y\}
\end{equation}

The equivalence relation defined by $\F$ on $\cT$ is the set
$$\cRF = \{(x,y)  \mid x \in \cT, y \in L_x \cap \cT\}$$

We note a fundamental convention used throughout this paper. Given maps $f \colon U \to V$ and $g \colon V \to W$, we write their composition as $g \circ f \colon U \to W$. Thus, $\GF$ can be thought of as an ``opposite functor'' from the category of admissible strings to maps
$$\{(i_1, i_2), (i_2,i_3), \ldots , (i_{k-1},i_k)\} \mapsto \gamma_{i_k, i_{k-1}} \circ \cdots \circ \gamma_{i_3 , i_2} \circ \gamma_{i_2, i_1}$$
Thus, for example, composition gives a map $\GF^{y,z} \circ \GF^{x,y}   \to \GF^{x,z}$. 
This convention does not arise when considering compositions of elements of $\cRF$, but is fundamental  when considering the local actions of maps in $\cGF$ on $\cT$.

  The objects $\cGF^r$,  $\GF^r$ and $\cRF$ associated to $\F$ are the primary sources of our understanding of both the topological and dynamical classification of foliations.

\section{Topological dynamics} \label{sec-top}
 
We recall here some of the basic concepts of topological dynamics, as applied to the case of the pseudogroup $\cGF$ acting on the complete transversal $\cT$. These  ideas   play a fundamental part in our understanding of the relationship between the secondary classes of $\F$ and its properties as a dynamical system. 

Some properties of topological dynamics require the full pseudogroup $\cGF$ for their definition and study,  while others are inherent to the equivalence relation $\cRF$. The focus on properties of the equivalence relation can be found in the early works of Dye \cite{Dye1959} and Mackey \cite{Mackey1963,Mackey1966},  reached its full development in the works of Krieger \cite{Krieger1976}, Feldman-Moore 
\cite{FeldmanMoore1977a} and Connes \cite{Connes1994}, and continues very actively in the study of Borel equivalence relations today (see for example, \cite{DJK1994,HjorthKechris2005}). Though the theme of this paper is   about  a   classification scheme based on the role of approximations to the holonomy maps in $\cGF$ by their derivatives, underlying many of the results are fundamental structure theorems for the Borel equivalence relation $\cRF$. We  recall below the  decomposition of  $\cRF$ into its Murray-von~Neumann types.

 For $x \in M$, let $L_x$ denote the leaf of $\F$ containing $x$.
 
  For  $x \in \cT$, the \emph{orbit of $x$} is the set
 \begin{equation}
 \cO(x) = L_x \cap \cT  = \GF \cdot x \equiv \{y = r(\gamma)  \mid  \gamma \in \GF^x \}
\end{equation}

A subset $E \subset \cT$ is \emph{saturated} if $x \in E$ implies $\cO(x) \subset E$; that is, $L_x \cap \cT \subset E$. 

Given a subset $E \subset \cT$ we define   the saturation of $E$, either in $M$ or in $\cT$, 
$$E_{\F} = \bigcup_{x \in E} ~ L_x ~ \subset M ~ {\rm and} ~ E_{\cR} = \bigcup_{x \in E} ~ \cO(x) ~ \subset \cT$$
Note that if $E$ is a Borel subset of $\cT$, then $E_{\F}$ is a Borel subset of $M$, and $E_{\cR}$ is a Borel subset of $\cT$.  Also, for Lebesgue measure on $M$ and $\cT$, the assumption that $\cGF$ is finitely generated by $C^1$ maps implies that if $E$ has Lebesgue measure 0 in $\cT$, then $E_{\F}$ has Lebesgue measure 0 in $M$, and $E_{\cR}$ likewise in $\cT$.

\begin{definition} $\cBF$ is the $\sigma$-algebra of saturated, Borel subsets of $\cT$. 
\end{definition}
Note that $\cBF$ defines a Borel structure on the ``quotient space'' of $\cT$ by the action of $\cGF$.
  Given $E \in \cBF$, the \emph{full   sub-equivalence relation}  on $E$ is 
\begin{equation}\label{eq-restriction}
\cRF^E   =   \{(x,y) \in \cRF \mid x \in E\}
\end{equation}
 
For a homeomorphism $f \colon N \to N$ of some space $N$, it is clear how to define fixed-points and periodic points. For groupoid dynamics, this is not so clear. We   use the following as a workable definition:
\begin{definition}
$x \in \cT$ is a \emph{periodic point} for $\cGF$ if there exists $g \in \cGF$ with $g(x) = x$ and   $[g]_x$ is non-trivial; that is, the leaf $L_x$ admits a non-trivial element of germinal holonomy at $x$.  Hence, every element of $\cO(x)$ is also a periodic point. 
\end{definition}
It is clear how to define transitive points:
 \begin{definition} Let $E \in \cBF$. Then
$x \in E \subset \cT$ is a \emph{transitive point for $E$} if  $\cO(x) \subset E$ is dense in $E$. That is, their closures in $\cT$ are equal: $\overline{\cO(x)} = \overline{E}$.  We say that $x$ is a transitive point in the case when  $E = \cT$.
\end{definition}
Note that it is possible for a point $x \in \cT$ to be both transitive and periodic; in the   Roussarie example \cite{GodbillonVey1971}, every leaf with holonomy has both properties.

  There are several notions of ``minimal sets'' used in the study of the dynamics of codimension-one foliations \cite{CanCon2000,CantwellConlon1981b,Hector1983,HecHir1981}.  For higher codimension foliations, Marzougui and Salhi introduced a notion of local minimal set in \cite{MarzouguiSalhi2003}. Recall that a compact saturated subset $Z \subset M$  is \emph{minimal} if it  admits no proper closed saturated subset.  Clearly, every leaf $L \subset Z$ must then be dense. An open saturated subset  $U \subset M$ is \emph{locally minimal} if  for every leaf $L \subset U$, the closure $\overline{L} = \overline{U}$.  Local minimal sets play a fundamental role in the study of the structure for codimension-one foliations    \cite{CantwellConlon1981b,CantwellConlon1984,CantwellConlon1988b,Hector1974,Hector1983}.
The following definition combines these notions for general Borel saturated subsets of $\cT$:
 \begin{definition}
  $E \in \cBF$ is   \emph{minimal} if every $x \in E$ is transitive in $E$. 
\end{definition}
If   $E \in \cBF$ is relatively compact, then $E \subset \cT$ minimal implies that $E_{\F}$ is minimal in the usual sense. If  $E \in \cBF$ is open, then $E_{\F}$ is a local minimal set. However, for general $E$, the notation ``minimal'' is an abuse of notation, as the condition does not satisfy any descending chain condition. To see this, simply observe that if   $E \in \cBF$ is minimal, and $x \in E$, then $E - \cO(x)$ is again minimal in the above sense. A better notation might be to call such sets ``totally transitive'', but this notation is already in use in the dynamics literature (cf. page 768, \cite{HasselblattKatok2002}).

The notions of wandering and non-wandering points are easily generalized:
 \begin{definition}
 $x \in \cT$ is   \emph{wandering} if there exists an open set $x \in U_x \subset \cT$ such that for all $g \in \cGF$ with $x \in D(g)$ and $[g]_x \ne Id$, then $g(U_x \cap D(g)) \cap U_x = \emptyset$. 

 The \emph{wandering set} 
$\WF = \{ x \in \cT \mid x  ~ {\rm wandering}\}$
 is an open saturated subset. 
 
 The \emph{non-wandering set} $\NWF = \cT \setminus \WF$ is a closed invariant set.
\end{definition}

Let $E$ be a compact minimal set, which is not a single orbit. Then  $E \subset  \NWF$.

Finally, we define the $\omega$-limit set of an orbit. 

\begin{definition}\label{def-omega}
The $\omega$-limit set of a point $x \in \cT$ is the  relatively  compact saturated subset 
$$ \omega(x) ~ = ~ \bigcap_{\stackrel{S \subset \cO(x)}{\#S < \infty}} ~ \overline{\cO(x) - S} ~ \subset ~ \cT$$
\end{definition}
Observe that if $E$ is compact minimal set, then $\omega(x) = E$ for all $x \in E$.
\begin{definition}\label{def-proper}
An orbit $\cO(x)$ is \emph{proper} if $\omega(x) \cap \cO(x) = \emptyset$. A leaf $L_x$ of $\F$ is proper if its orbit is proper. An orbit which is not proper, is said to be \emph{recurrent}. 
\end{definition}
A point in a minimal set  which is not a single orbit is always recurrent.

The notion of distal actions and proximal orbits also have natural generalizations to pseudogroup dynamics.
\begin{definition}\label{def-proximal}
We say that a pair $x \ne y \in \cT$ is \emph{proximal} for $\cGF$ if for all $\epsilon > 0$, there exists $g \in \cGF$ such that $x, y \in D(g)$ and 
$d_{\cT}(g(x), g(y)) < \epsilon$.
\end{definition}

\begin{definition}\label{def-distal}
We say that a pair $x \ne y \in \cT$ is \emph{distal} for $\cGF$  if there exists $\epsilon_{x,y} > 0$, so that for every  $g \in \cGF$ such that $x, y \in D(g)$ then  $d_{\cT}(g(x), g(y)) \geq \epsilon_{x,y}$. The pseudogroup $\cGF$ is said to be distal if   every $x\ne y \in \cT$ is distal.
\end{definition}

Given a     set $Z$, let $\chi_Z$ denote its characteristic function.  
The Riemannian metric on $M$ defines a volume form $dvol$. Define   the ``Lebesgue measure'' $\mu_L(Z)$ of $Z$ as 
$$\mu_L(Z) = ~  \int_Z ~ dvol ~ = ~ \int_M ~ \chi_Z ~ dvol$$ 
The transversal $\cT$ is identified with a subset of $\mR^q$ via the foliation coordinate charts. Let $d\vec{x}$  denote  the Lebesgue measure   on $\mR^q$. For $E \in \cBF$, define   
$$\mu_L(E)  ~ = ~ \int_E ~ d\vec{x} ~ = ~ \int_{\cT} ~ \chi_E ~ d\vec{x}$$
Of course, $\mu_L(E)$ in general depends upon the choice of the foliation coordinate charts, but the property $\mu_L(E) = 0$ is independent of the charts used.

A key concept for the study of the dynamics of a Borel map  $f \colon N \to N$ is the existence of  invariant and quasi-invariant measures  on $N$.   For foliation groupoids, there is a similar concept.
\begin{definition}
Let $\mu$ be a Borel measure on $\cT$ which is finite on compact   subsets. Then we say:
\begin{enumerate}
\item $\mu$ is $\cGF$-quasi-invariant if for all $g \in \cGF$ and Borel   subsets $E \subset D(g)$,  $\mu(g(E)) = 0$ if and only if $\mu(E) = 0$.
\item $\mu$ is $\cGF$-invariant if for all $g \in \cGF$ and Borel   subsets $E \subset D(g)$, then $\mu(g(E)) = \mu(E)$.
\end{enumerate}
\end{definition}
The following is immediate:
\begin{lemma}
 Let $\mu$ be a  $\cGF$-invariant , Borel measure on $\cT$ which is finite on compact   sets. 
Then  $\mu(\WF) =0$.
\end{lemma}

A set $E \in \cBF$ is $\mu$-null if $\mu(E) = 0$, and $\mu$-conull if $\mu(\cT\setminus E) = 0$.

We can now formulate the notion of ergodic sets for $\cGF$ (or equivalently, for $\cRF$) which is the measure-theoretic version of minimal sets.
 \begin{definition} Let $\mu$ be a  $\cGF$-quasi-invariant, $\sigma$-finite, Borel measure on $\cT$ which is finite on compact   sets. We say that  
 $E \in \cBF$ with    $\mu(E) > 0$ is \emph{ergodic with respect to $\mu$},   if for any  $E' \in \cBF$ with   $E' \subset E$,   either $\mu(E') = 0$ or $\mu(E \setminus E') = 0$.
 
 In the case where $\mu = \mu_L$ is Lebesgue measure, then we just say that $E$ is \emph{ergodic}.
\end{definition}
The following is then immediate from the definitions.
\begin{proposition}
Let $E \in \cBF$ be    ergodic with respect to $\mu$. For any open set $U \subset \cT$ with $\mu(U \cap E) > 0$, then $(U \cap E)_{\cR}$ is a set of full $\mu$-measure in $E$.  
\end{proposition}

Finally, we recall  the  decomposition of $\cRF$ into its Murray-von~Neumann types. (See Nghiem \cite{Nghiem1973,Nghiem1975}, Krieger \cite{Krieger1976}, Feldman-Moore \cite{FeldmanMoore1975,FeldmanMoore1977a, Moore1988}, Katznelson-Weiss \cite{KatznelsonWeiss1991}, or   Section 4.$\gamma$, pages 50--59 of Connes \cite{Connes1994}.)   Assume    $E \in \cBF$ is ergodic with respect to the  $\cGF$-quasi-invariant  measure $\mu$ on $\cT$. Recall from (\ref{eq-restriction}) the  full   sub-equivalence relation $\cRF^E$  on $E$. Then the measured equivalence relation  $(\cRF^E, \mu)$ can be classified into one of three broad categories:

$\bullet$ $(\cRF^E, \mu)$  has  ``Type I'' if  there exists a Borel   subset $E_0 \subset E$   such that  for $\mu$-almost every  $x \in E$, the orbit $\cO(x) \cap E_0$ contains precisely one point. 

$\bullet$ $(\cRF^E, \mu)$  has    ``Type II'' if  it is not Type I, and  the measure $\mu$ is  $\cGF$-invariant. 

$\bullet$ $(\cRF^E, \mu)$  has   ``Type III''  there is no     $\sigma$-finite, $\cGF$-invariant measure $\mu'$ on $E$ which is   absolutely continuous relative to $\mu$, other than  the zero measure. 

\medskip

For the equivalence relation  $\cRF$ on $\cT$ associated to the $C^1$-foliation $\F$, equipped with the standard Lebesgue measure on $\cT$, there is an ergodic decomposition into ergodic subequivalence relations, with induced absolutely continuous measures on each factor. Each ergodic factor can then be classified into one of the above three types, resulting in a \emph{measurable} decomposition of $\cT$, which is well defined up to sets of measure zero:
\begin{equation}\label{eq-MvN}
\cT ~ = ~ \cT_{I} ~ \cup ~  \cT_{II} ~  \cup ~ \cT_{III}  
\end{equation}
Let  $M_I, M_{II}, M_{III}$ denote the saturations of $\cT_I , \cT_{II} , \cT_{III}$ respectively, then we have the corresponding decomposition of $M$ into measurable, $\F$-saturated subsets
\begin{equation}\label{eq-MvN2}
 M ~ = ~ M_I ~  \cup ~ M_{II} ~  \cup ~  M_{III}
\end{equation}
 The Murray-von~Neumann type decomposition of $\cT$    thus obtained satisfies:
 
 The ``Type I'' component $\cT_{I}$   is the  largest saturated Borel measurable subset $E \subset \cT$ for which $\cRF^E$ is dissipative: That is, there exists a Borel   subset $E_0 \subset E$ such that for {a.e.} $x \in E$, the orbit $\cO(x) \cap E_0$ contains precisely one point. Moreover, the quotient space $E/\cRF$ is a standard, non-atomic Borel space.  For example, $\WF \subset \cT_I$.  Let $\BR$ denote the union of the finite orbits of $\cRF$, then $\BR \subset \cT_{I}$.
  
 The ``Type II'' component   $\cT_{II}$  is the largest saturated Borel measurable subset  $E \subset \cT \setminus  \cT_{I}$ such that no ergodic component is Type I, and $\cRF^E$ admits an absolutely continuous,  $\sigma$-finite, $\cGF$-invariant measure $\mu$ with almost every orbit $\cO(x)$ for $x \in E$ being $\mu$-essential (i.e., for every open neighborhood $x \in U$, we have  $\mu(U \cap E) > 0$.) 
 
The ``Type III'' component $\cT_{III}$  is the complement of $\cT_{I} \cup \cT_{II}$. Thus, $\cT_{III}/\cRF$  is a completely singular Borel measure space, and  the only absolutely continuous,  $\sigma$-finite, $\cGF$-invariant measure $\mu$ on $\cT_{III}$ is the zero measure.
 
There are further, finer partitions of the Murray-von~Neumann decomposition. 
For example, $\cT_{I}$ is the union of the  finite orbits and the  infinite orbits. Also, we   have
\begin{proposition}[\cite{Glimm1961}; Proposition~1.9, \cite{HK1987}, \cite{Millett1987}]\label{prop-proper}
Let $E \in \cBF$, and suppose that for every $x \in E$ the orbit $\cO(x)$ is proper. Then $E \subset \cT_{I}$.
\end{proposition}

The   set  $\cT_{III}$ can be further decomposed into the ergodic types of its flow of weights. 

The relations between the Murray-von~Neumann types of a foliation and its secondary classes were first  considered in the paper \cite{HK1987}. At that time,   it was hoped that the Type decomposition (\ref{eq-MvN2}) would provide an effective means to ``classify'' foliations via their ergodic theory properties. Unfortunately,  little progress has been made towards this goal, partly because the Type II and Type III components for foliations are so difficult to characterize.

\section{Classifying spaces}\label{sec-classifying}

 Let $\Gamma^r_q$ denote the universal groupoid defined by all $C^r$ local diffeomorphisms of open sets of $\mR^q$ to open sets of $\mR^q$. This has a classifying space denoted by 
   $B\Gamma^r_q$ which was first  introduced by  Andr\'{e}  Haefliger in 1970 \cite{Haefliger1970,Haefliger1971}. 
   
   Recall that two    codimension-$q$, $C^r$-foliations $\F_0$ and $\F$ on $M$ are \emph{concordant} is there is a 
   codimension-$q$ $C^r$-foliation $\F$ on the product space $M \times [0,1]$ which is transverse to the boundary $M \times \{0,1\}$ such that $\F \mid M \times \{i\} = \F_i$ for $i=0,1$. Concordance forms an equivalence relation on foliations. The first main result of the homotopy classification theory of foliations states: 
\begin{theorem} [Haefliger \cite{Haefliger1970}] \label{haefliger1} Let   $\F$ be a $C^r$-foliation of codimension-$q$ on a manifold $M$ without boundary. Then there exists  a well-defined, functorial  map $h_{\F} \colon M \to B\Gamma^r_q$ whose homotopy class is uniquely defined by $\F$. Moreover, the homotopy class of $h_{\F}$ depends only on the foliated concordance class of $\F$.
 \end{theorem}
In other words, the set of homotopy classes of maps $[M, B\Gamma^r_q]$ ``classifies'' the concordance classes of codimension-$q$  $C^r$ foliations on $M$.  The monograph by Lawson \cite{Lawson1975} gives an excellent overview of this theory.  Haefliger's works \cite{Haefliger1970, Haefliger1971, Haefliger1984} offer deeper insights into the construction and properties of the spaces $B\G_q^r$.

The tangent bundles to the leaves of $\F$  define a subbundle $F = T\F \subset TM$. The normal bundle to $\F$ is the  orthogonal complement $Q =  F^{\perp} \subset TM$. Thus, each foliation defines a splitting $TM = F \oplus Q$.

 The derivative of a $C^r$ germ  $[g]_x$ gives an element $D[g]_x \in GL(\mR^q)$. This yields a natural transformation from $\GF$ to $GL(\mR^q)$ and induces a map of classifying spaces, $\nu \colon B\GF^r \to BGL(\mR^q) \simeq BO(q)$. The composition $\nu_{Q} = \nu \circ h_{\F} \colon M \to BO(q)$ classifies the normal bundle to $\F$. A more precise statement of Theorem~\ref{haefliger1} is that $\F$ defines a lifting of the classifying maps for the vector bundles  $F$ and $Q$:
  
\begin{picture}(100,80)
\put(180,60){$BO(p) \times B\G^r_q$}
\put(100,45){$\nu_F \times h_{\F}$}
\put(130,33){\vector(3,2){30}}
\put(220,40){$id \times \nu$}
\put(210,50){\vector(0,-1){17}}
\put(110,20){$M$}
\put(140,27){$\nu_F \times \nu_Q$}
\put(140,22){$\vector(1,0){30}$}
\put(180,20){$BO(p) \times BO(q)$}
\end{picture}

   The second main result of the homotopy classification theory of foliations is Thurston's celebrated  converse to   Theorem~\ref{haefliger1}: 
\begin{theorem}[Thurston \cite{Thurston1974b,Thurston1976}] \label{thm-thurston1}
A lifting  $\nu_F \times h_{\F} \colon M \to BO_p \times B\Gamma^r_q$ of $\nu_F \times \nu_Q$ yields    a $C^r$-foliation $\F$ on $M$ with   concordance class    determined by $h_{\F}$.
 \end{theorem}
 
 Suppose that $F \subset TM$ a codimension-one subbundle with oriented normal bundle, then the map $\nu_Q$ is homotopic to a constant, so always admits a lift. Hence, $M$ admits a foliation $\F$ whose tangent bundle $T\F$ is homotopic to $F$. This is one  of the well-known implications of   Theorem~\ref{thm-thurston1}.  The method of proof of the existence of $\F$ gives few insights  as to its geometric or dynamical properties.

 The classifying map of the universal normal bundle, $\nu \colon B\G_q^r \to BO(q)$, has a homotopy fiber denoted   by $F\G_q^r$, or sometimes by  $B\overline{\G_q^r}$ in the literature.   The space $F\G_q^r$ classifies the  codimension-$q$, $C^r$-foliations equipped with a  given (homotopy type of) framing of the normal bundle.  
One of the   ``milestone'' results from the 1970's gives a partial understanding of the homotopy type of this fiber.
\begin{theorem}[Mather-Thurston, Haefliger \cite{Haefliger1970,Mather1973,Mather1975a,Mather1979,Thurston1974a}]\label{thm-MT}
For $r  \ne q+1$, the space $F\G_q^r$ is $q+1$-connected.
\end{theorem}
The outstanding problem  is to show: 
\begin{conjecture}
For $r \geq 2$, the space $F\G_q^r$ is $2q$-connected.
\end{conjecture}
Tusboi provided a complete solution to this conjecture  in the case of $C^1$-foliations:
\begin{theorem}[Tsuboi \cite{Tsuboi1989a,Tsuboi1989b}] \label{thm-tsuboi}
For $q \geq 1$, the space $F\G_q^1$ is weakly homotopic to a point. That is, the map 
  $\nu \colon B\G_q^1 \to BO(q)$ is a weak homotopy equivalence.
\end{theorem}
Theorem~\ref{thm-tsuboi}  is one of the most beautiful results of the 1980's in foliation theory, both in the simplicity and strength of its conclusion, and the methods of proof which combined dynamical systems results with sophisticated methods of the study of classifying spaces, building on the earlier works of Mather \cite{Mather1971a,Mather1973,Mather1979} and Tsuboi \cite{Tsuboi1984a,Tsuboi1985,Tsuboi1987}.

To make the   the homotopy classification theory of foliations effective, we must understand the homotopy type of the spaces $B\Gamma^r_q$ and $F\Gamma^r_q$ for $r > 1$. This remains   the most important    open problem in the field, after almost 40 years. In the next section, we discuss some of the results about the homotopy theory of $F\Gamma^r_q$ for $r \geq 2$, obtained using the theory  of secondary characteristic classes of foliations and constructions of explicit examples.

\section{Characteristic classes of foliations}\label{sec-charclasses}

The normal bundle $Q$ to a $C^r$-foliation $\F$, when restricted to a leaf $L_x$ of $\F$, has  a natural flat connection $\nabla^{L_x}$ defined on $Q \mid L_x \to L_x$. The  collection of these leafwise flat connections  define the  \emph{Bott connection}  $\nabla^{\F}$ on $Q \to M$, which need not be flat over $M$.  The connection data provided by $\nabla^{\F}$ can be thought of as a  ``linearization'' of the normal structure to $\F$ along the leaves, which varies $C^{r-1}$ in the transverse coordinates.  Thus, $\nabla^{\F}$  captures aspects of the data provided by the Haefliger groupoid $\GF^r$ of $\F$ --  it is a ``partial linearization'' of the highly nonlinear data which defines the homotopy type of $B\GF$. (The discussion in section~4.48 of \cite{KT1975a} provides some more insight on this point of view.) In this section, we discuss the applications of this partial linearization  to the study of the space $B\G^r_q$. 

The seminal observation was made by Bott around 1970. The cohomology ring $H^*(BO(q); \mR) \cong \mR[p_1 , \ldots , p_k]$ where $2k \leq q$, and $p_{j}$ has graded degree $4 j$.

 \begin{theorem}[Bott Vanishing \cite{Bott1970}] \label{BVT} Let $\F$ be a codimension-$q$, $C^2$-foliation.   Let  $\nu_{Q}   \colon M \to BO(q)$ be the classifying map for the normal bundle $Q$. Then $\nu_Q^* \colon H^{\ell}(BO(q); \mR) \to H^{\ell}(M; \mR)$ is the trivial map for $\ell > 2q$.
\end{theorem}
\sop Let $\nabla^{\F}$ denote a Bott connection on $Q$ with curvature 2-form $\Omega^{\F}$. Then the restriction of $\Omega^{\F}$ to each leaf $L$ of $\F$ vanishes, as the restricted connection is flat. Hence, the entries of the matrix of 2-forms $\Omega^{\F}$ must lie in the ideal of the de~Rham complex of forms, $\cI^*(M,\F)$,  generated by the 1-forms which vanish when restricted to leaves. It follows that all powers $(\Omega^{\F})^{\ell} = 0$ when $\ell > q$ as this is true for $\cI^*(M,\F)^{\ell}$. Now, by Chern-Weil theory we can calculate a de~Rham representative for each $\nu_Q^*(p_j) \in H^{4j}(M, \mR)$ in terms of the curvature matrix $\Omega^{\F}$, so the image of  $\nu_Q^*$   must vanish in degrees above $2q$. \hfill $\eop$

\medskip
To the best of the author's knowledge,  there is   no explicit construction    of a foliation for which 
$\nu_Q^* \colon H^{\ell}(BO(q); \mR) \to H^{\ell}(M; \mR)$ is non-trivial in the range   $q < \ell \leq 2q$. Morita observed in \cite{Morita1977} that  there exists a codimension-$2$, $C^2$-foliation  on a closed 4-manifold $M$  for which 
$p_1(Q) = \nu_Q^*(p_1) \in H^4(M; \mR)$ is non-zero. The existence is based on some  of the deepest results of Mather and Thurston \cite{Thurston1974a}.  (See Morita \cite{Morita1977}, and also \S4 of Hurder \cite{Hurder1981b}, for applications of this remark.  
Open problems related to this example  are discussed in \S15 of \cite{Hurder2003}). It is a measure of our lack of understanding of the geometry of foliations that no more concrete constructions have been obtained to illustrate Bott's Theorem in a positive direction. This is another example of the lack of understanding of the spaces $B\G^r_q$ for $r \geq 2$.
 
It is a remarkable observation that Theorem~\ref{BVT} is false for integral coefficients, and the counter-examples are provided by quite explicit foliations:
\begin{theorem}[Bott-Heitsch \cite{BottHeitsch1972}] \label{BH} The universal normal bundle map, 
\begin{equation}
\nu^* \colon H^{\ell}(BO(q); \mZ) \to H^{\ell}(B\G_q^r; \mZ)
\end{equation}
 is injective for all $\ell \geq 0$.
\end{theorem}
 \sop Let $\mT^k \subset SO(q)$ be a maximal compact torus. Let $\mT^k_{\delta}$ denote this continuous group considered with the discrete topology. There is a natural map $B\mT^k_{\delta} \to B\G_q$, where a cycle $f \colon N \to  B\mT^k_{\delta}$ corresponds to a flat-bundle foliation over $N$ via the natural action of $O(q)$ on $\mR^q$, hence we obtain $h_f \colon N \to B\G_q^r$.
 
 Cheeger-Simons Character Theory \cite{CheegerSimons1985,ChernSimons1974}  then implies that  the composition 
 $$\nu^* \colon H^{\ell}(BO(q); \mZ) \to H^{\ell}(B\G_q^r; \mZ) \to H^{\ell}(B\mT^k_{\delta}; \mZ)$$
is injective.  \hfill $\eop$
     
      The author has recently given a construction of smooth foliated manifolds of compact manifolds  which ``realize'' these cohomology classes \cite{Hurder2008a}. 

Theorem~\ref{BH} implies that for $q \geq 2$ and $r \geq 2$, the space $F\G_q^r$ does not have the homotopy type of a finite type CW complex. More is true, that a CW model for $F\G_q^r$ must infinitely many cells in all   dimensions $4\ell -1$ for $2\ell > q$. 
One may ask what properties of foliations do these cells ``classify''? This is unknown.

\section{Secondary characteristic classes}\label{sec-scc}

We turn now to the theory of the secondary characteristic classes for $C^2$-foliations. There are a variety of expositions on this topic; we select a few aspects of the theory to discuss, based on the author's preferences. The reader can confer with any of the following general references 
for more details:
\cite{BernshteinRozenfeld1972,BernshteinRozenfeld1973,Bott1972a,Bott1975c,Bott1976,Bott1978,BottHaefliger1972,CanCon2003,Fuks1973,Fuks1976,GF1968,GF1969b,KT1974c,KT1975a,KT1978a,Lawson1975,Morita2001,Pittie1976a}.
The reader is cautioned that notation  in these papers is not particularly consistent between various authors.

First, consider the case of   codimension-one foliations with  oriented normal   bundle. Then the foliation $\F$ is defined by a non-vanishing $1$-form $\omega$, so that 
$$T\F = \{\vec{v} \in TM \mid \omega(\vec{v}) = 0\}$$
For $r \geq 2$, the integrability of the distribution $T\F$ is equivalent by  the Froebenius Theorem to the condition that $\omega \wedge d\omega = 0$. That is, $d\omega =  \eta \wedge \omega$ for some $1$-form $\eta$. 

Given a vector field $\vec{Y}$  on $M$, let $\cL_{\vec{Y}}$ denote the Lie derivation operator. Applied to a vector field   $\vec{Z}$ we have $\cL_{\vec{Y}} \vec{Z} = [\vec{Y}, \vec{Z}]$. 
Recall also the Cartan formula: let $\psi$ be a 1-form on $M$. 
Then $\cL_{\vec{Y}}  = \iota_{\vec{Y}}  \circ d  + d \circ \iota_{\vec{Y}}$, or 
$$ \cL_{\vec{Y}} (\psi)(\vec{Z}) =   d\psi(\vec{Y}, \vec{Z}) + \cL_{\vec{Y}} (\psi(\vec{Z})) $$
The Bott connection $\nabla^{\F}$  for $\F$ has a direct interpretation in terms of the Lie operator. 
Let $\vec{X}$ be the vector field on $M$ with values in $Q = T\F^{\perp}$ such that $\omega(\vec{X}) = 1$.  
Let $\vec{Y}$ be a vector field on $M$ with values in $T\F$. Then 
\begin{equation}
\nabla^{\F}_{\vec{Y}} (\vec{X}) = \omega(\cL_{\vec{Y}}(\vec{X})) \cdot \vec{X} = \omega([\vec{Y}, \vec{X}]) \cdot \vec{X}
\end{equation}
That is, the Bott connection acts along leaves as the Lie derivative operator, projected to the normal bundle. The vanishing of the curvature of the Bott connection along leaves is due to the Jacobi identity for vector fields tangent to the leaves. 
$$  d\omega(\vec{Y}, \vec{X}) = \eta \wedge \omega(\vec{Y}, \vec{X}) = \eta(\vec{Y}) $$
$$ 0 = \cL_{\vec{Y}} (\omega(\vec{X})) =  \cL_{\vec{Y}} (\omega)(\vec{X}) -  \omega( \cL_{\vec{Y}} ( \vec{X}))$$
so that by the Cartan formula, 
$$   d\omega(\vec{Y}, \vec{X}) = \cL_{\vec{Y}} (\omega)(\vec{X}) - \cL_{\vec{Y}} (\omega(\vec{X}) )   = \cL_{\vec{Y}} (\omega)(\vec{X}) =  \omega([\vec{Y}, \vec{X}]) $$
hence $\eta(\vec{Y}) = \omega([\vec{Y}, \vec{X}])$.
Thus, the $1$-form $\eta$   is the Bott connection 1-form for the   normal line bundle to $\F$,  with respect to the framing of the normal bundle $Q$ defined by the section $\vec{X}$. Furthermore,  
 $\eta(\vec{Y})$ measures the normal expansion of the normal field $\vec{X}$ under parallel transport by the leafwise vector field $\vec{Y}$.   The curvature of the connection $\eta$ is   $d\eta$. 
 
 A similar ``naive'' interpretation of the Bott connection and its curvature can be given in arbitrary codimension (see Shulman and Tischler \cite{ShulmanTischler1976}).

Define $h_1 = \frac{1}{2 \pi} \eta \in \Omega^1(M)$ and $c_1 = \frac{1}{2\pi} d\eta \in \Omega^2(M)$.
   
\begin{theorem} [Godbillon-Vey \cite{GodbillonVey1971}]   Let $\F$ be a codimension-one, $C^2$ foliation on $M$ with trivial normal bundle. Then the  3-form $h_1 \wedge c_1$ is closed, and the cohomology class  
 $GV(\F) = [h_1 \wedge c_1] = \frac{1}{4\pi^2} [\eta \wedge d\eta] \in H^{3}(M; \mR)$ is independent of all choices. Moreover, $GV(\F)$ depends only on the concordance class of $\F$.
 \end{theorem}
 
The   paper \cite{GodbillonVey1971} also included an example by Roussarie to show this class was non-zero. As is well-known, soon afterwards, Thurston gave  a construction of families of examples of foliations on the 3-sphere $\mS^3$, such that the $GV(\F) \in H^3(\mS^3; \mR) \cong \mR$  assumed a continuous range of real values. As a consequence, he deduced:

\begin{theorem} [Thurston \cite{Thurston1972}]   \label{thm-thurstonq=1}
For $r \geq 2$,  $\pi_3(F\Gamma^r_1)$ surjects onto $\mR$. 
\end{theorem}

The appendix by Brooks to   \cite{Bott1978} gives a clear and concise explanation of the examples Thurston constructed in his very brief paper \cite{Thurston1972}.

As the space $F\Gamma^{\infty}_1$ is known to be $2$-connected by Mather \cite{Mather1971a}, Theorem~\ref{thm-thurstonq=1} implies that a CW-complex model for $F\Gamma^r_1$ must have an uncountable number of $3$-cells. Morita  asked in \cite{Morita1985} whether the cup product map $H^3(F\Gamma^{\infty}_1) \otimes H^3(F\Gamma^{\infty}_1) \to H^6(F\Gamma^{\infty}_1)$ is non-zero; the answer is not known.  Tsuboi further studied this problem in \cite{Tsuboi1998} in the PL setting. Here is another simple question of this type:
    \begin{question}
  Suppose that   $M^n$ is a closed manifold of dimension $n > 3$ and admits a non-vanishing vector field, and $H^3(M; \mR)$ is non-trivial. Does $M$ admit a codimension-one, $C^2$-foliation $\F$ with $GV(F) \ne 0$? 
    \end{question}
    This is simply a question of whether we can find a continuous map $h_{\F} \colon M \to F\Gamma^{\infty}_1$ such that $h_{\F}^*(GV) \ne 0$. Nothing is known about this question, unless strong assumptions are made about the topological type of $M$.

  In the case where the normal bundle $Q$ is not trivial,     the defining $1$-form $\omega$ is still well-defined ``up to sign'', hence the $1$-form $\eta$ is well-defined up to sign, and hence the product $GV(\F) = \frac{1}{4\pi^2} [\eta \wedge d\eta] \in H^{3}(M; \mR)$ is well-defined.

  The Seminaire Bourbaki article by Ghys \cite{Ghys1989} is a basic reference for the properties of the Godbillon-Vey class; the author's survey  \cite{Hurder2002a} is a more recent update. There are many  ``classic'' papers on the subject: \cite{
 Brooks1979,
 BrooksGoldman1978,
 CantwellConlon1984,
 Duminy1982a,
 Duminy1982b,
 Ghys1987a,
 Herman1977,
 HK1990,
 MMT1981,
 MMT1983,
 Mitsumatsu1985b,
 Thurston1972,
 Tsuboi1990,
 Tsuboi1992a,
 Tsuboi1992b}.

For  foliations of codimension $q \geq 2$, there are two related theories of secondary classes, corresponding to the cases where the normal bundle $Q$ is trivial, or not. If $Q$ is trivial, then $\F$ along with a choice of framing $s$ of $Q$, defines a classifying map $h_{\F, s} \colon M \to F\G_q^r$. When $Q$ is not trivial,   the classifying map is denoted by  $h_{\F} \colon M \to B\G_q^r$. The additional data of a framing for $Q$ yields more secondary invariants, while implying that all of the Pontrjagin characteristic classes   vanish.

The    construction  of secondary classes for foliations followed several paths during the rapid  development of the subject in the  early 1970's. Kamber and Tondeur \cite{KT1974c,KT1975a}    constructed the secondary classes     in terms of the  truncated Weil algebra $W^*(\mathfrak{gl}(q,\mR), O(q))_q$. (As a doctoral student of Kamber, the author has an innate respect for the power of this more formal approach.)   For simplicity of exposition, we will take the approach in \cite{Bott1972a,BottHaefliger1972,Godbillon1974}, which defines the secondary classes directly using the differential graded algebra (or DGA) model    $WO_q  \subset W^*(\mathfrak{gl}(q,\mR), O(q))_q$.

Denote by $I(\mathfrak{gl}(q, \mR))$ the graded ring of adjoint-invariant polynomials on the Lie algebra $\mathfrak{gl}(q, \mR)$ of the real general linear group $GL(q, \mR)$.  As a ring, $I(\mathfrak{gl}(q, \mR)) \cong \mR[c_1, c_2, \ldots , c_q]$  is a polynomial algebra on $q$ generators, where the $i^{th}$-Chern polynomial $c_i$ (with polynomial degree $i$ and graded degree $2i$) is defined by the relation
$$\det ( t \cdot Id - \frac{1}{2\pi} \cdot A  ) = \sum_{i=1}^q~ t^{q-i} c_i(A)$$
for $A \in \mathfrak{gl}(q,\mR)$ where  $Id$ the identity matrix. Let $I(\mathfrak{gl}(q, \mR))^{(q+1)}$ denote the ideal of polynomials of degree greater than $q$, and introduce    the quotient ring,
$$I(\mathfrak{gl}(q, \mR))_{q} = I(\mathfrak{gl}(q, \mR))/I(\mathfrak{gl}(q, \mR))^{(q+1)} \cong \mR[c_1, c_2, \ldots , c_q]_{2q}$$  isomorphic to the   polynomial ring truncated in graded degrees larger than $2q$.  

 The  Lie algebra cohomology of the Lie algebra $\mathfrak{gl}(q,\mR)$ has DGA model
$$H^*(\mathfrak{gl}(q, \mR)) \cong \Lambda(h_1, h_2, \ldots, h_{q}) $$
where $h_i$ denotes the ``transgression class'' of the Chern polynomial $c_i$, so has graded degree $2i-1$. It satisfies the differential identity 
 $d(h_i \otimes 1) = 1 \otimes c_i$ as noted below.
 
The    Lie algebra cohomology relative to the   group $O(q)$ of orthogonal matrices    has   DGA model 
$$H^*(\mathfrak{gl}(q, \mR), O(q)) \cong \Lambda(h_1, h_3, \ldots, h_{q'}) $$
where   $q'$ is the greatest  odd  integer $\leq q$.

The secondary classes  for foliations (whose normal bundle is not assumed to be trivial) arise from the DGA complex
$$WO_q  = \Lambda(h_1, h_3, \ldots, h_{q'}) \otimes \mR[c_1, c_2, \ldots , c_q]_{2q}$$
where   $q'$ is the greatest odd  integer $\leq q$. The differential on  $WO_q$ is defined by the relations $d(h_i \otimes 1) = 1 \otimes c_i$ and $d(1\otimes c_i) = 0$.

The monomials 
$\ds h_I \wedge c_J =  h_{1_1} \wedge \cdots h_{i_{\ell}} \wedge c_1^{j_1} \cdots c_q^{j_q}$ such that 
\begin{equation}\label{eq-Veybasis}
i_1 <  \cdots<  i_{\ell} ~, ~ |J| = {j_1} + 2 j_2 + \cdots + q  {j_q} \leq q ~, ~i_1 + |J| > q
\end{equation}
are closed in $WO_q$, and they span  the cohomology $H^*(WO_q)$ in degrees greater than $2q$. The {\it Vey basis} of $H^*(WO_q)$ is a subset of these (cf. \cite{BottHaefliger1972,Godbillon1974,Heitsch1978,Lawson1974}).

Let $\nabla^{\F}$ denote a Bott connection on $Q$ with curvature 2-form $\Omega^{\F}$.  
Chern-Weil theory yields the characteristic DGA homomorphism $\Delta_{\F} \colon I(\mathfrak{gl}(q, \mR)) \to \Omega^*(M)$, where  $\Delta_{\F}(c_i) = c_i(\Omega^{\F}) \in \Omega^{2i}(M)$. For $i = 2j$,  the closed differential form  $\Delta_{\F}(c_{2i})$ of degree $4i$ is a representative of the     $i^{th}$ Pontrjagin class  $p_i(Q)$ of $Q$. 

Extend $\Delta_{\F} $ to a DGA homomorphism $\Delta_{\F} \colon WO_q \to \Omega^*(M)$ where $\Delta_{\F}(h_{2i-1}) \in \Omega^{4i-1}(M)$ is a transgression class for $\Delta_{\F}(c_{2i-1})$.  
The induced map in cohomology, $\ds \Delta^*_{\F} \colon H^*(WO_q) \to H^*(M; \mR)$,   depends only  on the concordance  class of $\F$.

For the case of codimension-one foliations, 
    $\Delta_{\F}(h_1) = \frac{1}{2\pi} \eta  \in \Omega^1(M)$ is the Reeb class introduced before, and so 
    $\Delta_{\F}(h_1\otimes c_1) = \frac{1}{4\pi^2} \eta \wedge d\eta \in \Omega^3(M)$ represents the Godbillon-Vey class.

 For   foliations of  codimension greater than one, the    secondary classes of $\F$ are spanned by the   images $\ds \Delta^*_{\F}(h_I \wedge c_J)$, where $h_I \wedge c_J$ satisfies (\ref{eq-Veybasis}). 
 

When  the normal bundle $Q$ is trivial, the choice of a framing, denoted by $s$, enables the definition of additional secondary classes.
Define the DGA complex 
$$W_q  = \Lambda(h_1, h_2, \ldots, h_{q}) \otimes \mR[c_1, c_2, \ldots , c_q]_{2q}$$
where $\mR[c_1, c_2, \ldots , c_q]_{2q}$  is the truncated polynomial algebra, truncated in graded degrees greater than $2q$.   The monomials 
$\ds h_I \wedge c_J =  h_{1_1} \wedge \cdots h_{i_{\ell}} \wedge c_1^{j_1} \cdots c_q^{j_q}$
satisfying (\ref{eq-Veybasis}) are closed, and they span  the cohomology $H^*(W_q)$ in degrees greater than $2q$.

The data $(\F, s,\nabla^{\F})$  determine a map of differential algebras  $ \Delta_{\F,s} \colon W_q \to \Omega^*(M)$.
The induced map in cohomology, $\ds \Delta_{\F,s}^{*} \colon H^*(W_q) \to H^*(M)$,   depends only on the homotopy class of the framing $s$ and the framed concordance class of $\F$.

A monomial $h_I \wedge c_J \in WO_q$ or $W_q$ is said to be {\it residual} if the degree of the Chern component  $c_J$ is $2q$. That is, if $|J| = q$.    These are the classes that define generalized measures on the $\sigma$-algebra $\cBF$ (see \S\ref{sec-weilmeas}). A special case of these are the {\it generalized Godbillon-Vey classes}, of the form 
$\Delta_{\F}^*(h_1 \wedge c_J)$.  The usual Godbillon-Vey class is $GV(\F) \equiv \Delta_{\F}^*(h_1 \wedge c_1^q) \in H^{2q+1}(M; \mR)$.

There are  natural restriction   maps $R \colon WO_{q+1} \to WO_q$ and $R \colon W_{q+1} \to W_q$. The images of these maps in cohomology with degree greater than $2q$ are called the \emph{rigid secondary classes}, so called because they are constant under $1$-parameter deformations of the given foliation \cite{Heitsch1973,Heitsch1978}. The only known examples of foliations with non-zero rigid classes are a set of examples constructed by the author   in  \cite{Hurder1981b,Hurder1985a} using homotopy methods. 
The examples realize classes in the image of 
$R^* \colon H^*(W_{q+1}) \to H^*(W_q)$ -- no examples are known of foliations for which the classes in the image of 
$R^* \colon H^*(WO_{q+1}) \to H^*(WO_q)$ are non-trivial.

The above constructions are ``functorial'', hence induce universal characteristic maps. This is described very nicely in Lawson \cite{Lawson1975}:

\begin{theorem}   Let $q \geq 1$ and $r \geq 2$. There are    well-defined   characteristic maps 
\begin{eqnarray*}
\Delta^*   \colon    H^*(WO_q) & \to &  H^*(B\G^r_q) \\
\Delta_s^*   \colon    H^*(W_q)  & \to &  H^*(F\G^r_q) 
\end{eqnarray*}
whose constructions   are   ``natural''. That is, given a codimension-$q$, $C^r$-foliation $\F$, the classifying map $\ds \Delta^*_{\F} \colon H^*(WO_q) \to H^*(M; \mR)$ satisfies the universal property:

\begin{picture}(100,60)\label{comm_diag}
\put(170,50){$H^*(B\G^r_q ; \mR) $}
\put(190,40){\vector(0,-1){17}}
\put(130,23){\vector(3,2){30}}
\put(200,30){$h_{\F}^*$}
\put(130,35){$\Delta^*$}
\put(150,18){$\Delta_{\F}^*$}
\put(100,10){$H^*(WO_q)  ~ \longrightarrow ~ H^*(M; \mR)  $}
\end{picture}

A similar conclusion holds for foliations with framed normal bundles.
\end{theorem}
 The study of the universal maps $\Delta^* $ and $\Delta_s^*$   has been the primary source of information, beyond Theorem~\ref{thm-MT}, about the (non-trivial) homotopy types of ~ $B\G^r_q$ and $F\G^r_q$  ~ for $r \geq 2$. The outstanding problem is to show:
 \begin{conjecture} \label{conj-injective}
 For  $q \geq 2$ and $r \geq 2$,    the   maps $\Delta^* $ and $\Delta_s^*$ are injective.
 \end{conjecture}

There was some hope that Conjecture~\ref{conj-injective} had been proved in 1977, based on a construction given by Fuks \cite{Fuks1977a,Fuks1977b,Fuks1978}. The basic idea was to start with the fact that the continuous cohomology version of the universal maps is known to be injective \cite{Bott1975b,Haefliger1979}, so one ``only needs to construct appropriate homology cycles'' to detect these continuous cohomology classes, in a fashion similar to the situation for the locally homogeneous  examples discussed below.   
It remains an open problem whether the method of proof sketched out in these notes by Fuks   can be filled in.

Next, we survey some of the explicit constructions of foliations for which the characteristic maps 
$\Delta_{\F}^*$ and $\Delta_{\F,s}^*$ are non-trivial.  
There are two general methods which have been employed, along with a few exceptional approaches.

The original example of Roussarie \cite{GodbillonVey1971}, and its extensions to   codimension $q > 1$, start with a semi-simple Lie group $G$.  Choose   closed subgroups $K \subset H \subset G$, with $K$ compact.  Then $G/K$ is foliated by the left cosets of $H/K$.
Choose a cocompact, torsion-free lattice $\G \subset G$, then the  foliation of $G/K$ descends to   a foliation $\F$ on the compact manifold  $M = \G\backslash G/K$ which is a    locally homogeneous space. The calculation of the secondary invariants for such foliations  then follows from explicit calculations in Lie algebra cohomology, using    Cartan's approach to the cohomology of homogeneous spaces.  Examples of this type are studied in 
\cite{Baker1978a,Fuks1976,KT1974b,KT1974d, KT1975a, KT1978a,KT1979, Pelletier1983,Pittie1976a, Pittie1976b, Tabachnikov1984, Tabachnikov1985, Tabachnikov1986, Yamato1975}. 
 For example,  Baker  shows in  \cite{Baker1978a}:
 \begin{theorem} Let $q = 2m > 4$. Then the set of classes 
$$\{h_1h_2h_{i_1}\cdots h_{i_k} \wedge c_1^{q} ; h_1h_2h_{i_1}\cdots h_{i_k} \wedge c_1^{q-2} c_2   \mid 2<i_1<\cdots<i_k\leq m\}$$
in $H^*(W_{q})$ map under $\Delta^*$ to linearly independent classes in $H^*(F\G_q; \mR)$.
 \end{theorem}
 The non-vanishing results of Kamber and Tondeur follow a similar format, but are more extensive, as given in Theorem~7.95 of \cite{KT1975a} for example. The conclusion of all these approaches is to show that the    universal maps $\Delta^* $ and $\Delta_s^*$ are injective on various subspaces of $H^*(WO_q)$  and $H^*(W_q)$.

 The second approach to constructing foliations with non-trivial secondary classes uses the method of ``residues''. 
The concept of a residue dates back to Grothendieck; its application to foliations began with vector-field residue theorems of Bott \cite{Bott1967a, Bott1967b},  Baum and Cheeger \cite{BaumCheeger1969},  and Baum and Bott \cite{BaumBott1970}. 
  Heitsch developed the residue theory for smooth foliations \cite{Heitsch1977, Heitsch1978, Heitsch1981, Heitsch1983} which was essential to his calculations of the non-vanishing of the secondary classes for codimension $q \geq 3$.

 We describe the basic idea of the    construction  of the Heitsch examples, simplifying somewhat (see \cite{Heitsch1988} for a nice description of these examples.) Again,  start with a  semi-simple Lie group $G$,   a   compact subgroup $K  \subset G$, and   a cocompact, torsion-free lattice $\G \subset G$. We require one more piece of additional data: a representation $\rho \colon \G \to GL(q+1, \mR)$ such that the induced action of $\G$ on $\mR^{q+1}$ commutes with the flow $\varphi_{\lambda}$ on $\mR^{q+1}$ of a ``radial'' vector field $\vec{v}_{\lambda}$ on $\mR^{q+1}$ where $\lambda$ is some multi-dimensional parameter. The simplest example might be to let $\vec{v}_{\lambda} = r \partial/\partial r$ be the standard Euclidean radial vector field, but any vector field which commutes with the action, vanishes at the origin, and  the quotient space $\mR^{q+1}/ \varphi_{\lambda} \cong \mS^{q}$ will work. 
  
  Form the associated flat bundle $\mE = (G/K \times \mR^{q+1})/\G$, which has a foliation $\F_{\rho}$ whose leaves are coverings of $G/K$. The codimension of this foliation is $q+1$.  Now     form the   quotient manifold $M = \mE/\varphi_{\lambda}$ which is diffeomorphic to an $\mS^q$-bundle over $B = \G \backslash G/K$. The foliation $\F_{\rho}$ descends to a foliation denoted by $\F_{\lambda}$ on $M$.  
  
  The diffeomorphism class of the quotient manifold $M$ is independent of $\vec{v}_{\lambda}$, but the foliations $\F_{\lambda}$ need not be. The secondary classes of $\F_{\lambda}$ are calculated using the residues at the zero set of the zero-section of $\mE \to B$ for the induced vector field $\vec{v}_{\lambda}$ on $\mE$. Note that by assumption, this zero set equals the vanishing locus of $\vec{v}_{\lambda}$.   By various clever choices of the vector field $\vec{v}_{\lambda}$ and groups $K \subset G$, one then obtains that the secondary classes of $\F_{\lambda}$ are non-zero.   Moreover, and perhaps the fundamental point, is that collections of secondary classes for the family $\F_{\lambda}$ of foliations  can vary independently with the multi-variable parameter $\lambda$.   One note about this construction, is that it works starting with codimension $q \geq 3$.
  
 For the case of $q =2$,   Rasmussen \cite{Rasmussen1980} modified  the construction by Thurston of codimension-one foliations with varying Godbillon-Vey class \cite{Thurston1972}. Thurston's construction used the weak-stable foliation of the geodesic flow on a compact Riemann surface (with boundary) of constant negative curvature. Rasmussen extended these ideas to the case of   compact hyperbolic 3-manifolds.  He showed there exists families $\{\F_{\lambda} \mid \lambda \in \mR\}$  of smooth foliations in codimension-$2$ for which the secondary    classes $\Delta^{\F_{\lambda}}(h_1 \wedge c_1^2),  \Delta^{\F_{\lambda}}(h_1 \wedge c_2) \in H^5(M; \mR)$ vary continuously and independently.

  Together, these examples yield that for a fixed collection of classes in the image of $\Delta^*$, there are continuous families of cycles $h_{\F_{\lambda}} \colon M \to B\G_q$ such that the evaluation of these fixed secondary classes on the cycles defined by the $\F_{\lambda}$ varies continuously. Thus, $H_*(B\G_q; \mZ)$ must be a truly enormous integral homology group! 
 
 The Thurston examples  give $1$-parameter family of foliations $\F_{\lambda}$ on $\mS^3$,  for which the evaluation map  $GV(\F_{\lambda}) \colon H^3(\mS^3; \mZ) \to \mR$ has continuous image.  However, the   constructions for codimension $q > 1$ discussed above yield   information on the groups $H_*(B\G_q; \mZ)$, but no direct information about the homotopy groups $\pi_*(B\G_q)$.

 There is a second approach to constructing examples of foliations with non-trivial secondary classes,   developed by the author \cite{Hurder1981a,Hurder1981b,Hurder1985a, HLe1990}. The idea is  to use knowledge of the homotopy theory of $F\G^r_q$ to deduce from  the construction of examples above, the existence of   ``classifying maps'' $h_{\F,s} \colon \mS^n \to F\G^r_q$ for which $\Delta_{\F,s}$ has  ``continuously varying'' non-trivial images in $\pi_*(B\G_q)$. 
Then by Thurston's Existence Theorem~\ref{thm-thurston1}, one  concludes that $h_{\F,s}$ is homotopic to the  classifying map of some foliation $\F$ on $M$. Hence, one obtains  for codimension $q > 1$,  that there  exists   families of foliations $\{\F_{\lambda}\}$ on spheres $\mS^n$, where $n \geq 2q+1$, with continuously varying secondary classes, and also continuously varying Whitehead products of these classes \cite{Haefliger1978b, Hurder1981a, SchweitzerWhitman1978}. This method provides absolutely no insight into the geometry or dynamics  of the foliations $\F_{\lambda}$ so obtained.
 
Here are two  typical results:  let $\{h_1 \wedge c_1^2, h_1 \wedge c_2\} \subset H^5(W_2)$ denote the Vey basis for degree $5$.
\begin{theorem} [Theorem~2.5,  \cite{Hurder1985a}]\label{thm-hurder1}
Evaluation of the universal cohomology classes $\{ \Delta^*(h_1 \wedge c_1^2), \Delta^*(h_1 \wedge c_2)\}  \subset H^5(B\G_2 ; \mR)$ on the image $\pi_5(B\G_2) \to H_5(B\G_2 ; \mZ)$ defines a surjection of abelian groups
\begin{equation}
\{ \Delta*(h_1 \wedge c_1^2), \Delta*(h_1 \wedge c_2) \} \colon \pi_5(B\G_2) \to \mR \oplus \mR
\end{equation}
\end{theorem}

 For each $q > 1$, there exists sequences of non-negative integers $\{v_{q,\ell}\}$ (defined in \S2.8 of \cite{Hurder1985a}) with the properties:
 \begin{enumerate}
\item For $q=2$, $\ds \lim_{k\to \infty} v_{2, 4k+1} = \infty$, with $v_{2, 4k+1} > 0$ for all $k > 0$;
\item For $q=3$, $\ds \lim_{k\to \infty} v_{3, 3k+1} = \infty$, with $v_{3, 3k+1} > 0$ for all $k > 0$;
\item For $q>3$, $\ds \lim_{\ell \to \infty} v_{q, \ell} = \infty$.
\end{enumerate}
\begin{theorem} [Hurder \cite{Hurder1981a, Hurder1985a}]
For $q > 1$ and $\ell > 2q$, there is an epimorphism of abelian groups 
\begin{equation}
h^* \colon \pi_{n}(B\Gamma^r_q) \to \mR^{v_{q, \ell}}
\end{equation}
The maps $h^*$ are defined in terms of the dual homotopy invariants of \cite{Hurder1981a}, and the integers 
$\{v_{q,\ell}\}$ are the ranks of various free graded Lie algebras in the minimal model for the truncated polynomial ideal $ \mR[c_1, c_2, \ldots , c_q]_{2q}$. 
\end{theorem}
 One conclusion of all these results is   that the secondary  classes   measure some     uncountable aspect    of  foliation geometry. What that is, remains to be determined; but without a doubt, the homotopy groups $\pi_{*}(B\Gamma^r_q)$, $r \geq 2$,  are fantastically large.
 
  What is striking, looking back at the roughly 10 years between 1972 and 1982 during which this problem was actively researched, is how limited the types of examples discovered proved to be.  
 All of the ``explicit constructions'' of foliations in the literature with non-trivial secondary characteristic classes are either locally homogeneous, or deformations of locally homogeneous actions;  essentially, they  are all   generalizations and/or modifications of the original example of Roussarie for the Godbillon-Vey class in codimension-one.

Certain special classes of foliations  have their own theory of secondary invariants, as well as dynamical properties. For example, 
 a foliation $\F$ is said to be \emph{Riemannian} \cite{Haefliger1971,Reinhart1959, Reinhart1983,Molino1988,Molino1994} if there is a Riemannian metric on the transversal $\cT$ of section~2 for which all of the holonomy transformations $h_{j,i}$ are isometries.  All of the   secondary classes in $H^*(W_q)$ vanish for Riemannian foliations with framed normal bundles. However, if the truncation degree used in the construction of the DGS complex     $W_q$  is reduced from $2q$ to the exact codimension, $q$, then there is a modified construction of secondary classes for Riemannian foliations with framed normal bundles, yielding a map of complexes $\Delta \colon H^*(RW_q) \to H^*(M;\mR)$ \cite{KT1974c, KT1975a,LazarovPasternack1976a,LazarovPasternack1976b,Morita1979,Pasternack1975,Yamato1979a,Yamato1981}. 
The author proved  that the associated universal map of secondary invariants for Riemannian foliations is injective in  \cite{Hurder1981c}.  The work of the author with T\"{o}ben   \cite{HT2008} establishes relations between the values of these classes and the dynamical properties of Riemannian foliations.

 A foliation $\F$ is said to be \emph{transversally holomorphic}  if there is an integrable complex structure   on the transversal $\cT$ of section~2 for which all of the holonomy transformations $h_{j,i}$ are holomorphic. This class of foliations was introduced by Haefliger in \cite{Haefliger1971}, and    properties of their classifying spaces have been studied    by Adachi \cite{Adachi1991},
  Haefliger {\it et al}  \cite{HaeSid1982,HaeSun1985}, and Landweber \cite{Landweber1974}.   The theory of secondary classes for transversally holomorphic foliations is much richer   than for ``real'' foliations, as pointed out by Kamber-Tondeur \cite{KT1975a}. Rasmussen \cite{Rasmussen1978} and the author \cite{Hurder1982} gave non-vanishing results for various subsets of their secondary invariants. More recently, Asuke has studied their secondary classes in much greater depth, and also related to the values of certain of the secondary invariants to the dynamical properties of the foliations \cite{Asuke2000,Asuke2001,Asuke2003a,Asuke2003b,Asuke2003b,Asuke2004}. 

Transversally conformal and transversely projective  foliations provide yet another subclass of foliations, whose characteristic classes have been investigated, along with  their specialized dynamical properties 
\cite{Abe1979,Asuke1996,Asuke1998,BensonEllis1985, Nishikawa1981,NishikawaTakeuchi1978,  Takeuchi1980,Tarquini2004,Vaisman1979,Yamato1979b}.

\section{Localization and the Weil measures} \label{sec-weilmeas}

Localization is a property of the residual secondary classes,  apparently unique to the theory of  characteristic classes for foliations, and a distinct phenomenon  from residue theory.  In essence, it states that if $h_I \wedge c_J \in WO_q$  is a residual class  and $E \in \cBF$, then there is a well-defined restriction $\Delta^*_{\F}(h_i \wedge c_J)|E \in H^*(M; \mR)$. Moreover, this restriction is countably additive, and  vanishes if $E$ has Lebesgue measure zero. Finally,  the value of $\Delta^*_{\F}(h_i \wedge c_J)|E$ can be estimated using the dynamical and ergodic theory properties of $\F | E$.  Localization principles first appeared in the study  of the Godbillon-Vey class for codimension-one foliations in the 1970's, and was a key point in Duminy's proof of Theorem~\ref{thm-duminy}.  See \S2 of \cite{Hurder2002a} for a survey of its development for codimension-one foliations; details are in  Heitsch-Hurder  \cite{HH1984}.

Assume that both $M$ and the normal bundle $Q$ are oriented, so there exists   a positively oriented,  decomposable $q$-form $\omega$ on $M$ which defines $\F$. Let $A^*(M,\F)$ denote the   ideal in $\Omega^*(M)$ generated by $\omega$ : $A^{q+\ell}(M,\F) = \{\omega \wedge \psi \mid  \psi \in \Omega^{\ell}(M)\}$.

By the Frobenius Theorem,    $A^*(M,\F)$ is a differential ideal, whose  cohomology   is denoted by
$ H^{*}(M,\F) = H^*(A^*(M, \F), d)$. 
 Let $[M] \in H_n(M; \mZ)$ denote  the  fundamental class of $M$. 

Recall from \S\ref{sec-scc} that the secondary classes arise from the cohomology of the complex 
$$WO_q  = \Lambda(h_1, h_3, \ldots, h_{q'}) \otimes \mR[c_1, c_2, \ldots , c_q]_{2q}$$
   The choice of a Bott connection on $Q$ defines a DGA map $\Delta_{\F} \colon WO_q \to  \Omega^*(M)$.
The idea of the Godbillon and Weil  functionals is to separate the roles of the forms $\Delta_{\F}(h_i)$ and $\Delta_{\F}(c_i)$ in the definition of $\Delta(h_I \wedge c_J)$, and then study   the special properties of  the forms $\Delta_{\F}(h_I)$.

The first basic result is as follows. Given $E \in \cBF$, let $\chi_E$ denote its characteristic function. 
For each   monomial 
$h_I \in  \Lambda(h_1, h_3, \ldots, h_{q'})$, with degree $\ell$,  let $\phi \in A^{n-\ell}(M,\F)$ and $E \in \cBF$, define
\begin{equation}\label{eq-weildef}
\chi_E(h_I) [\phi] = \int_E ~ \Delta_{\F}(h_I) \wedge \phi  \equiv  \int_M ~ \chi_E \cdot  \Delta_{\F}(h_I) \wedge \phi
\end{equation}

\begin{theorem}[Heitch-Hurder \cite{HH1984}] \label{thm-HH} Let $\F$ be a $C^r$ foliation for $r \geq 2$. Suppose that $d\phi = 0$. Then $\ds \chi_E(h_I) [\phi] $ depends only on the cohomology class $[\phi] \in H^{n-\ell}(M, \F)$, and is independent of the choice of the Bott connection $\nabla^{\F}$. 

That is, for each $E \in \cBF$ there is a well-defined continuous linear map
$$\chi_E(h_I) \colon H^{n-\ell}(M, \F) \to \mR$$

Moreover, the correspondence $E \mapsto \chi_E(h_I)$ defines a countably additive measure
\begin{equation}
\chi(h_I) \colon \cBF \to {\rm Hom}_{cont}(H^{n-\ell}(M, \F), \mR) \equiv H^{n-\ell}(M, \F)^*
\end{equation}
which vanishes on sets $E \in \cBF$ with Lebesgue measure zero.
\end{theorem}
The proof of Theorem~\ref{thm-HH} uses the Leafwise Stokes' Theorem, Proposition~2.6 of \cite{HH1984},  and basic techniques of  Chern-Weil   theory. 
  \begin{definition}
  $\chi(h_I)$ is called the \emph{Weil measure}  on $\cBF$ corresponding to $h_I \in WO_q$. The Godbillon measure on $\cBF$ is the functional $g_E = \chi_E(h_1)$.
  \end{definition}
 
 Theorem~\ref{thm-HH} enables us to define the localization of the residual secondary classes. Let $h_I \wedge c_J \in WO_q$ be a residual class, hence    $|c_J| = 2q$.  
 Then the closed form $\Delta_{\F}(c_J) \in A^{2q}(M, \F)$, and $\Delta_{\F}(h_I \wedge c_J) \in A^{2q+\ell}(M, \F)$ where    $h_I$ has degree $\ell$.
 
 Let $\psi \in \Omega^{n - 2q-\ell}(M)$ be a closed form, then $\chi_E \cdot \Delta_{\F}(c_J)  \wedge \psi$ is a closed Borel form of top degree on $M$, and its integral over the fundamental class $[M]$ equals 
 $\chi_E(h_I)[\Delta_{\F}(c_J)  \wedge \psi]$, which  depends only on $[\psi] \in H^{n - 2q-\ell}(M; \mR)$.

\begin{definition} \label{def-localization}
The localization of $\Delta^*_{\F}(h_I \wedge c_J) \in H^{2q+\ell}(M; \mR)$ to $E \in \cBF$ is the cohomology class  $\Delta^*_{\F}(h_I \wedge c_J)|E  \in H^{2q+\ell}(M; \mR)$ defined by Poincar\'e Duality for $M$ and the linear functional
$$ [\psi] \mapsto \chi_E(h_I)[\Delta_{\F}(c_J)  \wedge \psi] ~ , ~ [\psi] \in H^{n - 2q-\ell}(M; \mR)$$
\end{definition}

We can now state the localization principle for the residual secondary classes.
 \begin{corollary} \label{cor-localization}
 Suppose that $\{E_{\alpha} \in \cBF \mid \alpha \in \cA\}$ form a  disjoint countable decomposition of $M$ into foliated Borel subsets. Then for each residual class $h_i \wedge c_J$ we have
 \begin{equation}\label{eq-localization}
\Delta^*_{\F}(h_I \wedge c_J) = \sum_{\alpha \in \cA} ~ \Delta^*_{\F}(h_I \wedge c_J)|E_{\alpha}
\end{equation}
In particular, if   $\chi_{E_{\alpha}}(h_I)= 0$ for each $\alpha \in \cA$,  then $\Delta^*_{\F}(h_I \wedge c_J)  = 0$.
 \end{corollary}
 The formula (\ref{eq-localization}) demands that we ask, what determines the values of the terms $\Delta^*_{\F}(h_I \wedge c_J)|E_{\alpha}$ in the sum?  One partial answer  is that the dynamics of $\F |E_{\alpha}$ gives estimates for the values of the Weil measures $\chi_{E_{\alpha}}(h_I)$.

\section{Foliation  time and distance} \label{sec-time}

     A continuous dynamical system on a compact manifold $M$ is a flow $\varphi \colon M \times \mR \to M$, where the orbit $L_x = \{\varphi_t(x) = \varphi(x,t) \mid t \in \mR\}$ is thought of as the time trajectory of the point $x \in M$.  The trajectories   are points, circles or lines immersed in $M$, and ergodic theory is the study of their aggregate and statistical behavior.

In foliation dynamics, the concept of time-ordered trajectories is replaced by   multi-dimensional futures for points, the leaves of $\F$. Ergodic theory  of foliations  asks for properties of the    aggregate and statistical behavior of the collection of its leaves. One of the key points in the development of foliation dynamics in the 1970's  was  the  use of  leafwise distance  as a substitute for the ``time''  in a dynamical system defined by a flow or a map; distance along a leaf  measures how far we can get along a leaf in a given time, hence provides a  substitute for ``dynamical time''.

Recall that the     pseudogroup  $\cGF$ is generated by the set  $\cGF^{(1)} = \{h_{j,i} \mid (i,j) ~{\rm admissible} \}$.

 \begin{definition}\label{def-gamma}
 For $Id \ne \gamma \in \GF^x$,    the word length $\| \gamma \|_x$    is the least $m$ so that 
$$\gamma  = [g_{i_m} \circ \cdots \circ g_{i_1}]_x  ~ {\rm where} ~ {\rm each} ~ g_{i_{\ell}} \in \cGF^{(1)}$$
If $\gamma$ is the   germ of the identity map at $x$, then set $\| \gamma \|_x = 0$.
\end{definition}
Word length is a measure of the ``time'' required to get  from $x$ to a point $y \in \cO(x)$ following a path which has the same germinal holonomy at $x$ as  $\gamma$. Thus, even if $x = y$,  so that $x$ is a fixed-point for the action of some $g \in \cGF$ where $\gamma = [g]_x$, the ``time'' required to get from $x$ to $y$ need not be  zero.  

Let $\wtL_x$ denote the holonomy covering of the leaf $L_x$ through $x$, endowed with the Riemannian metric lifted from the induced metric on $L_x$.  Then the factorization $\gamma = [g_{i_m} \circ \cdots \circ g_{i_1}]_x$ as above defines a piecewise differentiable path in $\wtL_x$ from the lift $\wtx \in \wtL_x$ to a well-defined point $\wty \in \wtL_x$. 
Let  $\sigma_{\wtx, \wty} \colon [0,1] \to \wtL_x$ be the distance-minimizing geodesic   from $\wtx$ to $\wty$. Denote its   length   by 
$\| \sigma_{\wtx, \wty} \|$. 
\begin{proposition}\label{prop-dist1}
There exists constants $0 < C_1 \leq C_2$, which depend only on the choice of Riemannian metric on $M$ and the foliation covering $\cU$,  such that 
$$C_1 \cdot \| \sigma_{\wtx, \wty} \| ~ \leq ~ \| \gamma \|_x ~ \leq ~ C_2 \cdot  \| \sigma_{\wtx, \wty} \|$$
\end{proposition}
That is, the word length on    $\GF^x$ is quasi-isometric to  the geodesic length function on the holonomy cover $\wtL_x$.

There is an alternate notion of ``time'' for the leaves of a foliation, based on the distance function on leaves. 
Recall that the equivalence relation defined by $\F$ on $\cT$ is the set
$\cRF = \{(x,y)  \mid x \in \cT, y \in L_x \cap \cT\}$. 
In essence, the equivalence relation forgets the information of which leafwise path is taken from $x$  to $y$, and uses only that there is some path. There is a natural map $s \times r \colon \GF \to \cRF$, and the fiber $(s \times r)^{-1}(x,x) = \GF^{x,x}$ is  the holonomy group of the leaf $L_x$ at $x$, by definition.

 \begin{definition}\label{def-relation}
 For $(x,y) \in \cRF$,  the   distance $d_{\cR}(x, y)$    is the least $m$ so that 
$$y  = g_{i_m} \circ \cdots \circ g_{i_1}(x)  ~ {\rm where} ~ {\rm each} ~ g_{i_{\ell}} \in \cGF^{(1)}$$
Set $d_{\cR}(x, x) = 0$ for all $x \in \cT$.
\end{definition}

Given $(x,y) \in \cRF$ let $\sigma_{x,y} \to L_x$ be the distance-minimizing leafwise geodesic   from $x$ to $y$. Set $d_{\F}(x,y) = \|\sigma_{x,y}\|$. 
We then have an estimate in terms of the leafwise distance function $d_{\F}$, using the same constants as in Proposition~\ref{prop-dist1}.
\begin{proposition}\label{prop-dist2}
There exists constants $0 < C_1 \leq C_2$, which depend only on the choice of Riemannian metric on $M$ and the foliation covering $\cU$,  such that 
$$C_1 \cdot d_{\F}(x,y)  ~ \leq ~ d_{\cR}(x,y) ~ \leq ~ C_2 \cdot d_{\F}(x,y) $$
\end{proposition}

The   distance functions on $\GF$ and $\cRF$ are related by:
\begin{proposition}\label{prop-dist3} Let $(x,y) \in \cRF$, then
$$d_{\cR}(x,y) = \inf ~ \{ \| \gamma \| \mid \gamma \in \GF^{x,y} \}$$
In particular, if $L_x$ is a leaf without holonomy,  then $d_{\cR}(x,y) = \| \gamma\|_x$  for $\gamma \in \GF^{x,y}$.
\end{proposition}

Both notions of  ``foliation time'' -- Definitions~\ref{def-gamma} and  \ref{def-relation} -- appear  in the literature. Plante's definition of  the growth of leaves in foliations 
used  the     leafwise distance function  \cite{Plante1972a,Plante1973b,Plante1975a,Plante1975b}. The     leafwise distance function also is used in the work by   Connes, Feldman and Weiss \cite{CFW1981} on their study of   amenable equivalence relations, and  in the work   by the author with Katok \cite{HK1987}. 

On the other hand, the word metric on $\GF$ is crucial, for example,  in the study of amenable groupoids, as in Anantharaman-Delaroche and Renault \cite{ADR2000}. The groupoid $\GF$ appears naturally in the study of non-commutative geometry associated to a foliation \cite{Connes1994}, and here again the groupoid metric plays a fundamental role.

\section{Orbit growth and the F{\o}lner condition} \label{sec-growth}

One of the most basic invariants of foliation dynamics is the growth rates of orbits. References for this section are Plante's original article \cite{Plante1975b}, and   \S1.3 of  \cite{HK1987} for the properties of metric equivalence relations. The survey \cite{Hurder1994} discusses    quasi-isometry invariants for foliations in a much broader context.

Given $x \in \cT$ and a positive integer $\ell > 0$, let 
\begin{equation}
B_{\cR}(x, \ell) = \{ y \in \cO(x) \mid d_{\cR}(x,y) \leq \ell\}
\end{equation}
The first remark is that it does not matter here whether we use the distance function $d_{\cR}(x,y)$ or the the norm $\| \gamma \|_x$ on $\GF^{x,y}$ to define the ``balls'' $B_{\cR}(x, \ell)$; Proposition~\ref{prop-dist3} implies that we get the same sets. Thus, the notation with subscript $\cR$ is justified -- these sets are inherent for  the equivalence relation $\cRF$.

 \begin{definition}  \label{def-growthexp}
 ${\rm Gr}(\cR, x, \ell) = \# B_{\cR}(x, \ell) $ is 
the \emph{growth function}  of $x \in \cT$.

The \emph{growth rate function} on $\cT$ is defined by:
\begin{equation}
{\rm gr}(\cR, x) = \limsup_{\ell \to \infty} ~ \frac{\ln \{ {\rm Gr}(\cR, x, \ell) \}}{\ell}
\end{equation}
\end{definition}

Recall that 
$\cU = \{ \varphi_i \colon U_i \to \cI^n \mid 1 \leq i \leq k\}$ is our covering of $M$ by foliation charts. Introduce the number 
$${\cat}(\cGF) = \max_{1 \leq i \leq k} ~ \# \{U_j \mid   ~U_i \cap U_j \ne \emptyset \}$$
Note that   ${\rm Gr}(\cR, x, \ell) \leq {\cat}(\cGF)^{\ell}$, hence ${\rm gr}(\cR, x) \leq \ln \{{\cat}(\cGF)\}$. 

Also note that  ${\rm Gr}(\cR, x, \ell)$ is a bounded function exactly when $\cO$ is a finite set, in which case   ${\rm Gr}(\cR,y, \ell)$ is bounded for all $y \in \cO(x)$ and the leaf $L_x$ is   compact. 

 \begin{proposition} \label{prop-growth}  The growth rate function satisfies:
\begin{enumerate}
\item  For all $(x,y) \in \cRF$, ${\rm gr}(\cR, x) = {\rm gr}(\cR, y)$. 
\item ${\rm gr}(\cR, x)$ is a Borel function on $\cT$, hence is a.e. constant on ergodic subsets. 
\item There is a disjoint Borel decomposition
$$ \cT = \BR ~ \cup ~ \SR ~ \cup ~ \FR ~ {\rm where} $$
\begin{itemize}
\item $\BR = \{ x \in \cT \mid \# \cO_{\cR}  < \infty\}$  (bounded ~ orbits)
\item $\SR = \{ x \in \cT - \cB \mid {\rm gr}(\cR, x) = 0\}$  (slow~orbit ~ growth)
\item $\FR = \{ x \in \cT \mid {\rm gr}(\cR, x) > 0\}$  (fast~orbit ~ growth)\\
\end{itemize}
\item 
The saturations of these sets yields a disjoint Borel decomposition
$$M = \BF ~ \cup ~ \SF ~ \cup ~ \FF$$
which is independent of the choices of Riemannian metric on $M$, and the choice of foliation charts   for $\F$. Moreover:  $\BF$  consists of the compact leaves of $\F$;  $\SF $ consists of non-compact leaves with subexponential growth; 
 $\FF$  consists of leaves with (possibly non-uniform) exponential growth.
\end{enumerate}
\end{proposition}

Note  that none of the three Borel sets $\BR$,  $\SR$ or  $\FR$ need be closed. In fact, a key (and usually difficult) problem in foliation dynamics is   to understand how the closures of these sets intersect.

When the foliation $\F$ is defined by a flow,   every leaf either is a closed orbit, hence of bounded growth, or is a line so has linear, hence subexponential growth. The presence of leaves with exponential growth rates for foliations can only occur when the  leaves have dimension at least two.     It was observed already in various works in the 1960's that the existence of leaves of exponential growth for foliations make the study of foliation dynamics ``exceptional'' \cite{Hurder1991a,Plante1974,Plante1975a,RosenbergRoussarie1970c,Sacksteder1965,SackstederSchwartz1964}. 

Given two increasing  functions of the natural numbers, $f,g \colon \mN \to [0,\infty)$, we say that $f \lesssim g$ if  there exists 
$A > 0$ and $B \in \mN^+$ such that 
$$ f(n) \leq A \cdot g(B\cdot n) ~, ~ {\rm for ~ all} ~ n \in \mN$$
This defines an equivalence relation, and the equivalence class of $f$, denoted by $[f]$,  is called its \emph{growth type}. 

The \emph{growth type}   of $x \in \cT$ is the growth type of the function ${\rm Gr}(\cR, x, \ell)$. The function $x \mapsto [{\rm Gr}(\cR, x, \ell)]$ is  constant on each orbit $\cO(x)$. This invariant of orbits (and correspondingly of leaves of foliations)  was introduced by Hector  \cite{Hector1977} (see also Hector-Hirsch \cite{HecHir1981}). 
A remarkable construction of Hector \cite{Hector1977} yields a  foliation of codimension-one on a compact 3-manifold for which  the leaves    have a continuum of growth types!  

While the growth rate function ${\rm gr}(\cR, x)$ is constant on ergodic components of $\cR$,   the properties of the growth type function  are not well-understood. For example, given $x \in \cT$, one can ask how the growth type function behaves when restricted to   the closure of its orbit, $\overline{\cO(x)}$. 
For codimension-one, $C^2$-foliations, the relation between growth types and orbit closures has been  extensively studied  \cite{CantwellConlon1978,CantwellConlon1981b,CantwellConlon1982,CantwellConlon1983,Hector1977,Hector1978,Plante1973b,Plante1975a,PlanteThurston1976,SackstederSchwartz1965,Tsuchiya1979a,Tsuchiya1979b,Tsuchiya1980a,Tsuchiya1980b}.
 
Finally, we mention that there is yet another growth invariant one can associate to   the growth function,   the   \emph{polynomial growth rate} (see \S1.3 of \cite{HK1987}):
 
 \begin{definition}  \label{def-growthpoly} 
The \emph{polynomial growth rate function} on $\cT$ is defined by:
\begin{equation}
{\rm p}(\cR, x) = \limsup_{\ell \to \infty} ~ \frac{\ln \{ {\rm Gr}(\cR, x, \ell) \}}{\ln \{ \ell \}} \leq \infty
\end{equation}
The orbit $\cO(x)$ has \emph{polynomial growth rate} of degree ${\rm p}(\cR, x)$ if  ${\rm p}(\cR, x) < \infty$.
Let 
$$\bP(\cR,r) = \{x \in \cT \mid {\rm p}(\cR, x) \leq r\}$$
\end{definition}
 The polynomial growth rate function has properties similar to that of the growth rate function: it is a Borel function, constant on orbits, and hence $a.e.$ constant on ergodic components of $\cR$. 
 Clearly, if ${\rm p}(\cR, x) < \infty$ then $x \in \BR \cup \SR$. 
 If $\F$ is defined by the locally-free action of a connected, nilpotent  Lie group $G$ on $M$, then ${\rm p}(\cR, x) \leq {\rm p}(G)$ for all $x \in \cT$, where ${\rm p}(G)$ is the integer     polynomial growth rate of $G$ for a left-invariant metric. 
The   examples in \cite{Hector1977} have leaves of fractional polynomial growth rates, showing that the properties  of the leaves with polynomial growth rates   is  again a subtle subject. 
Foliations with orbits having  polynomial growth rates have been investigated by  Egashira  \cite{Egashira1993a,Egashira1993b,Egashira1996} and Badura \cite{Badura2005}.

Next consider the relation between the growth rates of orbits, the existence of invariant measures for $\cGF$, and the F{\o}lner condition. The seminal paper in this area   was the 1957 work by Schwartzman \cite{Schwartzman1957}, but owes its modern development to Plante \cite{Plante1975b} and Ruelle \& Sullivan \cite{RS1975}.
 
For each  $x \in \BR$, one defines   a  $\cGF$-invariant, Borel probability measure $\mu_x$ on $\cT$.  Given  a continuous function $f \colon \cT \to \mR$, set 
$$
\mu_x(f) = \frac{1}{\# \cO(x)} \sum_{y \in \cO(x)} ~ f(y)
$$
 It is obvious that $\mu_x = \mu_y$  for all  $y \in \cO(x)$, and so $\mu_x$ is  $\cGF$-invariant.
 
For each  $x \in \SR$,  one can also define    a  $\cGF$-invariant, Borel probability measure $\mu^*_x$ on $\cT$,  but this requires a more subtle     averaging process.   Given  a continuous function $f \colon \cT \to \mR$, set 
\begin{equation} \label{eq-average}
\mu^{\ell}_x(f) = \frac{1}{\# B_{\cR}(x, \ell)} \sum_{y \in B_{\cR}(x, \ell)} ~ f(y)
\end{equation}
Then $\mu^{\ell}_x$ is a  Borel probability measure on $\cT$, but need not be $\cGF$-invariant. 
\begin{theorem}[Plante \cite{Plante1975b}, Ruelle-Sullivan \cite{RS1975}]\label{thm-PS}
Let $\ds \mu_x^* $ be a weak-* limit of the sequence $\{\mu^{\ell}_x \mid \ell = 1, 2, \ldots\}$ in the unit ball of the dual of $C(\cT)$. Then $\mu_x^*$ 
is a  $\cGF$-invariant, Borel probability measure on $\cT$. Moreover, the support of $\mu_x^*$  is contained in the orbit  closure 
$\overline{\cO(x)}$. 
\end{theorem}

Note that in general, the weak-* limit of a bounded sequence in the dual space $C(\cT)^*$ need not be unique. Thus, associated to each $x \in \cT$ there may be more than one limiting invariant measure $\ds \mu_x^* $. For example , one may have  a leaf $L$ of subexponential growth rate of a foliation $\F$, so that the ends of $L$ limit to distinct compact leaves of $\F$, and that the choice of the converging subsequence of the sequence of probability measures $\{\mu^{\ell}_x\}$ yields invariant measures supported on distinct compact leaves. Of course, more sophisticated examples are possible, but this simplest case shows the basic idea.

Generically, one expects that the measure $\mu_x^*$ is singular with respect to Lebesgue measure on $\cT$. By definition, for a Borel subset $E \subset \cT$, the measure $\mu_x^*(E)$ depends on the rates of accumulation of points in the balls $B_{\cR}(x, \ell)$ near the set $E$.   In other words, $\mu_x^*(E)$ depends on the asymptotic properties of the ends of the leaf $L_x$ near the set $E_{\F}$. The problem to determine the support of a transverse, $\cGF$-invariant measure   has been studied for codimension-one foliations of 3-manifolds 
\cite{CantwellConlon1977,CantwellConlon1978,Tsuchiya1979a,Tsuchiya1980b}.   

 The averaging formula  (\ref{eq-average}) is a special case of the F{\o}lner condition. 
Given a subset $S \subset \cO(x)$, the ``boundary'' of $S$ is defined by
\begin{equation}
\partial S \equiv \{ x \in S \mid d_{\cR}(x, \cO(x) \setminus S) =1\} ~ \cup ~  \{ y \in \cO(x) \setminus S) \mid d_{\cR}(y,  S) =1\}
\end{equation}
For example, note that $\partial B_{\cR}(x, \ell) \subset B_{\cR}(x, \ell + 1)$ and $B_{\cR}(x, \ell - 1) \cap \partial B_{\cR}(x, \ell)    = \emptyset$. Hence,  
\begin{equation}\label{eq-folnerballs}
\partial B_{\cR}(x, \ell) \subset B_{\cR}(x, \ell + 1) \setminus B_{\cR}(x, \ell - 1)
\end{equation}

\begin{definition} A sequence of \emph{finite} subsets $\{ S_{\ell} \subset \cO(x) \mid \ell = 1, 2, \ldots \}$ is said to be a \emph{F{\o}lner sequence} for the orbit $\cO(x)$ if for all $\ell \geq 1$, 
\begin{equation}\label{eq-folner}
S_{\ell} \subset S_{\ell +1 } ~, ~   \cO(x) =  \bigcup_{\ell = 1}^{\infty} ~ S_{\ell} ~ , ~  \lim_{\ell \to \infty} ~ \frac{\# \partial S_{\ell}}{\# S_{\ell}} = 0
\end{equation}
\end{definition}

 It is an easy exercise using (\ref{eq-folnerballs}) to show: 
 \begin{proposition} Let $x \in \SR$, then the sequence $\{ B_{\cR}(x, \ell) \mid \ell = 1, 2, \ldots \}$ is a F{\o}lner sequence for $\cO(x)$.
 \end{proposition}
Note that it is   possible for $x \in \FR$ and yet $\cO(x)$   admits a F{\o}lner sequence.  The standard example if for the leaves of  a foliation defined by the locally free action of a non-unimodular solvable Lie group. 
 
Now, let  $\mathfrak{S} = \{ S_{\ell} \mid \ell = 1, 2, \ldots \}$ be a F{\o}lner sequence for $\cO(x)$, and given    a continuous function $f \colon \cT \to \mR$, set 
\begin{equation} \label{eq-folav}
\mu^{\ell}_{\mathfrak{S}}(f) = \frac{1}{\# S_{\ell}} \sum_{y \in S_{\ell}} ~ f(y)
\end{equation}
Then Theorem~\ref{thm-PS} admits a generalization.

\begin{theorem}[Goodman-Plante \cite{GoodmanPlante1979}]\label{thm-GP}
Let  $\mathfrak{S} = \{ S_{\ell} \mid \ell = 1, 2, \ldots \}$ be a F{\o}lner sequence for $\cO(x)$. 
Let  $\ds \mu_{\mathfrak{S}}^* $ be a weak-* limit of the sequence $\{\mu^{\ell}_{\mathfrak{S}} \mid \ell = 1, 2, \ldots\}$ of Borel probability measures on $\cT$. Then $\ds \mu_{\mathfrak{S}}^* $
is a  $\cGF$-invariant, Borel probability measure on $\cT$. Moreover, the support of $\ds \mu_{\mathfrak{S}}^* $ is contained in the closure 
$\overline{\cO(x)}$. 
\end{theorem}

 Finally, we give a generalization of these ideas, with   applications to the study of   dynamics of minimal sets.
 
 \begin{proposition}[\cite{Hurder2008b}]\label{prop-growthbound}
 Suppose that $K \subset \cT$ is a closed saturated subset such that for all $\ell > 0$, there exists $x_{\ell} \in K$ for which 
 ${\rm gr}(\cR, x_{\ell}) < 1/\ell$. Then there exists a  $\cGF$-invariant, Borel probability measure $\mu^*$ on $\cT$ with support in $K$. 
 \end{proposition}
 \begin{corollary}\label{cor-growthbound}
 Suppose that $K \subset \cT$ is a closed saturated subset such that there is no $\cGF$-invariant, Borel probability measure $\mu^*$ supported on $K$. Then there exists $\lambda_K > 0$ so that for all $x \in K$,  ${\rm gr}(\cR, x) \geq \lambda_K$. 
 \hfill $\eop$
 \end{corollary}
 {\bf Proof of Proposition~\ref{prop-growthbound}:} We are given that for each $\ell \geq  1$
 $${\rm gr}(\cR, x_{\ell}) = \limsup_{i \to \infty} ~ \frac{\ln \{ {\rm Gr}(\cR, x_{\ell}, i) \}}{i} < 1/\ell$$
Hence, there exists $i_{\ell}$ such that  
 $$\frac{\# \{ B_{\cR}(x_{\ell}, i_{\ell} + 1) \setminus B_{\cR}(x_{\ell}, i_{\ell} - 1)\}}{\#  B_{\cR}(x_{\ell}, i_{\ell})} \leq  1/\ell$$
 Then by (\ref{eq-folnerballs}) we have that 
 $$\frac{\#  \partial B_{\cR}(x_{\ell}, i_{\ell})}{\#  B_{\cR}(x_{\ell}, i_{\ell})} \leq  1/\ell$$
 Given    a continuous function $f \colon \cT \to \mR$, set 
  $$\mu^{\ell}(f) = \frac{1}{\#    B_{\cR}(x_{\ell}, i_{\ell})} \sum_{y \in \#  B_{\cR}(x_{\ell}, i_{\ell})} ~ f(y)$$
Let $\mu^*$ be a weak-* limit of the sequence  $\{\mu^{\ell} \mid \ell = 1, 2, \ldots\}$.  Then $\mu^*$ is a  $\cGF$-invariant, Borel probability measure  with support in $K$. \hfill $\eop$

\section{Cocycles over metric equivalence relations} \label{sec-cocycles}

In this section, we consider the properties of cocycles over the groupoid $\GF$ and the equivalence relation $\cRF$. The importance of this data for the secondary classes is seen intuitively, when we consider the  derivative cocycle $D \colon \GF \to GL(q,\mR)$ introduced in \S3. In essence, this is just the holonomy transport data for the Bott connection $\nabla^{\F}$ restricted to leaves of $\F$. The cocycle records the linear structure of the groupoid $\GF$ along paths in leaves. Thus, for $x \in \cT$, the cocycle data 
$D \colon \GF^x \to GL(q,\mR)$ 
gives a family of linear approximations  to  the action of $\GF$ along the orbit $\cO(x)$. The values of the Weil measures of a set $E \in \cBF$ can be estimated by these approximations.

We begin with some general considerations. Let $G$ be a Polish group, equipped with a metric $d_G$; in particular, let $G$ be a closed subgroup of the matrix group $GL(q, \mR)$ equipped with a left invariant metric. A Borel (measurable) $G$-cocycle over $\GF$ is a Borel (measurable) map $\phi \colon \GF \to G$ which satisfies the cocycle equation,
\begin{equation}
\phi(\gamma_2 \cdot \gamma_1) = \phi(\gamma_2) \cdot \phi(\gamma_1) ~, ~ {\rm for~all}~ \gamma_1 \in \GF^{x,y} ~, ~\gamma_2 \in \GF^{y,z} 
\end{equation}
In formal terms, a cocycle is   a functor from the small category $\GF$ to the category with one object, whose group of morphisms is $G$. 

We   also consider cocycles  $\phi \colon \cRF \to G$ over the equivalence relation $\cRF$. Via the map $(s \times r) \colon \GF \to \cRF$ each cocycle over $\cRF$ defines a cocycle over $\GF$. In fact, the cocycles obtained in this way are exactly those $\phi \colon \GF \to G$ which are trivial  on the holonomy groups $\GF^{x,x}$ for all $x \in \cT$.

 Recall that two cocycles, $\phi, \psi \colon \GF \to G$ are \emph{cohomologous} if there exists a Borel (measurable) map $f \colon \cT \to G$ such that 
\begin{equation}\label{eq-coboundary}
\psi(\gamma) = f(y)^{-1} \cdot \phi(\gamma) \cdot f(x)  ~, ~ {\rm for~all}~ \gamma \in \GF^{x,y}
\end{equation}
The map $f$ is called the coboundary (or transfer function) between the cocycles $\phi$ and $\psi$.
The equivalence class of a $G$-cocycle $\phi$ is denoted by $[\phi]$, and the set of equivalences classes is denoted  by $H^1(\GF ; G)$. If $G$ is abelian, this is a group. Note that the cohomology set $H^1(\GF ; G)$ usually depends strongly on the notion of equivalence, whether we allow Borel, measurable or only smooth coboundaries.

One of the fundamental    points of cocycle theory in dynamical systems is that because the 
   cocycles   are   Borel functions (or possibly measurable functions where sets of measure zero can be neglected) and the  coboundaries are of the same type, means that the  cohomology set $H^1(\GF ; G)$  is a dynamical (or ergodic) property of the system.   The goal is to find normal forms for the cocycle, and then deduce    dynamical implications from it.

The use of cocycles to study problems in ergodic theory was pioneered by Mackey in two seminal papers \cite{Mackey1963,Mackey1966}, where they were called ``virtual representations'', and used  to construct unitary representations for Lie groups. Subsequent applications of cocycle theory to the study of  group actions on   manifolds  lies behind some of the deepest results in the field; two of the original references remain the best, the books by Schmidt \cite{Schmidt1977} and  Zimmer  \cite{Zimmer1984};   see also Margulis \cite{Margulis1991}.
 
For example, cocycle theory is a fundamental  aspect of the Oseledets Theorem  \cite{Oseledets1968}, which given    a $C^{1+\alpha}$-diffeomorphism $f \colon N^q \to N^q$ and an invariant probability measure $\mu$, yields a ``diagonal'' canonical form for the derivative cocycle $Df$ $a.e.$ with respect to $\mu$. The diagonal entries of this normal form are called the Lyapunov spectrum of $f$ over $\mu$.  Pesin  theory then uses this normal form to deduce dynamical properties of the map $f$ (see \cite{Oseledets1968,Pesin1977,KatokMendoza1995}).   Katok's celebrated paper \cite{Katok1980} provides the model application of this   technique. For further discussion, see for example   the recent  text  Barreira and Pesin \cite{BarreiraPesin2007}. 
The use of Oseledets Theory and Pesin Theory for the analysis of the transverse  dynamics of foliations was introduced by the author in the paper \cite{Hurder1988}. 

A key technical aspect in Oseledets Theory  is the growth rate of a cocycle $\phi$ defined  over metric equivalence relations. 

Assume that $G \subset GL(m, \mR)$ is a linear group. Fix a left-invariant metric $\rho$ on $G$, or possibly a pseudo-metric $\rho$ such that the set $G_0 = \{g \in G \mid \rho(Id, g) = 0\}$ is compact. For example, define the pseudo-norm,  for $A \in GL(m, \mR)$, 
\begin{equation}\label{eq-norm}
 | A | = \max \{\ln \|A\|, \ln \|A^{-1}\|\}
\end{equation}
where $\|A\|$ denotes the usual matrix norm, $\ds \|A\| = \sup_{0 \ne \vec{v} \in \mR^m} ~  \|A \vec{v}\| / \|\vec{v}\|$. 

Define $\rho(A,B) = |A^{-1} \cdot B|$, which is   a left-invariant pseudo-metric.  
Note that $\rho(Id, A) = 0$ precisely when $\|A\| = 1 = \|A^{-1}\|$, which implies that $A$ is   orthogonal.
\begin{definition}\label{def-tempered}
Let $c > 0$  and $E \in \cBF$.
The $G$-cocycle $\phi \colon \GF^E \to G$ is \emph{$c$-tempered} on $E$,  if   for all $x \in E$, 
\begin{equation}\label{eq-tempered}
|\phi(\gamma)| \leq c \cdot \|\gamma\|_x  ~, ~ \gamma \in \GF^x
\end{equation}
\end{definition}

\begin{definition}\label{def-exptype}
Let $b > 0$  and $E \in \cBF$. 
 The  $G$-cocycle $\phi \colon \GF^E \to G$ has \emph{exponential type $b$} on $E$ if for all $x \in E$, 
\begin{equation}\label{eq-exptype}
\limsup_{\ell \to \infty} ~ \frac{\max \{|\phi(\gamma)| ~ {\rm for} ~ \gamma \in \GF^x ~ , ~ \|\gamma\|_x  \leq \ell \}  }{\ell} \leq b
\end{equation}
and $b$ is the least such $b \geq 0$ such that (\ref{eq-exptype}) holds. 
If $b = 0$, then we say that $\phi$ has \emph{subexponential type} on $E$.
\end{definition}
The following is an immediate consequence of the definitions. Let $E \in \cBF$.
\begin{lemma}
 If $\phi$  is \emph{$c$-tempered} on $E$, then    $\phi$   has exponential type $b \leq c$ on $E$.
\end{lemma}

We now turn the discussion to the properties of the derivative cocycle. Recall   there is given an identification $T_x\cT \cong \mR^q$. Given $\gamma \in \cT_{\F}^{x,y}$ choose $g \in \cGF$ such that $\gamma = [g]_x$. Set  $D\gamma = D_x(g) \colon T_x\cT \to T_y\cT$, which yields  $D\gamma \in GL(q, \mR^q)$. By the Chain Rule, the map $D \colon \GF \to GL(q, \mR^q)$ is a cocycle.

Let ${\rm Hol}_0(\F) \subset \cT$ denote the set of points without holonomy. That is, $x \in  {\rm Hol}_0(\F)$  if the set $\GF^{x,x}$ contains only the germ of the identity.  Epstein, Millett and  Tischler proved in \cite{EMT1977} that  ${\rm Hol}_0(\F)$ is a dense  $G_{\delta}$ in $\cT$. In particular, ${\rm Hol}_0(\F) \in \cBF$.
 
Let ${\rm Hol}_1(\F) \subset \cT$ denote the set of points with trivial linear holonomy. That is, $x \in  {\rm Hol}_1(\F)$  if for each  $\gamma \in \GF^{x,x}$ the derivative map $D\gamma$ is the identity.  
\begin{proposition}[Hurder-Katok, Proposition~7.1, \cite{HK1987}] \label{prop-linholo}
${\rm Hol}_1(\F)  \in \cBF$ is a set of full measure. That is, the set of points $x \in \cT$ for which there exists $\gamma \in \GF^{x,x}$ with $D\gamma \ne Id$   is a Borel set of Lebesgue measure zero.
\end{proposition}
\sop
Let $g \in \cGF$. Define ${\rm Fix}(g) = \{x \in D(g) \mid g(x) = x\} \subset \cT$.    
Let  $x \in {\rm Fix}(g)$ be  a point of Lebesgue density $1$. 
Then for any vector $\vec{v} \in T_x\mR^q$ there exists a sequence $\{y_{\ell} \mid y_{\ell} \ne x, y_{\ell} \in {\rm Fix}(g)\}$ such that $y_{\ell} \to x$ and limits from the direction of $\vec{v}$. 
As $g(y_{\ell}) = y_{\ell}$, this implies that $D_x(g)(\vec{v})  = \vec{v}$. Hence, the subset  
$\{x \in {\rm Fix}(g) \mid D_x(g) \ne Id\}$ has Lebesgue measure zero. As $\cGF$ is countably generated by the compositions of the elements of $\cGF^{(1)}$, this implies that the set of points in $\cT$ for which there is some $\gamma \in \GF^{x,x}$ with $D\gamma \ne Id$ has Lebesgue measure $0$. \hfill $\eop$

\begin{corollary}\label{cor-strict}
There exists  $Z \in \cBF$ with full measure such that for all $x \in Z$, the map $D \colon \GF^x \to GL(q,\mR)$ is determined  by the endpoints $(x,y) \in \cRF$.
\end{corollary}

Define a new cocycle $\whD \colon  \GF^x \to GL(q,\mR)$ by setting $\whD\gamma = D\gamma$ if $\gamma \in \GF^x$ and $x \in {\rm Hol}_1(\F)$, and $\whD\gamma = Id$ otherwise. Then $\whD$ is a cocycle which depends only on the value $(s \times r)(\gamma) \in \cRF$, hence can be considered as a cocycle over the equivalence relation. Note that $D$ and $\whD$ differ at most on a set of Lebesgue  measure $0$.

The point of this   modification procedure, is that there it is often possible to prove the derivative cocycle $D$ has exponential type $b$ on the subset $Z$ of Corollary~\ref{cor-strict}, or some other structure theorem holds for $D$ on  $Z$. Via the above modification, one can then assume this structure holds for all of $\cT$. 

The Radon-Nikodym  cocycle    is defined by  $\nu(\gamma) = \ln \{ \det (D\gamma)\}$. 
That is,  $\nu \colon \GF \to \mR$ is the additive cocycle obtained from the volume expansion of the transverse linear holonomy maps, and as such,  is an important aspect  of the dynamics of the foliation. A basic observation in \cite{Hurder1986} is that the properties of this cocycle connect together the Godbillon measure on $\cBF$ and the growth rates of the leaves of $\F$. The key technical result is as follows:

\begin{theorem}[Theorem~4.3, Hurder \cite{Hurder1986}] \label{thm-subexp}
There exists $E_{\nu} \in \cBF$ with  $E_{\nu} \subset \BR \cup \SR$ of full relative Lebesgue measure, such that 
$\nu(\gamma) = \ln \{ \det (D\gamma)\}$ has exponential type $0$ on $E_{\nu}$. In particular, if the set  $\FF$ of leaves with exponential growth type has Lebesgue measure zero, then $\nu(\gamma) $ has exponential type $0$ almost everywhere on $\cT$.
\end{theorem}

A version of this result can be found in the monograph \cite{Schmidt1977} for the case of a $\mZ$ action on a compact manifold. The proof for pseudogroup actions involves   greater technical complications, although the principle of the proof remains the same: if the orbits of a pseudogroup on a set $E \in \cBF$ with positive measure have subexponential growth rate, then there is ``no room at infinity'' for there to be a subset of $E$ with positive Lebesgue measure having exponential growth type $b > 0$.

\section{Amenable foliations} \label{sec-amenable}

The concept of amenable equivalence relations in the measurable setting was introduced by Zimmer \cite{Zimmer1978}, as an analog for the concept of a hyperfinite equivalence relation \cite{Dye1959,Krieger1976,Series1980a}. 
A celebrated result of Connes, Feldman and Weiss  \cite{CFW1981} proved that   an amenable equivalence relation is generated by a single transformation.     The properties of amenable groupoids  in the topological setting is a related but less well understood theory \cite{ADR2000}.
In this section, we recall results of Zimmer about the applications of amenability to the structure theory of measurable cocycles, and the applications of this theory to foliations.

Let $X$ be a compact convex space, and $G = {\rm Aut}(X)$ be the group of affine automorphisms of $X$.
Given $E \in \cBF$ and a $\cGF$-quasi-invariant  measure $\mu$ for which $\mu(E) = 1$,   the equivalence relation $\cRF^E$ is \emph{$\mu$-amenable} if, given any such $X$ and cocycle $\phi \colon \cRF^E \to {\rm Aut}(X)$, there exists a $\mu$-measurable section $h \colon E \to X$ such that for all $(x,y) \in \cRF^E$, 
$$\phi(x,y) \cdot h(x) = h(y)$$
That is, for the fiberwise action of $\cRF^E$ on the sections of the bundle $E \times X \to E$ defined by the cocycle $\phi$, there exists a ``global fixed-point'' -- the section $h$.

For a finitely generated group $\G$, the existence of a  F{\o}lner sequence  for   $\G$, equipped with the word metric, is equivalent to the group $\G$ being amenable \cite{Greenleaf1969,Pier1984}. Every continuous action of an amenable group on a compact metric space, 
$\varphi \colon \G \times Y \to Y$, generates an equivalence relation $\cR_{\varphi}$ which is $\mu$-amenable for any $\G$-quasi-invariant Borel probability measure on $Y$ (see \cite{ADR2000}). Moreover, a F{\o}lner sequence for $\G$ generates an orbit-wise ``uniform F{\o}lner sequence'' for $\cR_{\varphi}$.

Kaimanovich   observed that for a general equivalence relations,  amenability is no longer characterized 
by the existence of orbit-wise F{\o}lner sequences \cite{Kaimanovich2001}, although the two concepts are still closely related. (See also   \cite{AdamsLyons1991,Brooks1983,CarriereGhys1985,CFW1981}, and the recent work of Rechtman \cite{Rechtman2008}.)  
For the converse,   one needs a strong form of ``uniform F{\o}lner sequences'' on the orbits   $\cO(x)$ for $x \in E$ to obtain that $\cRF^E$ is amenable. For example,   if 
$E \subset  \BR \cup \SR$ then each orbit admits such  a uniform     F{\o}lner sequence.
\begin{proposition}[\cite{HK1987,Kaimanovich2001,Series1980a}] \label{prop-folner}
 Let $E \in \cBF$ with $E \subset  \BR \cup \SR$. Then $\cRF^E$ is amenable with respect to any $\cGF$-quasi-invariant measure $\mu$.  
\end{proposition}

The Roussarie foliation is amenable \cite{Bowen1977}, but every leaf has exponential growth.
More generally, suppose that $\F$ is defined by a locally free action of an amenable Lie group $H$ on $M$, then $\cRF$ is amenable.    If the parabolic subgroup $H$ has exponential growth type, then every leaf of $\F$ will have exponential growth type.

Suppose that $\F_{\varphi}$ is obtained via the suspension of a $C^1$-action $\varphi\colon \G \times N \to N$, where $\G$  is 
 a  finitely generated group and $N$ is a closed $q$-dimensional manifold. As remarked above, if  $\G$ is an amenable group, then  $\cR_{F_ {\varphi} }$ is always an amenable equivalence relation. For $\G$ a word-hyperbolic group, hence non-amenable, Adams proved in \cite{Adams1994} that the group action defines  an amenable equivalence relation if  $N$ is homeomorphic to the   boundary of  $\G$. The Roussarie example is also of this type. 
 
Let $H \subset GL(q,\mR)$ be an amenable subgroup. Then the linear action of $H$ on $\mR^q$ preserves a flag 
$$\{0\} = V_0 \subset V_1 \subset \cdots \subset V_k = \mR^q$$
such that the induced action on each quotient $V_i/V_{i-1}$ for $1 \leq i \leq k$ preserves a positive-definite inner product, up to similarity. 
It follows that $H$ is conjugate to a subgroup of a maximal parabolic subgroup of $GL(q, \mR)$, 
which is the maximal subgroup of $GL(q,\mR)$ preserving a ``standard flag'',  
$$\{0\}  \subset \mR^{m_1} \subset \mR^{m_2} \subset \cdots \subset \mR^{m_k} = \mR^q$$ where $m_i = \dim V_i$, and which acts via $\mR \times O(m_i - m_{i-1})$ on the quotient $\mR^{m_i}/\mR^{m_{i-1}}$. There are $2^q$ conjugacy classes  of maximal parabolic Lie  subgroups of $GL(q,\mR)$, which we label $H_i \subset GL(q, \mR)$ for $i \in \cI_q  = \{1, 2, \ldots, 2^q\}$.

One of the basic applications of amenability to cocycle theory is given by the following result of Zimmer:
\begin{theorem}[Zimmer \cite{Zimmer1982a,Zimmer1984}]\label{thm-zimmer}
Let $(X, \cR)$ be an ergodic amenable discrete equivalence relation, with respect to a quasi-invariant measure $\mu$. 
Let  $G$ be a real algebraic group. Then for every cocycle $\phi \colon \cR \to G$, there is an amenable subgroup $H_{\phi} \subset G$ and a cocycle $\psi \colon \cR \to H_{\phi}$ which is  $\mu$-measurably cohomologous to $\phi$ in $G$.
\end{theorem}

Given $E \in \cBF$, the restricted equivalence relation $ \cRF^E$ admits an ergodic decomposition $\ds E = \cup_{\alpha \in \cA} X_{\alpha}$. Theorem~\ref{thm-zimmer} can then be applied to each $X_{\alpha}$ to conclude that the restriction of   the derivative cocycle $\whD \colon \cRF^E \to GL(q,\mR)$ to   $X_{\alpha}$ is cohomologous to a cocycle with values in one of the canonical amenable subgroups $H_i \subset GL(q, \mR)$. This defines the measurable function   $\sigma \colon E \to \cI_q$, where $x \in X_{\alpha}$ implies that $\whD | X_{\alpha}$ is conjugate to a cocycle with values in $H_{\sigma(x)}$. These are the ingredients of the proof of the following result from \cite{HK1987}:

 \begin{theorem}[Corollary~3.3,  \cite{HK1987}] \label{thm-amenred}
Let $E \in \cBF$ and assume that $\cRF^E$ is an amenable equivalence relation with respect to Lebesgue measure.
Let $\whD \colon  \cRF^E   \to GL(q, \mR)$  be the restriction of the modified derivative cocycle.
Then there exists 
\begin{itemize}
\item a measurable selection function $\sigma \colon E \to \cI_q$
\item a cocycle $\psi \colon \cRF^E \to GL(q, \mR)$ such that  for $(x,y) \in \cRF^E$, $\psi(x,y) \in H_{\sigma(x)}$
\item a   cohomology $f \colon E \to GL(q, \mR)$ between  $\psi$ and $\whD \colon \cRF^E \to GL(q, \mR)$.
\end{itemize}
Moreover, it can be assumed that  $\psi$ is $c$-tempered for some $c  > 0$.
\end{theorem}

The consequences of Theorem~\ref{thm-amenred}   for the dynamical properties of $\F$ in $E$ are not well-understood. However, it is known to have  strong consequences for the values of the Weil measures of $E$. This follows from a basic result about the relation between Lie algebra cohomology, continuous cohomology, and  $H^*(\mathfrak{gl}(q, \mR), O(q))$:  

 \begin{proposition}[Proposition~3.8, \cite{Hurder1985b}] \label{prop-leafvan}
 Let 
 $$\chi \colon H^{\ell}(\mathfrak{gl}(q, \mR), O(q)) \cong \Lambda(h_1, h_3, \ldots, h_{q'}) \to H^{\ell}_{cont}(H)$$
 be the characteristic homomorphism induced from the inclusion $H \subset GL(q, \mR^q)$, where $H^*_{cont}(H)$ denotes the continuous cohomology of $H$. If $H$  is amenable, the $\chi$ is the zero map for $\ell > 1$.
 \end{proposition}
 From this we conclude: 
 \begin{theorem}[Theorem~3.5, \cite{HK1987}] \label{thm-amenvan}
 Let $E \in \cBF$ be such that $\cRF^E$ is an amenable equivalence relation with respect to Lebesgue measure. Then for all classes 
 $h_I \in H^{\ell}(\mathfrak{gl}(q, \mR), O(q))$ with $\ell > 1$, the Weil measure $\chi_E(h_I) = 0$. 
 \end{theorem}
 
This suggest another fundamental decomposition of the space $\cT$, defined $a.e.$ with respect to Lebesgue measure:
\begin{equation}\label{eq-amendecomp}
\cT = \AR \cup \KR
\end{equation}
where $\AR$ is the maximal subset $E \in \cBF$ such that $\cRF^E$ is amenable, and $\KR$ is the complement. In terms of the von~Neumann algebra $\cM^*(\cRF)$ associated to $\cRF$ the decomposition (\ref{eq-amendecomp}) corresponds to a decomposition of $\cM^*(\cRF)$ into its injective and non-injective summands \cite{Connes1994}.

It would be very interesting to understand better the foliation dynamics in both invariant subsets of  (\ref{eq-amendecomp}). For example, 
the main result of \cite{Hurder1985b} states:
\begin{theorem}[Theorem~1, \cite{Hurder1985b}]
Let $\F$ be a $C^r$-foliation, for $r \geq 2$. Suppose there exists $x \in \cT$ such that the image $D \colon \GF^{x,x} \to GL(q,\mR)$ is not an amenable group. Then $L_x$ is in the closure of the set $\FF$ of leaves with fast growth.
\end{theorem}

Does a similar conclusion hold for subsets $E \in \cBF$ such that the range of $D \colon \cRF^E \to GL(q,\mR)$ is not amenable?
This is related to the question whether Theorem~\ref{thm-zimmer} is true for Borel sets, and not just as a measurable decomposition.

\section{Transverse infinitesimal expansion} \label{sec-expansion}

We  introduce a basic invariant of the derivative cocycle, which measures the degree of asymptotic infinitesimal expansion of $D \colon \GF \to GL(q, \mR)$ along an orbit. 
\begin{definition}  
The   transverse expansion rate  function for $\GF$ at $x$ is
 \begin{equation}\label{ep-atgn}
 e(\GF , \ell, x) =  
   \max_{\stackrel{\gamma \in \GF^x}{ \|\gamma \|_x \leq \ell}}~ \left\{  \frac{ \ln \left(  \max \{ \| D\gamma  \| , \| (D \gamma)^{-1} \|\} \right) }{\ell }\right\}
\end{equation}
  \end{definition}
 Note that $ e(\GF , \ell, x)$ is a Borel function on $\cT$, as each norm function $\| D \gamma  \|$ is continuous for $x \in D(g)$ where $g \in \cGF$ with $\gamma = [g]_x$ and the maximum of Borel functions is Borel.
  
\begin{definition}  
The asymptotic transverse expansion rate at $x$ is
 \begin{equation}\label{eq-atg}
e(x)  = e(\GF, x) = \limsup_{\ell \to \infty} ~  e(\G_{F}, \ell, x) ~ \geq ~ 0
  \end{equation}
\end{definition}  
The limit of Borel functions is Borel, and each $ e(\G_{\F}, n, x)$ is Borel, hence $e(\GF, x)$ is Borel. 
Note that (\ref{eq-atg})  is just a pointwise version of the estimate (\ref{eq-exptype}) which appears in Definition~\ref{def-exptype}. The value $e(x)$ can be thought of as the maximal Lyapunov exponent of the transverse holonomy at $x$.

\begin{lemma}\label{lem-orbitinv}
Suppose that $(x,y) \in \cRF$, then $e(x) = e(y)$.
\end{lemma}
\sop Let $\gamma \in \GF^{z,y}$ and $\gamma' \in \GF^{x,z}$ with $\|\gamma' \|_x = 1$. Note that  
$$ \|\gamma\|_z  -1 \leq \| \gamma' \cdot \gamma \|_x \leq \|\gamma\|_z + 1$$
and use the estimates, for $A, B \in GL(q,\mR)$,  $\|A\| \cdot \|B\|^{-1} \leq \|B \cdot A\| \leq \|A \| \cdot \|B\|$, 
where $A = (D_z\gamma)^{\pm 1}$ and $B = (D  \gamma')^{\pm 1}$. \hfill $\eop$

\medskip

The following is the precise statement of Theorem~\ref{thm-main1} of the Introduction.
\begin{theorem} [Hurder \cite{Hurder2008b}]\label{thm-decomp} 
  Let $\F$ be a $C^1$-foliation on a compact manifold $M$. Then $M$ has a decomposition into disjoint subsets of $\cBF$, 
$M = \EF \cup \PF \cup \HF$, which are the saturations of the sets
defined by:
\begin{enumerate}
\item Elliptic:  $\ER   = \{ x \in \cT \mid \forall ~ n \geq 0, ~  e(\GF , m, x)  \leq \kappa(x) \}$
\item Parabolic:  $\PR     = \{ x \in \cT \setminus  \ER  \mid  e(\GF, x) = 0 \}$
\item Hyperbolic:  $\HR     = \{ x \in \cT \mid  e(\GF, x) > 0 \}$
\end{enumerate}
\end{theorem}
Note that $x \in \ER$ means that  the holonomy homomorphism $D \colon \GF^x \to GL(q, \mR)$ has bounded image. The constant $\kappa(x) = \sup \{ \|D\gamma \|    ~{\rm for} ~  \gamma \in \GF^x \}$.

Note that the dynamics of a matrix $A \in GL(q,\mR)$ acting on $\mR^q$ is divided into three types: elliptic (or isometric);  parabolic (or distal); and partially-hyperbolic, when there is a non-unitary eigenvalue.  
The nomenclature in Theorem~\ref{thm-decomp} reflects this: The elliptic points are the points of $\cT$  where the infinitesimal holonomy transport ``preserves ellipses up to bounded distortion'';  that is, it is measurably isometric. The parabolic points are where the    infinitesimal holonomy acts similarly to that of a parabolic subgroup of $GL(q,\mR)$; for example, the action is infinitesimally distal. The   hyperbolic points are where the  the infinitesimal holonomy has some degree of exponential expansion. Perhaps more properly, the set $\HR$   should be called ``non-uniform, partially-hyperbolic'' points;  but hyperbolic is suggestive enough.

\section{Foliation  entropy} \label{sec-entropy}
 
 The geometric entropy of the foliation pseudogroup, $h(\cGF)$,   was introduced by Ghys, Langevin and Walczak \cite{GLW1988} in 1988.  Their definition is a generalization of    the Bowen definition of topological entropy  for a map \cite{Bowen1971}. The study of the  properties of $h(\cGF)$ has been one of the major motivating concepts in foliation dynamics, much as  topological entropy has been a central focus for the study of the dynamics of a diffeomorphism   $f \colon N \to N$ for the past 50 years. We will also introduce a new concept, called the \emph{local entropy} of the pseudogroup $\cGF$ and relate it to both $h(\cGF)$ and to the dynamical invariants introduced previously.

Let $X \subset \cT$. We say that   $S = \{x_1, \ldots , x_{\ell}\} \subset X$ is $(k, \epsilon)$-separated for $\cGF$ and $X$ if 
$$\forall ~ x_i \not= x_j ~, ~ \exists ~ g \in \cG | X ~ \mbox{such that} ~ d_{\cT}(g(x_i), g(x_j)) \geq \epsilon, ~ {\rm where}~  \| g \| \leq k $$
Here, $\|g\| \leq k$ means that $g$ can be written as a composition of at most $k$ elements of the generating set $\cGF^{(1)}$ of $\cGF$.
Then set
\begin{equation}\label{eq-separated}
h(\cGF, X,k,\epsilon) = \max ~ \# \{ S \mid S \subset X \text{ is } (k,\epsilon)-\text{separated} \}
\end{equation}
When $X = \cT$,     set $h(\cGF,k,\epsilon)  = h(\cGF, \cT,k,\epsilon)$.  

\begin{definition}[Ghys, Langevin and Walczak \cite{GLW1988}] Let $\cGF$ be a $C^r$-pseudogroup, for $r \geq 1$. 
The \emph{geometric entropy}  of $\cGF$ on $X \subset \cT$ is 
 \begin{equation}\label{eq-entropy}
h(\cG, X) =   \lim_{\epsilon \to 0} \Big\{  \limsup_{k \to \infty} \frac{\ln \{h(\cGF,X,k,\epsilon)\}}{k}\Big\}
\end{equation}
The \emph{geometric entropy}  of $\F$ is   defined to be $h(\F) \equiv  h(\cGF, \cT)$. 
\end{definition}
 One of the fundamental points about the definition of $h(\cGF)$ is that the definition  of $\epsilon$-separated sets in (\ref{eq-separated}) is based on the groupoid  distance function, and the normalizing denominator is the groupoid distance, rather than  the number of points in a ball of radius $k$ as was used in other approaches to defining topological entropy for group actions.  With this modification of the usual definition, $h(\cGF) < \infty$ when $\cGF$ is a $C^1$-pseudogroup. Moreover,     $h(\cGF) $ is positive for many examples;  it reflects expansive or chaotic behavior in the dynamics of $\cGF$.

For a given foliation $\F$,  the value of $h(\cGF)$ for the associated pseudogroup $\cGF$ depends strongly on the choice of the covering of $M$ by foliation charts. Ghys, Langevin and Walczak show in \cite{GLW1988} that the property that $h(\cGF) > 0$ is independent of the choice of covering, so we may speak of a foliation $\F$ with positive geometric entropy. 

One of the main results in \cite{GLW1988}, Theorem~6.1,  is   a characterization  of   the $C^2$-foliations of codimension-one with $h(\F) > 0$ -- such foliations must have a resilient leaf. Recall that $x \in \cT$ is resilient if   there exists $g \in \cGF$ with $x \in D(g)$,  $g(x) = x$, and $g$ is a local one-sided contraction at $x$, such that the orbit $\cO(x)$ intersects the domain of the contraction in some point other than $x$.  The proof of Theorem~6.1 in \cite{GLW1988} relies on subtle properties of the theory of levels for $C^2$-foliations; this is discussed further in the text by Candel and Conlon \cite{CanCon2000}.

Ghys, Langevin and Walczak show in \cite{GLW1988} another remarkable result:
\begin{theorem}[Theorem~5.1, \cite{GLW1988}] Let $K \in \cBF$ be a closed subset, with 
 $h(\cGF,K) =0$. Then there exists a $\cGF$-invariant Borel probability measure $\mu$ with support in $K$.
\end{theorem}
The proof is similar to that of Proposition~\ref{prop-growthbound}, but  more delicate.

One of the basic  results about the topological entropy $h(f)$ of a $C^2$-diffeomorphism $f \colon N \to N$, is  the entropy formula of Margulis \cite{Margulis2004}, Pesin \cite{Pesin1977} and Ma{\~n}{\'e} \cite{Mane1987}. These formulas estimate $h(f)$ in terms of the Lyapunov spectrum of $f$ with respect to a ``sufficiently regular'' ergodic  invariant measure $\mu$ for $f$.  No analogs of these results have been established for the geometric entropy $h(\cGF)$; the author's paper \cite{Hurder1988} sketched an approach to proving such formulas, but this aspect of the study of foliation dynamics via cocycles remains to be proven. Proving this formula appears  difficult, based on the currently available techniques. 

 The  \emph{local geometric entropy} of $\cGF$ is  a variant of $h(\cGF,K)$.  Brin \& Katok introduced in  \cite{BrinKatok1983} the concept    of \emph{local measure-theoretic entropy} for maps. The concept of local entropy, as adapted to geometric entropy,  is very useful for the study of foliation dynamics. For example, it is used to establish the    relationship     between $h(\cGF)$ and the transverse Lyapunov spectrum of ergodic  invariant measures for the leafwise geodesic flow (see Theorem~\ref{thm-lyapspectrum} below.)

Recall that in the definition (\ref{eq-separated}) of $(k,\epsilon)$-separated sets, the separated points can be restricted to a given subset $X \subset \cT$, where the set $X$ need not be assumed  saturated.  If we take $X =  B(x,\delta) \subset \cT$, the open $\delta$--ball about $x \in \cT$, then we obtain a measure of the amount of ``expansion'' by the pseudogroup in an open neighborhood of $x$. Perform the same double limit process as used to define $h(\cGF)$ for the sets $B(x,\delta)$, but then also let the radius of the balls tend to zero, to obtain:
\begin{definition}[Hurder, \cite{Hurder2008b}] 
The \emph{local geometric entropy} of $\cGF$ at $x$ is
\begin{equation}\label{eq-localentropy}
h_{loc}(\cGF, x) = \lim_{\delta \to 0} \Big\{  \lim_{\epsilon \to 0} \Big\{  \limsup_{n \to \infty} \frac{\ln \{h(\cGF, B(x,\delta),k,\epsilon)\}}{k}\Big\} \Big\}
\end{equation}
\end{definition}

The local entropy has some very useful properties, which are elementary to show. For example, we have:
\begin{proposition}[\cite{Hurder2008b}] \label{prop-locentprop}
Let $\cGF$   a $C^1$-pseudogroup. Then $h_{loc}(\cGF,x)$ is a Borel function of $x \in \cT$, and $h_{loc}(\cGF,x) = h_{loc}(\cGF,y)$ if $(x,y) \in \cRF$. Moreover,
\begin{equation}\label{eq-locentest}
 h(\cGF)  ~  = ~  \sup_{x \in \cT} ~ h_{loc}(\cGF,x)
\end{equation}
\end{proposition}
It follows that there is a disjoint Borel decomposition into $\cGF$-saturated subsets
\begin{equation}\label{eq-locentdecomp}
\cT = \bZ_{\cR} \cap \bC_{\cR}
\end{equation}
where $\bC_{\cR} = \{x \in \cT \mid h(\cGF, x) > 0\}$ consists of the ``chaotic'' points for the action,  and $\bZ_{\cR} = \{x \in \cT \mid h(\cGF, x) = 0\}$ are the tame points.  
Here is a corollary of Proposition~\ref{prop-locentprop}:
\begin{corollary}  $h(\cGF) > 0$   if and only if $\bC_{\cR} \ne \emptyset$. 
\end{corollary}

 Here are two  typical results from \cite{Hurder2008b} about local entropy:

 \begin{theorem}
  Let $K \in \cBF$ be a $\cGF$-minimal set such that $h_{loc}(\cGF,x) > 0$ for some $x \in K$. Then   $ K \cap \HR \ne \emptyset$. 
\end{theorem}

\begin{theorem}
 Let $\cGF$   a $C^1$-pseudogroup. Then 
 \begin{equation}\label{eq-locentest2}
 h(\cGF)  ~  = ~  \sup_{x \in \NWF} ~ h_{loc}(\cGF,x)
\end{equation}
In particular,  $h(\cGF) = h(\cGF , \NWF)$.  
\end{theorem}

 One of the basic problems is to characterize the set $ \bC_{\cR}$ of chaotic points. For example,  
 suppose $E \in \cBF$ and $E \subset \HR$. If $E$ has    positive Lebesgue measure,  must $E \cap \bC_{\cR}$ be a large set? It seems likely that the closure $\overline{E}$ must contain chaotic points, with $ h_{loc}(\cGF,x) > 0$.  Must such points be dense in  $\overline{E}$?

\section{Tempering cocycles} \label{sec-tempering}

The results of the last sections discussed how dynamical properties of a foliation  $\F$  translate into properties of the normal derivative cocycle, $D \colon \GF \to GL(q, \mR)$. Some hypotheses on the dynamics, such as amenability or wandering on a subset $E \in \cBF$,  imply that the cocycle $\whD \colon \cRF^E \to GL(q, \mR)$ is cohomologous to a tempered cocycle with values in a subgroup $H \subset GL(q,\mR)$, where $H$ is amenable, or even the trivial subgroup. Other hypotheses  on the dynamics, such as subexponential growth or zero transverse expansion rate  on a subset $E \in \cBF$,  imply that the cocycle $\whD \colon \cRF^E \to GL(q, \mR)$ has estimates on its asymptotic growth.  

Recall that the cocycle $D \colon  \cRF^E \to GL(q, \mR)$ is the leafwise data contained in the Bott connection for the normal bundle of $\F$ restricted to $E$, and so   has applications to calculating the secondary invariants of $\F$ derived from the Bott connection. The first step in making these applications was the observation in the work of Heitsch and the author \cite{HH1984}, that the proof of   Theorem~\ref{thm-HH} shows that the Weil measures depend only on the \emph{measurable} cohomology class of the normal derivative cocycle.  This means that in calculating the  transgression  factor in the residual secondary classes of $\F$, we are allowed to use a measurable (but leafwise smooth!) normal framing of $Q$ that is adapted to the dynamics. For example, when $\F$ is amenable on $E$, we can change the normal framing of $Q|E$ so that $\whD | E$ takes values in an amenable subgroup of $GL(q,\mR)$ and hence the higher Weil measures must vanish. This is the basis for the proof of Theorem~\ref{thm-amenvan}.

The second technical observation is that if for all $\epsilon > 0$, a measurable framing of $Q|E$ can be chosen so that $\whD | E$  is equivalent to a cocycle with uniformly small norm less that $\epsilon$, then the Weil measures must be zero on $E$. This was shown for the Radon-Nikodym cocycle in \cite{HH1984}, and for the full normal cocycle in \cite{HK1987}. The problem then becomes to use asymptotic information about the norms of a cocycle to obtain a cohomologous cocycle which satisfies  uniform estimates.  This is done via \emph{tempering procedures}.

As mentioned in  section~\ref{sec-cocycles}, the  Pesin Theory  of a diffeomorphism $f \colon N^q \to N^q$ uses tempering to convert the asymptotic information obtained via the Oseledets Theorem applied to the derivative cocycle of $Df \colon N \times \mZ \to GL(q,\mR)$, into uniform estimates on the derivative of $f$. The work of Katok and the author  \cite{HK1987} introduced analogous   tempering procedures for measured equivalence relations. 

\begin{definition}\label{def-welltemp} 
A cocycle $\phi \colon \GF^E \to GL(q,\mR)$ is \emph{well-tempered} on $E \in \cBF$ if   for all $\epsilon > 0$,  there is a measurable coboundary 
$f_{\epsilon} \colon E \to GL(q,\mR)$ such that the cohomologous cocycle  $\psi_{\epsilon}= f_{\epsilon}^{-1} \cdot \phi \cdot f_{\epsilon}$ satisfies 
\begin{equation}\label{eq-wt}
|\psi_{\epsilon}(\gamma)| \leq \epsilon \cdot \|\gamma\|_x ~, ~ {\rm for ~all} ~ x \in E  ~, ~ \gamma \in \GF^x 
\end{equation}
\end{definition}

The tempering procedure for cocycles introduced in \cite{HK1987} applied to the orbits of subexponential growth  yields:
\begin{theorem}\label{thm-wt1}
Let $E  \subset \BR \cup \SR$. Suppose that  $\phi \colon \GF^E \to GL(q,\mR)$ has exponential type $0$. Then $\phi$ is well-tempered. 
\end{theorem}
This result combined with Theorem~\ref{thm-subexp} yields 
\begin{corollary}\label{cor-subexp}
Let $E_{\nu} \subset \BR \cup \SR$ be the set of full relative measure as in Theorem~\ref{thm-subexp}. 
Then the Radon-Nikodym cocycle $\nu \colon \cRF^{E_{\nu}} \to \mR$ is well-tempered. 
\end{corollary}

The tempering procedure in \cite{HK1987} (see also Stuck \cite{Stuck1991b}) can at best yield $c$-tempered cocycles where $c$ is greater than the growth rates of the leaves in $E$, which is why there is the restriction $E  \subset \BR \cup \SR$ in Theorem~\ref{thm-wt1}. 
A new tempering procedure was introduced in the work of Langevin and the author  \cite{HLa2000} which removed this restriction for cocycles with values in $\mR$. One can show using a combination of the tempering methods of  Hurder \& Katok \cite{HK1987} and Hurder \& Langevin  \cite{HLa2000} the following general result:
  
  \begin{theorem}\label{thm-wt2}\cite{Hurder2008b}
Suppose that  $\phi \colon \GF^E \to GL(q,\mR)$ has exponential type $0$, for $E \in \cBF$. Then $\phi$ is well-tempered. 
\end{theorem}

\section{Secondary classes and dynamics} \label{sec-secclasses2}

 We have now exhibited six dynamically-defined decompositions of $M$:
 
\begin{table}[htdp]
\begin{center}
\begin{tabular}{cclll}
$M$ & = & $M_I \cup M_{II} \cup M_{III}$ & --   M-vN type & (Equation~\ref{eq-MvN})\\
& = & $\WF \cup \NWF$ &  --  wandering & (Section~\ref{sec-top})\\
  & = & $\BF \cup \SF \cup \FF$ &  -- growth & (Theorem~\ref{prop-growth}.4) \\
  & = & $\AF \cup \KF  $ &  -- amenable & (Equation~\ref{eq-amendecomp}) \\
  & = & $\EF \cup \PF \cup \HF$ &  -- expansion & (Theorem~\ref{thm-decomp})\\
  & = & $\ZF \cup \CF$ &  -- local entropy & (Equation~\ref{eq-locentdecomp})\\
\end{tabular}
\end{center}
\end{table}
Of course, these are not independent, and there are multiple relations between the sets in these decompositions, some of which have been discussed previously. 

  In this   section, we discuss the localizations of the residual secondary classes to the   sets appearing  in these decomposition schemes.  This gives a summary of the results to date regarding how the secondary classes of a $C^2$-foliation are ``determined'' by its dynamics.  
 One of the main points of this paper is that these decomposition schemes   provide a framework for classifying foliations by their dynamics.

Note that all of the following results are true for arbitrary codimension $q \geq 1$, and the dynamical aspects of the conclusions generally hold for $C^1$-groupoids; or for $C^r$-groupoids with $r > 1$ if the particular  proof relies upon the existence of ``stable manifolds'' via Pesin Theory. 

\begin{proposition} [\cite{HH1984,Hurder1986}] \label{thm-trivialvan}
Let $E \in \cBF$. Suppose that the restriction of the modified derivative cocycle $\whD \colon \cRF^E \to GL(q, \mR)$ is cohomologous to the identity cocycle. Then $\chi_E(h_I) = 0$ for all $h_I$.
\end{proposition}
The idea of the proof is that if $\whD$ is cohomologous to a identity cocycle, then there is a measurable framing of the normal bundle to $\F$, which is smooth long leaves, such that the parallel transport of the Bott connection becomes the identity map. Hence, the leafwise flat classes of the Bott connection are trivial, which implies that the Weil measures are zero. 

Here is an application of this result:
\begin{theorem} 
For each monomial $h_I$ the Weil measure $h_E(h_I) = 0$ for all $E \subset \WF$. Hence, the localizations of the residual secondary classes to the wandering set of $\F$ always vanish.
\end{theorem}
{\bf Proof:}  $\whD \colon \cRF^E \to GL(q, \mR)$ is cohomologous to the identity for all $E \subset \WF$. That is,  the normal bundle $Q | \WF$ admits a measurable framing which is parallel for the Bott connection. Now apply Theorem~\ref{thm-trivialvan}. \hfill $\eop$

It follows that residual secondary classes for $\F$, such as the generalized Godbillon-Vey classes,  are supported on the non-wandering set of $\F$, hence we obtain:
\begin{corollary}
If  $\F$ has some residual secondary class $\Delta_{\F}^*(h_I \wedge c_J) \ne 0$, then  the non-wandering set $\NWF$ has positive Lebesgue measure.
\end{corollary}

Consider the Radon-Nikodym   cocycle,    $\nu \colon \GF \to \mR$,  where $\nu(\gamma) = \ln \{ \det (D\gamma)\}$. 
 \begin{theorem} [Theorem~4.1, \cite{Hurder1986}]  \label{thm-RNdecay}
 Let $E \in \cBF$, and suppose the restriction  $\nu \colon \cRF^E \to \mR$ is well-tempered. Then the Godbillon measure $g_E = \chi_E(h_1)$ vanishes.
 \end{theorem}
 In particular, this holds for the set  $E_{\nu} \subset \BR \cup \SR$ of Theorem~\ref{thm-subexp} by   Corollary~\ref{cor-subexp}.  We thus have a generalization of the Moussu-Pelletier and Sullivan Conjecture to     foliations of codimension $q \geq 1$:
  
 \begin{theorem} [\cite{HLa2000,Hurder2008b}] \label{thm-gv} Let $\F$ be a $C^1$-foliation of codimension $q\geq 1$. 
If  $E \in \cBF$ has non-zero Weil measure, $\chi_E(h_I) \ne 0$ for some $h_I$, then the intersection $E \cap \HR$ has  positive Lebesgue measure.
In particular, if some  generalized Godbillon-Vey class $\Delta_{\F}^*(h_1 \wedge c_J) \in H^{2q+1}(M; \mR)$ is non-zero, then the  leaves $\FF$  with exponential growth for $\F$ must have positive Lebesgue measure.  
\end{theorem}

 \begin{corollary} \label{thm-weilnonvan}
 Let $\F$ be a $C^2$-foliation of codimension $q$.
If some  residual secondary class $\Delta_{\F}^*(h_I \wedge c_J) \in H^*(M; \mR)$ is non-zero, then the set of leaves $\HF$ with non-trivial asymptotic expansion  has positive Lebesgue measure.  
 \end{corollary}
The conclusion that  $E \subset \HF$ has positive Lebesgue measure has strong consequences for the dynamics of the orbits in $E$. This will be discussed in \S\ref{sec-hyperbolic}.    

Finally, consider the relation between amenability and the values of the secondary classes. The following  is an ``integrated'' version of Theorem~\ref{thm-amenvan}:
  \begin{theorem} [Theorem~3.5, \cite{HK1987}] 
Let $E \in \cBF$ and suppose that the restricted equivalence relation   $\cRF^E$ is   amenable with respect to Lebesgue measure. For each monomial $h_I$ of degree $\ell > 1$, the Weil measure $\chi_E(h_I) = 0$. In particular, if some  residual  secondary class   $\Delta_{\F}^*(h_I \wedge c_J) \in H^{2q+\ell}(M; \mR)$ is non-zero, for $\ell > 1$, then there must exist a set $E \in \cBF$ with positive Lebesgue measure, such that $\cRF^E$ is non-amenable.  Hence, the von~Neumann algebra $\cM^*(\cRF)$ of the equivalence relation $\cRF$ contains a factor which is  not injective. 
 \end{theorem}
The hyperfiniteness condition is briefly discussed in    \S1.9 of \cite{HK1987}, and much more thoroughly by   Connes in Section 4.$\gamma$, pages 50--59 of \cite{Connes1994}.

We summarize these various results. First we discuss the  the Godbillon-Vey classes for codimension $q \geq 1$, so that the following is the most general answer to the Conjecture of Moussu \& Pelletier and Sullivan:

\begin{theorem} \label{thm-main3}
Suppose that $\F$ is a $C^2$-foliation with non-trivial generalized Godbillon-Vey class $\Delta_{\F}^*(h_1 \wedge c_J) \in H^{2q+1}(M; \mR)$. Then: 
\begin{enumerate}
\item the non-wandering set  $\NWF$ has positive Lebesgue measure;
\item the set of leaves with fast growth type, $\FF$, has positive   measure;
\item the set   $\HF$  has positive   measure; 
\item $\cM^*(\cR)$ contains a factor of type III.
\end{enumerate}
In fact, $\Delta_{\F}^*(h_1 \wedge c_J)$ is supported on the intersection
$\NWF \cap \FF\cap \HF \cap M_{III}$ which therefore must have positive Lebesgue measure. Moreover, for every point 
$x \in \NWF \cap \FF \cap \HF \cap M_{III}$ with positive Lebesgue density and every open neighborhood $x \in U \subset \cT$, the closure $\overline{U}$ contains points $ y \in \overline{U}$ with positive local entropy, $h(\cGF, y) > 0$.
\end{theorem}
   
   Next, we discuss the dynamical implications of the existence of non-trivial residual classes of higher degree, for $q > 1$.
   
   \begin{theorem} \label{thm-main4}
Suppose that $\F$ is a $C^2$-foliation with non-trivial residual secondary class $\Delta_{\F}^*(h_I \wedge c_J) \in H^{2q+\ell}(M; \mR)$ for $\ell > 1$.  Then: 
\begin{enumerate}
\item  $\NWF$ has positive   measure;
\item   $\HF$  has positive   measure;
\item the non-amenable component, $\KF$, has positive   measure; 
\item $\cM^*(\cR)$ contains a factor of which is not injective.
\end{enumerate}
In fact, $\Delta_{\F}^*(h_I \wedge c_J)$ is supported on the intersection
$\NWF \cap \HF \cap \AF$ which therefore must have positive Lebesgue measure. 
\end{theorem}

Thus,  $\Delta_{\F}^*(h_I \wedge c_J) \ne 0$ implies  there is a  Borel subset $E \in \cBF$ with positive Lebesgue measure such that  $E$ all points in $E$ are non-wandering,   have positive asymptotic expansion, and the normal linear holonomy cocycle 
$D \colon \cRF^E \to GL(q, \mR)$ has non-amenable algebraic hull. Given this information, the following seems to be surely true:
\begin{conjecture}
Let  $\F$ be a $C^2$-foliation. If   $\Delta_{\F}^*(h_I \wedge c_J) \in H^{2q+\ell}(M; \mR)$ is non-zero, for some class $h_i \wedge c_J$ with  $\ell > 1$ and $|J| = q$, then $h(\F) > 0$.
\end{conjecture}

There are many classes of foliations for which the dynamics are not ``chaotic'', but they  are not known to be trivial as cycles in 
$B\Gamma^r_q$ or $F\G^r_q$ for $r > 1$.
For example, if $\F$ is a foliation of codimension $q > 1$   for which all leaves are compact, so that $M = \BF$ , and the normal bundle $Q$ is framed, then it is not even known if the classifying map $h_{\F,s} \colon M \to F\G^2_q$ is homotopically trivial.  More generally, if $\F$ is a foliation with all leaves proper, then all residual secondary invariants of $\F$ vanish, yet almost nothing is known of their classifying maps.  
For example, one  question is whether there are (as yet unknown) cohomology or other homotopically-defined  invariants   of foliations which   detect the homotopy class of the classifying maps for proper foliations.  

 In    following sections, we examine several classes of non-trivial  foliation dynamics, organized by the classification scheme of Theorem~\ref{thm-decomp}. The range of examples, results and questions   provide further  motivation for the study of foliation dynamics;  the examples themselves are fascinating,  and only partly   understood at present.

\section{Elliptic foliations} \label{sec-elliptic}

A foliation $\F$ is said to be \emph{elliptic}  if $M = \bE_{\F}$. Riemannian foliations provide the most obvious examples of elliptic foliations, but there are many other types of examples.   All residual secondary classes vanish for elliptic foliations.

Let $\cS_q \subset GL(q, \mR)$ denote the convex cone of symmetric, positive definite matrices. Given $E \in \cBF$, a  Borel inner product on  $Q | E \to E$ is a Borel map $S \colon E \to \cS_q$. That is, 
 for $x \in E$,  the inner product on $T_x\cT \cong \mR^q$ is given by    $\langle \vec{v}, \vec{w} \rangle_x = \vec{v}^t \cdot S_x \cdot \vec{w}$.

We say that $\F$ is a \emph{Riemannian foliation} if there is a smooth Riemannian metric on $\cT$ so that the action of the pseudogroup $\cGF$ on $\cT$ consists of local isometries with respect to this metric \cite{Haefliger1971,Reinhart1959, Reinhart1983}. That is, there is a smooth map  $S \colon \cT \to \cS_q$ such that for all $\gamma \in \GF^{x,y}$ and all $\vec{v}, \vec{w} \in T_x\cT$ we have
\begin{equation}\label{eq-inner}
\langle \vec{v}, \vec{w} \rangle_x =  \langle D\gamma(\vec{v}), D\gamma(\vec{w}) \rangle_y ~ , ~ {\rm for~all} ~ \vec{v}, \vec{w} \in T_x\cT
\end{equation}
The choice of a continuous orthonormal framing $f \colon \cT \to GL(q, \mR)$ for the inner products $S_x$ defines a bounded cohomology of $D\colon \GF \to GL(q, \mR)$ to a cocycle with values in $O(q)$. Thus, if $\F$ is a Riemannian foliation, then  $\bE_{\F} = \cT$.

 Given $E \in \cBF$, we say that  $\F$ is a   \emph{Borel Riemannian foliation on $E$} if there a   Borel map $S \colon E \to \cS_q$ such that for $x \in E$ and all  $\gamma \in \GF^{x,y}$ then (\ref{eq-inner}) is satisfied.   This is equivalent to saying that the restricted cocycle $D^E \colon \GF^E \to GL(q,\mR)$ is Borel cohomologous to a cocycle with image in the orthogonal group $O(q)$. If  $\F$ is a   Borel Riemannian foliation on $E$, then clearly $E \subset \bE_{\F}$.   For $E = \cT$,  we say that 
 $\F$ is a  \emph{Borel Riemannian foliation}.

 We say that $\F$ is a   \emph{measurable Riemannian foliation} if there is a conull set $E\in \cBF$ such that  $\F$ is a Borel Riemannian foliation on $E$.    That is, there exists a saturated subset $X \subset M$ of Lebesgue measure zero, such that the restriction of $Q$ to $M-X$ admits a holonomy invariant Borel inner product. We mention examples later in this section of how such foliations arise in dynamical systems.

\begin{theorem}\label{thm-elliptic}
 $\F$   is a  Borel Riemannian foliation on $\bE_{\F}$.
\end{theorem}
The proof is based on a combination of standard techniques, which we briefly recall.

 Let $E \in \cBF$, and let  $S \colon E \to \cS_q$ be a Borel family of inner products.
 
 Given a linear map $L \colon T_x\cT \cong \mR^q \to T_y\cT \cong \mR^q$, represented by a matrix $A \in GL(q, \mR)$, we get an induced inner product on $T_x\cT$, 
 $$\langle \vec{v}, \vec{w} \rangle'_x = \langle A \vec{v},  A \vec{w} \rangle_y    =    \vec{v}^t A^t \cdot S_y \cdot A \vec{w}  =  \vec{v}^t \cdot S'_x \cdot \vec{w}$$
 where $S_x' = A^t S_y A$. The map $A \mapsto A^t S_x A$ defines a right action of $GL(q, \mR)$ on the symmetric matrices.
  
  Given  $\gamma  \in \GF^{x,y}$ then $D\gamma \colon T_x\cT \to T_y\cT$  which induces 
  $$D\gamma^*(S_y) = (D\gamma)^t S_y (D\gamma)  \in \cS_q$$
 Now let $S$ be the standard inner product on $\mR^q$, so that   $S_x = Id$ for all $x \in E$. 
 Then $D\gamma^*(Id)  =  (D\gamma)^t  (D\gamma)$.
  
  The assumption $x \in \bE_{\F}$ implies the set  $\{ D\gamma   \mid   \gamma \in \GF^{x}\} \subset GL(q,\mR)$
  is bounded for the norm defined by (\ref{eq-norm}) in \S\ref{sec-cocycles},  hence 
  $$\cM_x =  \{ (D\gamma)^t  (D\gamma)  \mid    \gamma \in \GF^{x} \} \subset \cS_q$$
is a bounded subset of the convex space $\cS_q$.  Let $C(x) \subset \cS_q$ denote the compact convex hull of $\cM_x$  and $S'_x\in \cS_q$ the center of mass for $C(x)$. Note that $S'_x$ depends continuously on the hull $\cM_x$ which in turn is a Borel function of $x$, as this was assumed for the given inner product $S$.

For   $\delta \in \GF^{z,x}$ and  $\gamma \in \GF^{x,y}$, then $\gamma \circ \delta \in \GF^{z,y}$, and we  calculate 
$$
D\delta^*((D\gamma)^t  (D\gamma))  =   (D\delta)^t(D\gamma)^t  (D\gamma) (D\delta) 
   =   (D(\gamma \circ \delta))^t  (D(\gamma \circ \delta))
$$
 Thus, 
$D\delta^* C(x) = C(z)$, hence $D\delta^*S_x' = S_z'$ so that $S' \colon E \to \cS_q$ satisfies  (\ref{eq-inner}). \hfill $\eop$

\medskip

One of the remarkable accomplishments in foliation theory during the 1980's was the almost complete classification of Riemannian foliations of closed manifolds  -- at least in principle \cite{Haefliger1989,Molino1988,Molino1994}. 

In contrast,  the more general class of elliptic  foliations is not well understood. First, we note that for  neither  case, when   $\F$  is a Borel nor a measurable Riemannian foliation,   does it necessarily follow that $\F$ is Riemannian foliation. The  examples below illustrate this.

There are   many examples of foliations for which there exists $E \in \cBF$ such that $\F$ is a Borel Riemannian foliation on $E$. For example, suppose that  $\WF \ne \emptyset$. Then for each $x \in \WF \cap \cT$, there is an open ball $U = B_{\cT}(x, \epsilon) \subset \WF$ such that   all $g \in \cGF$ with $x \in D(g)$ and $[g]_x \ne Id$, then $g(U  \cap D(g)) \cap U = \emptyset$. Select an inner product  $S_{U} \colon U \to \cS_q$, then its translates via the relation (\ref{eq-inner}) defines an inner product on the saturation, $S \colon U_{\cR} \to \cS_q$. This yields   a holonomy invariant Riemannian metric on $U_{\F}$, so that $U_{\F} \subset \EF$. The metric on $U_{\cR}$ extends to a Borel family of inner products on $\cT$, which if the orbit $\cO(x)$ is infinite,  cannot be extended or modified to yield   a smooth Riemannian metric on $\cT$ if $\cO(x)$ is infinite. To see this, note that  a invariant  continuous    Riemannian metric yields   a $\cGF$-invariant continuous volume form on $\cT$,  for which $\cT$ has finite total volume. But if $\cO(x)$ is an infinite orbit, then    $U_{\cR}$ has infinite volume for any $\cGF$-invariant continuous volume form, which would yield a contradiction. 

A  Denjoy foliation $\F$ of $\mT^2$ provides a concrete illustration of the above remarks.  (See  page 108, \cite{CanCon2000} or page 26,  \cite{Tamura1992} for constructions and properties of the Denjoy examples.) We suppose $\F$ has an exceptional minimal set $K \subset \mT^2$, whose complement $U = \mT^2 \setminus K$ is an open foliated product. The leaves in the complement $U$ admit a cross-section, hence lie in $\bE_{\F}$. On the other hand, the leaves in the minimal set must also lie in $\bE_{\F}$, as otherwise there would exists a fixed-point for $\cGF$ with linearly contracting holonomy   \cite{SackstederSchwartz1964,Sacksteder1965,Hurder1991a}. Thus, the Denjoy examples are elliptic foliations, but obviously are not Riemannian foliations.
 
The Reeb foliation $\F$ of $\mS^3$ is an elliptic foliation.  There is one compact toral leaf, $L \cong \mT^2$,  all of whose holonomy is one sided. This implies that the linear holonomy of the compact leaf is trivial, hence $L \subset \bE_{\F}$. The leaves of $\F$ in each component of the  complement  $ \mS^3 \setminus L$ fiber over a closed transversal, hence by the above remarks, they also belong to $\bE_{\F}$. 

On the other hand, one can construct  Reeb foliations of $\mT^2$, for which there are closed  leaves of $\F$ with two-sided holonomy given by linear contractions. Such circle leaves with contracting linear holonomy clearly lie in $\HF$, and in fact, their union is all of $\HF$.  Their  complement consists of   the  proper leaves, along with the  closed leaves with trivial linear holonomy, which are all contained in $\EF$.
 
A key point about the decomposition in Theorem~\ref{thm-decomp} is that  the components are Borel sets, but not necessarily closed.  The study of the leaves in the boundaries of the sets $\EF$, $\PF$ or $\HF$ reveals key aspects of the dynamics of $\F$.

A foliation $\F$ is said to be    \emph{compact} if   every leaf of $\F$ is compact \cite{EMS1977,Epstein1972,Epstein1976,EV1978}. 
 For codimension $2$, Epstein proved that the leaf space of     a compact foliation is a compact orbifold \cite{Epstein1972}. On the other hand, for codimension $q > 2$, Sullivan \cite{Sullivan1975b}, Epstein-Vogt \cite{EV1978}, and Vogt \cite{Vogt1976,Vogt1977a,Vogt1977b}    constructed compact foliations for which   the leaf space is not  Hausdorff. 
 The exceptional set $E(\F) \subset M$ of a compact foliation $\F$ is defined to be the union of all leaves of $\F$ whose holonomy groups are \emph{infinite}: that is,  $L_x \subset E(\F)$ if and only if $\GF^{x,x}$ is infinite. The exceptional set is a closed Borel subset without interior. The  complement, $G(\F) = M \setminus E(\F)$, is   the \emph{good set} of $\F$, and the restriction of $\F$ to $G(\F)$ is a Riemannian foliation, hence $G(\F) \subset \bE_{\F}$. We also note from the definitions and the fact that every $L \subset E(\F)$ has non-trivial holonomy,  $\cW(\F) \subset  G(\F)$.
 
Suppose  $\F$ is a compact foliation such that every leaf $L_x \subset E(\F)$ has finite linear holonomy, then $\F$ is elliptic. 
Otherwise, suppose that $L_x \subset E(\F)$ be such that $D \colon \GF^{x,x} \to GL(q, \mR)$ has infinite image, and let  $\gamma \in \GF^{x,x}$ such that  $D\gamma \in GL(q, \mR)$    is non-trivial . Then eigenvalues of   $D\gamma$ must be norm one, or else by Theorem~\ref{thm-hypfixpt2} below, there exists an attracting orbit for $\gamma$ and so this orbit is not a compact leaf.  Thus,  
$D\gamma$  is conjugate to a parabolic matrix. It follows that  the exceptional set $E(\F) \subset \bE_{\F} \cup \bP_{\F}$.
Moreover,  by Proposition~\ref{prop-linholo}, the set of leaves in $E(\F)$ with non-trivial  linear holonomy must have Lebesgue measure zero.  Thus,     we have:
 \begin{proposition} \label{prop-compact}
 Let $\F$ be a compact foliation of a closed manifold $M$, whose leaf space $M/\F$ is not Hausdorff. Then $\F$ is parabolic, and moreover $\F$ is a measurable Riemannian foliation, which is not Riemannian.
 \end{proposition}
Measurable Riemannian foliations were studied by Zimmer \cite{Zimmer1982b,Zimmer1991} in the context of a smooth action of a lattice group $\G$ on a closed manifold, $\alpha \colon \G \times N^q \to N^q$. 
The suspension of the action yields a foliation $\F_{\alpha}$, whose  pseudogroup $\cG_{\F_{\alpha}}$ is equivalent to that defined by the action of $\G$ on $N$.  
  Assume that the action $\alpha$  preserves some smooth measure on $N$. Then, with additional hypotheses, such as assuming that the real-rank of $\G$ is sufficiently large and every $\gamma \in \G$ defines a diffeomorphism with zero topological entropy, Zimmer proved that  there exists a measurable Riemannian metric on $TN$ which is invariant under the group action  \cite{Zimmer1982b,Zimmer1985,Zimmer1991}. Thus, $\F_{\alpha}$ is a measurable Riemannian foliation. 
Zimmer also gave conditions for when a  measurable Riemannian foliation admits  an invariant smooth metric   \cite{Zimmer1985}
(see also \cite{Benveniste2000, FisherZimmer2002}.)  An alternate approach to this result, which applies more generally  groups $\G$ which have Property T, was given by   Fisher and   Margulis  \cite{FisherMargulis2005}.    It remains an open problem to determine whether there can exist a smooth action of a higher-rank lattice, which preserves a measurable Riemannian metric but no smooth Riemannian metric.

\section{Parabolic foliations} \label{sec-parabolic}

A foliation $\F$ is said to be \emph{parabolic}  if $M = \EF \cup \PF$ with $\PF \ne \emptyset$. Parabolic foliations are \emph{almost-isometric}. Every distal foliation is either elliptic or parabolic.     There are a variety of constructions of parabolic foliations, but little is known of their classification. All residual secondary classes vanish for parabolic foliations.

For   $E \in \cBF$, we  say that $\F$ is \emph{almost isometric} on the saturation $E_{\F}$  if the restricted cocycle 
 $D \colon \GF^E \to GL(q,\mR)$ is \emph{well-tempered} (see Definition~\ref{def-welltemp}). 

It   follows directly from the definition of the sets $\ER$ and $\PR$ that for $E \subset \ER \cup \PR$,  $D \colon \GF^E \to GL(q,\mR)$ has exponential type $0$,  so by Theorem~\ref{thm-wt2} we have:
\begin{theorem}[Hurder \cite{Hurder2008b}] \label{thm-almostisoE}
If  $E \in \cBF$ with $E \subset \ER \cup \PR$, then $\F$ is almost isometric on $E_{\F}$.
\end{theorem}
 
 \begin{corollary}  \label{cor-almostiso2}
Let $E \in \cBF$ with $E \subset \ER \cup \PR$.  Then for each Weil measure,  $\chi_E(h_I) = 0$.
 Thus, given any residual secondary class $h_I \wedge c_J$, the localization   $\Delta_{\F}^*(h_I \wedge c_J)|E  \in H^*(M; \mR)$ is zero.
In particular, if $\F$ is parabolic, then all residual secondary classes for $\F$ vanish.
 \end{corollary}

Corollary~\ref{cor-almostiso2} is a generalization of an early result of Michel Herman, who showed that the Godbillon-Vey class vanishes for a $C^2$-foliation $\F$ on $\mT^3$ by planes 
in \cite{Herman1977}. His method was to observe that by Sacksteder \cite{Sacksteder1965}, such a foliation is equivalent to a suspension of a $\mZ^2$-action on the circle $\alpha \colon \mZ^2 \times \mS^1 \to \mS^1$. He then showed that while the action need not admit a smooth invariant measure on $\mS^1$, for all $\epsilon > 0$  there always exists an $\epsilon$-invariant, absolutely continuous  measure $\mu_{\epsilon}$ on $\mS^1$ which is equivalent to Lebesgue measure, and hence $GV(\F) = 0$. The measure $\mu_{\epsilon}$ defines an   $\epsilon$-invariant metric for the 1-dimensional normal bundle to $\F$, so this result follows from the Corollary. In fact, Herman's    method   foreshadowed the entire development of vanishing theorems for the Godbillon-Vey  classes. 

 Another generalization of the Herman Vanishing result  is the following:
 \begin{theorem}\label{thm-conjugate}
Let $\F'$ be a Riemannian  foliation of a closed manifold $M'$. Given a $C^r$-foliation $\F$ of codimension-$q$ on $M$, suppose there exists  a homeomorphism $h \colon M \to M'$ which maps the leaves of $\F$ to the leaves of $\F'$. If either $q =1$ and $r \geq 1$, or $q > 1$ and $r > 1$, then    $\F$ is parabolic. 
\end{theorem}
 
 The proof of Theorem~\ref{thm-conjugate}  is a consequence of a more general result:
  \begin{theorem}\label{thm-distal}
Let $\F$ be   a $C^r$-foliation $\F$ of codimension-$q$ on $M$, with  either $q =1$ and $r \geq 1$, or $q > 1$ and $r > 1$. If $\cGF$ is distal, then     $\F$ is parabolic. 
 \end{theorem}

 The proof follows from Corollary~\ref{cor-hypfixpt} in the next section, which implies that if $\HR \ne \emptyset$, and either $q =1$ and $r \geq 1$, or $q > 1$ and $r > 1$, then $\cGF$ has a proximal orbit. A distal foliation cannot have a proximal orbit. A foliation conjugate to a Riemannian foliation is distal, which yields Theorem~\ref{thm-conjugate}.
 
 The second class of examples of parabolic foliations is obtained from Lie group actions. 
 Let $G \subset GL(m, \mR)$ be a closed, connected Lie subgroup and $\Lambda \subset G$ a discrete, torsion-free cocompact subgroup. Thus $M = G/\Lambda $ is a closed manifold. Let $H \subset G$ be a connected subgroup such that, as a subgroup of $GL(m,\mR)$,
  every $A \in H$ is a matrix with all eigenvalues of modulus $1$. That is, $H$ is a parabolic subgroup of $GL(m,\mR)$. 
 Let $\F$ be the foliation on $M$ whose leaves are the orbits of  the left action of $H$.
The dynamical properties of this class of locally homogeneous  foliations have been extensively studied, as they are exactly the foliations which arise in the work of Ratner  \cite{Ghys1992,Morris2005}.

 Let $\mathfrak{g}$ denote the Lie algebra of $G$ of left-invariant vector fields on $G$. Let $\mathfrak{h} \subset \mathfrak{g}$ be the Lie subalgebra corresponding to $H$. Let $\mathfrak{m} \subset \mathfrak{g}$ be a complementary subspace to $\mathfrak{h}$, of dimension $q$,  which we identify with the quotient space $\mathfrak{g}/\mathfrak{h}$.  
 The Adjoint action of $H$ on $\mathfrak{g}$ leaves the subspace $\mathfrak{h}$ invariant, hence induces a unipotent representation $h \colon H \to {\rm Aut}(\mathfrak{m})$. The holonomy cocycle $D \colon \cGF \to GL(q, \mR) \cong {\rm Aut}(\mathfrak{m})$ for $\F$  is conjugate to a cocycle with values in the range of $h \colon H \to {\rm Aut}(\mathfrak{m})$, have all elements of holonomy are unipotent. Hence $\HF = \emptyset$. 
The terminology ``parabolic foliation'' is motivated by such  examples.
 
Compact foliations provide a  third class of examples  of parabolic foliations (see Proposition~\ref{prop-compact}).
 
Finally, we describe one further  class of examples, based on  an explicit construction  of distal foliations which are not homogeneous. 
Hirsch constructed  in     \cite{Hirsch1975}  an analytic foliation $\F$ of codimension-one with an exceptional minimal set on a closed $3$-manifold $M$, starting from a familiar method in dynamical systems to construct   diffeomorphisms of compact manifolds with expanding isolated invariant sets which are solenoids  \cite{Williams1967,Williams1974}. 
        The construction of the Hirsch foliation was generalized in   \cite{BHS2006} to   codimension $q > 1$. 
        
        As the   Hirsch foliation $\F$ has codimension-one, there is a transverse vector field $\vec{X}$ which defines a foliation transverse to $\F$. A surprising fact is that one can always choose this vector field so that it defines a parabolic foliation, which contains a solenoidal minimal set for the flow \cite{CH2008}.  This construction is part of a more general method to embed solenoids as minimal sets for foliations. Here is a typical result:
 \begin{theorem}[Clark-Hurder\cite{CH2008}] \label{thm-solenoids}
 Let $\F$ be a $C^1$-foliation with codimension $q > 1$.  
 Suppose that $\F$ has a compact leaf $L$ with $H^1(L, \mR) \ne 0$, and there is a saturated open neighborhood $L \subset U$ such that $\F \mid U$ is a product foliation. Then there is an arbitrarily small smooth perturbation $\F'$ of $\F$ such that $\F'$ has a solenoidal minimal set $\K \subset U$, where the leaves of $\F' \mid \K$ all cover $L$. Moreover, if $\F$ is distal, then $\F'$ is distal.
 \end{theorem}
 Thus, one can introduce solenoidal minimal sets into a wide variety of foliations, starting for example with product foliations, to obtain parabolic foliations.  It seems plausible that parts of the Williams classification theory for expanding attractors
 \cite{Williams1967,Williams1974}
can be used to construct even more general families of parabolic foliations, as varied as are the types of these attractors. 

In general, the class of parabolic foliations seems quite broad, and quite unknown.

\section{Hyperbolic foliations} \label{sec-hyperbolic}

The set $\HF$ consists  of the leaves of $\F$ with whose holonomy exhibits at least partial   hyperbolicity. For example, as remarked previously, an attracting closed leaf for a Reeb foliation of $\mT^2$ is in $\HF$. The non-closed leaves of a Reeb foliation on $\mT^2$ are proper, so   $\HF$ non-empty does not necessarily imply   ``chaos'' in the dynamics of $\F$. On the other hand, we have seen that:
\begin{theorem} Let $\F$ be a $C^2$-foliation of   a closed manifold $M$. If some residual secondary class $\Delta_{\F}(h_I \wedge c_J) \ne 0$, then  $\HF$ has positive   Lebesgue measure.   
\end{theorem}

This emphasizes the need  understand the dynamics of foliations which have non-uniformly, partially   hyperbolic behavior on a set of positive measure. This problem is wide-open.
 For codimension-one foliations, there are many partial results, as  discussed in \S\ref{sec-q1}.
   In  this section, we  give some general results about the relation between $\HF$ and the dynamics of $\F$, valid in arbitrary codimension.

Let $x \in \HF$ with $e(x) = e(\GF , x) > 0$, then there exists $0 <  \lambda \leq e(x)$ and some sequence of elements $\gamma_m \in \GF^x$ with $\|\gamma_m\|_x \to \infty$, and  such that either 
$$ \frac{\ln \{\|D\gamma_m\|\}}{\|\gamma_m\|_x} \to \lambda ~ , ~ {\rm or} ~ \frac{\ln \{\|D\gamma_m^{-1} \|\}}{\|\gamma_m\|_x} \to  - \lambda$$

The  Pesin Theory approach to analyzing hyperbolic behavior of a diffeomorphism requires the existence of an invariant measure, in order to ensure  uniform recurrence for ``typical''  hyperbolic orbits in the support of the measure.   This is fundamental for using 
  the infinitesimal information from the derivative cocycle  to obtain    dynamical conclusions. 
There are several difficulties with applying this method   to foliations with leaf dimension greater than one, among them the fact that $\cGF$  need not have any $\cGF$-invariant   Borel probability measures. Even if such a measure $\mu$ exists, one does not know how to assure that it is typical for the transverse expansion; that is, with $e(x) > 0$ for $\mu$-almost every $x$. One solution to this difficulty was introduced in \cite{Hurder1988}, based on the leafwise geodesic flow.

The choice of a Riemannian metric on $TM$ induces a Riemannian metric on each leaf $L$ of $\F$. We assume (without loss of generality)  that the leaves of $\F$ are smoothly immersed submanifolds, even if $\F$ is only transversally $C^r$ for some $r \geq 1$.   
Then for each $L$, there is a leafwise geodesic map, $\exp_L \colon  TL \to L$, associating to $\vec{v} \in T_xL$ the point     
$\exp_L(\vec{v}) \in L$. The map $\exp_L$ depends continuously on the choice of the leaf, so yields   a flow on $T\F$ by
\begin{equation}\label{eq-geodesic}
\varphi^{\F}_{t_0}(x,\vec{v}) =  (\exp_{L_x}(t_0 \cdot \vec{v}) , \frac{d}{dt} \exp_{L_x} (t \cdot \vec{v}) |_{t=t_0}
\end{equation}
Let $\whM \subset TF$ denote the unit sphere bundle for the tangent vectors to the leaves. Let $\wF$ denote the foliation on $\whM$ obtained by pulling back the leaves of $\F$ via the bundle projection, $\pi \colon \whM \to M$. That is, a leaf of $\wF$ is simply the unit  tangent bundle to a leaf of $\F$. Since the speed of a geodesic is constant, the restricted  flow $\varphi^{\F} \colon \mR \times \whM \to \whM$, 
$\varphi^{\F}_t \colon \whM \to \whM$ is well-defined.

The leafwise geodesic flow $\varphi^{\F}_t$ was used by Walczak \cite{Walczak1988} to study the curvature and Lyapunov spectrum of the leaves of a foliation. 
Our interest is derived from the following basic observation:
\begin{proposition} [Proposition~5.1, \cite{Hurder1988}] 
The derivative cocycle   lifts to a cocycle over the leafwise geodesic flow,  
\begin{equation}
D^{\nu} \equiv (D\varphi^{\F})^{\perp} \colon \mR \times \whM \to GL(q,\mR)
\end{equation}
\end{proposition}
The flow $\varphi^{\F}_t$ preserves the leaves of $\wF$,  which  implies that it induces an action on the normal bundle to $\wF$, which is identified with the pull-back $\pi^*Q$. 

Given  $\gamma \in \GF^{x,y}$ we can write it as a product of generators, 
$\gamma = g_{i_m} \circ \cdots \circ g_{i_1}$ as in  Definition~\ref{def-gamma}. The corresponding plaque-chain then defines a homotopy class of leafwise paths from $x$ to $y$, so there is  a (not necessarily unique)  leafwise geodesic $\exp_{L_x}(t \cdot\vec{v})$ for some unit vector $\vec{v_y} \in T_xL_x$
 which satisfies   $y = \exp_{L_x}(t_y \cdot \vec{v_y})$ and is endpoint homotopic to  the chosen plaque-chain.  The germinal transverse holonomy along $\exp_{L_x}(t \cdot \vec{v_y}) \colon [0,t_y] \to L_x$ 
  equals $\gamma$, hence    $D\gamma = D^{\nu} (t_0, \vec{v_0})$.

The advantage of this construction is that $D^{\nu}$ is a linear cocycle over the flow $\varphi_{\F}$ so we can apply the usual methods of the Oseledets Theory \cite{BarreiraPesin2007,Katok1980,KatokHasselblatt1995,Oseledets1968,Pesin1977}.

\begin{definition} [Theorem~5.2, \cite{Hurder1988}] \label{thm-lyapspectrum}
Let $\mu$ be an ergodic invariant    measure on $\whM$ for the flow $\varphi_{\F}$. The \emph{transverse Lyapunov spectrum} of $\mu$ is the set of exponents (logs of generalized eigenvalues) for the cocycle $D^{\nu}$ with respect to $\mu$:
\begin{equation}\label{eq-lyapunov}
\Lambda^{\mu} \equiv \{\lambda_1^{(\mu)}  < \lambda_2^{(\mu)} < \cdots < \lambda_s^{(\mu)}\}  
\end{equation}
   Note that the integer $s$ satisfies $1 \leq s \leq q$, and   the set $\Lambda^{\mu}$  depends on the  choice of  ergodic invariant measure $\mu$ as indicated. 
\end{definition}
The numbers $\{\lambda_i^{(\mu)} \mid 1 \leq i \leq s\}$ are called  the \emph{transverse Lyapunov exponents} for $\varphi^{\F}_t$ with respect to $\mu$. They are the transverse rates of expansion of the linear holonomy for  $\F$ in the direction of a generic leafwise geodesic in the support of $\mu$. As such, they clearly have a strong influence on the dynamics of $\cGF$ acting on $\cT$.

The following result justifies calling $\HF$ the set of partially hyperbolic leaves:

\begin{theorem} [Hurder \cite{Hurder1988,Hurder2000b, Hurder2008b}] \label{thm-stable}
Let $\F$ be a $C^1$-foliation. Given  $x \in \HF$,   there exists an ergodic invariant measure $\mu$ for $\varphi_{\F}$ supported on the set $\widehat{\cO(x)} = \pi^{-1}(\overline{\cO(x)})$ such that the largest transverse Lyapunov exponent $\lambda_s^{(\mu)} > 0$. 

Conversely, if $\mu$ is an ergodic invariant measure for $\varphi^{\F}_t$ with some $\lambda_i^{(\mu)} \ne 0$, then for every $(x, \vec{v}) \in \whM$ in the support of $\mu$, the point $x \in \HF$. 
\end{theorem}

Moreover, if $\F$ is $C^r$ for some $r > 1$, then there always exists a transverse stable or unstable manifold for this maximal exponent. More generally, there are stable/unstable  manifolds for all transverse exponents $\lambda_i^{(\mu)} \ne 0$. The collection of all such ergodic invariant measures $\mu$ and stable/unstable manifolds for the leafwise geodesic flow gives deep information about the dynamical properties of $\F$ on $\HF$. 
The problem is how to use this data to prove particular dynamical properties of $\F$. Let us describe  two cases where this has been achieved for foliation dynamics.

An ergodic invariant measure $\mu$ for $\varphi^{\F}_t$ is said to be \emph{normally hyperbolic}  if every transverse exponent 
$\lambda_i^{(\mu)} \ne 0$, $1 \leq i \leq s$.   
For example, given  a diffeomorphism $f \colon \Sigma_g \to \Sigma_g$ of a closed Riemann surface,  
Katok   \cite{Katok1980} used  the Pesin Theory of normally hyperbolic invariant measures to prove that the topological entropy of $f$ equals the growth rate of its periodic orbits. Similarly, for   foliation dynamics, Pesin Theory applied to a normally hyperbolic measure for the foliation geodesic flow $\varphi^{\F}_t$ can be used to construct periodic orbits;  the problem is to show the existence of normally hyperbolic ergodic  measures.

 For $1$-dimensional dynamical systems,    stable manifolds are just attracting domains, so that the maps are only required to be $C^1$.  For foliations of codimension-one, a $\varphi^{\F}_t$-invariant  measure is normally hyperbolic is it has non-zero exponent. In this case,  the author showed in  \cite{Hurder1988,Hurder1991a,Hurder2000b}:
 
\begin{theorem} [\cite{Hurder1991a}]  \label{thm-hypfixpt1} 
Let $\F$ be a $C^1$-foliation of codimension-one.  Suppose that $e(x) > 0$. Then for all $0 < \lambda < e(x)$,  there exists $y \in \omega(x)$ and $\gamma \in \GF^{y,y}$ such that $\ln | D_y\gamma | > \lambda$. That is,  $L_y$ is a leaf with expanding  linear holonomy.
\end{theorem}

An invariant measure $\mu$ for the foliation geodesic flow $\varphi^{\F}_t$ is said to be \emph{partially hyperbolic} is at least one of the exponents $\lambda_i^{(\mu)} \ne 0$. Even in this case, one can obtain dynamical consequences from the infinitesimal data. Using  the prevalence of  typical points for a partially hyperbolic measure $\mu$,   one constructs a discrete model for a   leafwise geodesic ray $\exp^{\F}(t \cdot \vec{v})$, $t \geq 0$,  along which some Lyapunov exponent is positive.  

\begin{definition}
An \emph{orbit ray} at $y$ for $\cGF$ is a mapping $g_r \colon \mN \to \cGF^{(1)}$ so that for each $\ell > 0$, 
  the composition $\varphi_r(\ell) = g_{r(\ell)} \circ g_{r(\ell -1)} \circ \cdots \circ g_{r(1)}$ is defined, with $y \in D(\varphi_r(\ell))$, and 
  $\| [\varphi_r(\ell-1)]_y \|_y \leq \| [\varphi_r(\ell)]_y \|_y \leq \ell$. 
 \end{definition}
 
This is the combinatorial version of saying that the maps $\varphi_r(\ell)$ are tracing out the flow-boxes crossed  by a leafwise geodesic ray $\exp^{\F}(t \cdot \vec{v})$, $t \geq 0$.

\begin{theorem}[\cite{Hurder2000b}] \label{thm-hypfixpt2} 
Let $\F$ be a $C^r$-foliation of codimension $q$, for $r \geq 1$. (If $q > 1$ then we require that $r > 1$.) Suppose that $e(x) > 0$. Then for all $0 < \lambda < e(x)$,  there exists $y \in \omega(x)$,  an   orbit ray at $y$,   $g_r \colon \mN \to \cGF^{(1)}$,  and a $C^1$-curve $\sigma \colon (-\delta, \delta) \to D(\varphi_r(1))$ with $0 \ne   \vec{X}  = \sigma'(0) \in T_y\cT$ such  that  
$\varphi_{r}(\ell -1)(\sigma(-\delta, \delta)) \subset D(\varphi_r(\ell))$ for all $\ell > 0$, and each $ \varphi_r(\ell)$ is a contraction along $\sigma$ with 
$$ \frac{\ln \| D_y \varphi_r(\ell)(\vec{X}) \|}{\ell} \to - \lambda $$
\end{theorem}
For codimension $q > 1$, we require that  $\F$ be $C^r$ for some $r > 1$, in order to apply the stable manifold theory to the flow $\varphi^{\F}_t$.

The conclusion of Theorem~\ref{thm-hypfixpt2} is simply that $\sigma$ is a partial stable manifold of the collection of maps $\{\varphi_r(\ell) \mid \ell = 1, 2, \ldots\}$.  The technical nature of the statement is due to the fact that the domains of the compositions $\varphi_r(\ell)$ are not necessarily fixed, and may in fact be shrinking as $\ell \to \infty$.  However, for all $\ell \geq 1$ the arc image of $\sigma$ is contained in their domains, and along this arc the maps are linear contractions at the basepoint $y$.

Together, Theorems~\ref{thm-hypfixpt1} and \ref{thm-hypfixpt2}   yield: 

\begin{corollary}\label{cor-hypfixpt} Let $x \in \HF$. Suppose that 
$q =1$ and $r \geq 1$, or $q > 1$ and $r > 1$, then $\cGF$ has a   proximal pair $y,z \in \cT$ with $y \in \omega(x) \subset \overline{\HF}$.
\end{corollary}

We conclude this section with a discussion of   examples   which have $M = \HF$.

The  Roussarie foliation $\F$ is the weak-stable foliation of the usual geodesic flow for a compact Riemann surface $\Sigma_g$ of genus $g > 1$ (see Lawson  \cite{Lawson1975} for a nice discussion of this example.) Every leafwise geodesic for $\F$ includes in its limit set  a transversely hyperbolic invariant measure for the associated leafwise geodesic flow. More generally, if $B$ is a closed manifold of dimension $m=q+1$ with a metric of uniformly negative sectional curvatures, then the usual geodesic flow defines a weak-stable foliation $\F$ on the unit tangent bundle $M = T^1B$.  The foliation $\F$ is always $C^1$ by the transverse Stable Manifold Theorem of  Hirsch, Pugh and Shub \cite{HPS1977}. It is standard that every  leafwise geodesic for $\F$ includes in its limit set  a transversely hyperbolic invariant measure for the associated leafwise geodesic flow. In this case, the transverse Lyapunov spectrum is just the transverse  part of the Lyapunov spectrum for the usual geodesic flow on $M$. The negative curvature hypothesis implies $M = \HF$.

Another class of examples is obtained from the actions of higher-rank lattices. Let $\G$ be such a group, and suppose there exists a closed Riemannian manifold $N$ and volume preserving smooth action $\alpha \colon \G \times N \to N$. The suspension of $\alpha$   yields a foliation denoted by $\F_{\alpha}$, whose pseudogroup $\cG_{\F_{\alpha}}$ is equivalent to the pseudogroup defined by the action  $\alpha$. Then there is the very strong dichotomy:
\begin{theorem}[Zimmer,  \cite{Zimmer1991,FisherZimmer2002}]  Either $\alpha \colon \G \times N \to N$ preserves a measurable Riemannian metric on $TN$, or the action has non-trivial Lyapunov spectrum almost everywhere. Hence, either $\F_{\alpha}$ is measurably Riemannian, or $\HF$ has full measure. 
\end{theorem}
 
 There are a variety of actions of lattices on manifolds with non-zero Lyapunov exponents, whose suspension foliations satisfy $M = \HF$ -- all of the examples in \S7 of \cite{Hurder1992} are of this type.

These examples suggest a general problem. Let $\F$ be a $C^1$-foliation, and suppose $E \in \cBF$ is ergodic for $\cRF$. 
Then Zimmer proved in \cite{Zimmer1990} that there exists a minimal algebraic subgroup $H = H(\F,  E) \subset GL(q,\mR)$ such that the strict normal derivative cocycle 
$\whD \colon \cRF^E \to GL(q, \mR)$ is cohomologous to a cocycle 
$$\phi_H \colon \cRF^E \to H \subset GL(q, \mR)$$ 
The algebraic subgroup $H(\F,  E)$ is called the \emph{algebraic hull} of $\whD$ on $E$, which  is well-defined up to conjugacy in $GL(q, \mR)$.  For example, if $H(\F,E)$ is compact, then $\F$ is measurably Riemannian on $E_{\F}$. At the other extreme, one can ask,
if the   algebraic hull $H(\F,E) = GL(q,\mR)$, must the closure $\overline{E_{\F}}$   contain the support of a transversely hyperbolic measure for $\varphi^{\F}_t$? 

We conclude with one more family of examples. 
Let $N$ be a closed $q$-dimensional manifold. Suppose there exists a collection of smooth maps $\mathfrak{F} = \{f_i \colon N \to N \mid 1 \leq i \leq k\}$ such that each $f_i \colon N \to N$ is a covering map. We call this a \emph{system of \'{e}tale correspondences} in \cite{BHS2006}. Then there exists a codimension-$q$ foliation $\F_{\mathfrak{F}}$ of a closed manifold $M$,  such that its holonomy pseudogroup $\cR_{\F_{\mathfrak{F}}}$ is equivalent to that generated by the collection of maps $\mathfrak{F}$. The foliation $\F_{\mathfrak{F}}$ is constructed using the generalized suspension construction, as described in \S5, \cite{BHS2006}.

For example, if $N = \mS^1$ and $f_1 \colon \mS^1 \to \mS^1$ is a covering map of degree $2$, then this yields the Hirsch foliations constructed in \cite{Hirsch1975}. For this reason, the foliations $\F_{\mathfrak{F}}$ are called \emph{generalized Hirsch foliations}. 

If at least one of the maps $f_i$ in the collection $\mathfrak{F}$ is expanding, then $M = \bH_{\F_{\mathfrak{F}}}$ and there are many transversely hyperbolic measures for the leafwise geodesic flow. The hypothesis that some $f_i$ is expanding implies that the fundamental group of $N$ admits a nilpotent subgroup of finite index \cite{Gromov1981}. More generally, any $C^1$-perturbation of such an expanding map remains expanding, so in this way, one obtains   a wide variety of examples with transversely hyperbolic measures. Perturbations of this type  are discussed in \cite{BHS2006}, and the theory of semi-Markovian minimal sets is developed  in \cite{BisHurder2008}, of which these examples are typical. The semi-Markovian minimal sets are a generalization of the examples studied by Matsumoto in \cite{Matsumoto1988}.

\section{Foliation dynamics in codimension one} \label{sec-q1}

The topological dynamics of codimension-one,   $C^2$-foliations  have been studied for almost    50 years (or in the case of flows in the plane, their topological study was started by Poincar\'e more than 100 years ago \cite{Poincare1881}.)   In this section, we give   applications of the methods in this paper to the codimension-one case, and discuss some of the  new insights it yields into this well-developed theory.

  The construction of foliations with exceptional minimal sets by Sacksteder and Schwartz  \cite{SackstederSchwartz1964}, and Sacksteder's famous paper \cite{Sacksteder1965} on the existence of resilient leaves  in exceptional minimal sets  \cite{Sacksteder1965} showed  that the dynamics of codimension-one $C^2$-foliations can exhibit robustly chaotic  dynamics. 
  Rosenberg and Roussarie \cite{RosenbergRoussarie1970c} gave   constructions of  analytic foliations with exceptional minimal sets, a result which seems unremarkable now, but pointed the field towards the study of this dynamical phenomenon.
  
There followed during the 1970's a period of rapid development. Themes included    the deeper understanding  of foliations almost without holonomy, which are natural generalizations of the dynamics of the Reeb foliation;  the understanding of asymptotic properties of leaves, which generalized the Poincar\'{e}-Bendixson theory of flows in the plane;  existence and consequences of non-trivial holonomy for leaves, and the relation between the growth of leaves and foliation dynamics, the phenomenon first seen in the properties of resilient leaves. 

Notable advances included   Hector's work on classification and examples, starting with his Thesis and subsequent developments of its themes
\cite{Hector1972a,Hector1972b,Hector1976,Hector1983};  Lamoureux's work on holonomy and ``captured leaves'' 
\cite{Lamoureux1973,Lamoureux1974b,Lamoureux1976b,Lamoureux1977}; Moussu's study of foliations almost without holonomy
\cite{Moussu1971,MoussuPelletier1974}; Nishimori's study of the asymptotic properties  and growth of leaves in foliations \cite{Nishimori1975a,Nishimori1975b,Nishimori1977,Nishimori1979}; and Plante's study of the relation between growth of leaves and the fundamental groups of the ambient manifolds.  

The study of leaves at finite level by Cantwell and Conlon    \cite{CantwellConlon1978,CantwellConlon1982,CantwellConlon1983,Inaba1979} and Tsuchiya 
\cite{Tsuchiya1979a,Tsuchiya1979b,Tsuchiya1980a,Tsuchiya1980b} explored   the relation between the hierarchy of leaf closures and their growth rates.  
This study reached its culmination in the Poincar\'{e}-Bendixson Theory of levels for $C^2$-foliations  developed by Cantwell and Conlon \cite{CantwellConlon1981b,CantwellConlon1988b} and Hector \cite{Hector1983}. 
  Poincar\'{e}-Bendixson Theory gives a framework for understanding the dynamics of  $C^2$- foliations    in codimension-one. 

This extensive list of works is certainly not complete, but gives a sense of the research activities in the field. It would take a separate survey to do justice to all of the works of this era, which are discussed in detail in the  books by Hector and Hirsch \cite{HecHir1981}, Godbillon \cite{Godbillon1991},  Tamura \cite{Tamura1992} and Candel and Conlon \cite{CanCon2000,CanCon2003} 
which are the current references for this subject. 
 
There were two notable conjectures in the field which remained unsolved by the end of the decade of the 1970's, and whose (partial) solutions in the 1980's directed research towards techniques that included  methods of ergodic theory as well as topological dynamics.   First, Hector posed the following in his thesis:

\begin{conjecture}[Hector \cite{Hector1972a,Schweitzer1978}]\label{conj-miniset}
Let $\F$ be a codimension-one, $C^2$-foliation of a closed manifold $M$. If $Z$ is an exceptional minimal set for $\F$,  then $Z$ has Lebesgue measure zero. 
\end{conjecture}
This problem  remains open in this generality, although many partial results are now known. The other conjecture was stated already, although we recall it here.

\begin{conjecture}[ Moussu-Pelletier \cite{MoussuPelletier1974};  Sullivan \cite{Schweitzer1978}] \label{conj-GV}
Suppose that $\F$ is a  codimension-one, $C^2$-foliation of a closed manifold $M$, with $GV(\F) \in H^3(M; \mR)$ non-zero. Then the set  of leaves of $\F$ with exponential growth is non-empty (or better, has positive Lebesgue measure.)
\end{conjecture}

Key to the study of   both  conjectures is the understanding of leaves of $\F$ with attracting linear holonomy, hence to the properties of $\F$ in an open neighborhood of the closure of the hyperbolic set $\HF$. We recall a concept in topological dynamics which is fundamental to this study.
 
\begin{definition} \label{def-markov}
Let $\F$ be a $C^1$-foliation with codimension $q \geq 1$.  A \emph{Markov sub-pseudogroup} (or more simply, a \emph{Markov system}) for $\cGF$ is a sub-collection of maps 
\begin{equation}
\cM = \{ h_i \colon D(h_i) \to R(h_i) \mid 1 \leq i \leq m\} \subset \cGF
\end{equation}
such that 
\begin{enumerate}
\item  each $h_i \in \cM$ is the restriction of an element $\widetilde{h}_i \in   \cGF$ with    $\overline{D({h_i})} \subset D({\widetilde{h}_i})$
\item  $R({h_i}) \cap R({h_j}) = \emptyset$ for $i \not= j$  (Open Set Condition) 
\item  if $R({h_i}) \cap D({h_j}) \not= \emptyset$ then $R({h_i}) \subset D({h_j})$
\end{enumerate}
\end{definition}
If the maps $h_i$ are linear contractions, then this is an example of what is called an \emph{Iterated Function System} (IFS) in the dynamics  literature. The standard construction of Cantor sets in the unit interval is based on an IFS with two generators. For $q> 1$, computer-generated  simulations of the forward orbits of an IFS can yield a fantastical variety of compact, self-similar   (fractal) regions in $\mR^q$. (For example, see \cite{Falconer2003,FalconerOConnor2007,PJS2004}.) 
There are multiple variations on  Definition~\ref{def-markov}. For example, one can work with compact domains, and in place of the maps $h_i$ use their extensions
\begin{enumerate}
\item $\widetilde{h}_i \colon \overline{D({h_i})} \to \overline{R({h_i})}$.  
\end{enumerate}
The condition (\ref{def-markov}.2) allows  that  $\overline{R({h_i})} \cap \overline{R({h_j})} \ne \emptyset$  for $i \not= j$.  A \emph{Discrete Markov System} is one which satisfies 
\begin{enumerate}\setcounter{enumi}{1}
\item $\overline{R({h_i})} \cap \overline{R({h_j})} = \emptyset$  for $i \not= j$.
\end{enumerate}
This   property  corresponds to the \emph{Strong Open Set Condition} for Iterated Function Systems, which implies that the invariant minimal set for $\cM$ is totally disconnected. (This assumes the generators are strong contractions!) Discrete Markov Systems generate what is called a ``Tits alternative'' or a ``Ping-Pong'' game in the literature \cite{delaHarpe1983,Tits1972}. There is a huge literature on this topic in dynamics.

The \emph{transition matrix} $P_{\cM}$ for a Markov System $\cM$ is the $m \times m$ matrix with entries $\{0,1\}$ defined by 
$P_{ij} = 1$ if $R({h_j}) \subset D({h_i})$, and $0$ otherwise. We say that $\cM$ is \emph{chaotic}    if $P_{\cM}$ is irreducible and aperiodic, so  there  exists $\ell > 0$ such that $P_{\cM}^{\ell}$ is a matrix with all entries positive.  

\begin{definition}\label{def-mms}
A $\cGF$-invariant minimal set $K \in \cBF$ is \emph{Markov} if there is a chaotic Markov sub-pseudogroup $\cM$ such that 
$$K \subset \overline{R(h_1)} \cup \cdots \cup \overline{R(h_k)}$$
and every orbit of $\cM$ in $K$ is dense in $K$.
\end{definition}
Markov minimal sets for codimension-one foliations have been studied by 
Cantwell and Conlon \cite{CantwellConlon1988a,CantwellConlon1989}, 
Inaba {\it et al} \cite{Inaba1986,InabaMatsumoto1990,InabaTsuchiya1992,InabaWalczak96}, 
Matsumoto \cite{Matsumoto1987b,Matsumoto1988}, 
 Walczak \cite{Walczak1996,Walczak2004}, and most recently by 
 the author \cite{Hurder2000b,Hurder2004b} and 
 Rams {\it et al}  \cite{CR2006,GelfertRams2007a,GelfertRams2007b}.
 Together, these works give a partial answer to Conjecture~\ref{conj-miniset}:
 \begin{theorem} Let $\F$ be a codimension-one, $C^2$-foliation of a closed manifold $M$. 
 If  $K$ is a   Markov minimal set for   $\cGF$, then $K$ has Lebesgue measure zero.
 \end{theorem}

 Sacksteder's Theorem implies that an exceptional minimal set for a codimension-one, $C^2$-foliation must have a resilient leaf with linearly contracting holonomy. In particular, this holds for a Markov minimal set $Z = K_{\F} \subset M$, and hence $K \cap \HF \ne \emptyset$.   This conclusion was extended by the author in   \cite{Hurder1991a} to exceptional minimal sets with exponential growth type  for $C^1$-foliations. Moreover, the author proved  the following partial answer to Conjecture~\ref{conj-miniset} for $C^r$-foliations:
 \begin{theorem}[Hurder, \cite{Hurder2005}]\label{thm-excmeaszero}
  Let $\F$ be a codimension-one, $C^r$-foliation of a closed manifold $M$, for $r > 1$. 
 Let $K \in \cBF$ be an exception minimal set for   $\cGF$. Then $K \cap \HF$ has Lebesgue measure zero.
 \end{theorem}
 
 These results highlights the importance of the parabolic points in a minimal set, 
 $$K_{\cP}    = K   \cap (\EF \cup \PF) =  K \setminus (K \cap \HF)$$
  The general form of Conjecture~\ref{conj-miniset} for $C^r$-foliations is thus equivalent to showing:
 \begin{conjecture}
   Let $\F$ be a codimension-one, $C^r$-foliation, for $r > 1$. 
 Let $K \in \cBF$ be an exception minimal set for   $\cGF$. Then the set of parabolic points $K_{\cP}$ has Lebesgue measure zero.
 \end{conjecture}
 
There is a close connection between the existence of a Markov minimal set $\cM$ for $\cGF$ and $h(\cGF) > 0$. 
\begin{proposition}[Th\'{e}or\`{e}me~6.1, \cite{GLW1988}] 
   Let $\F$ be a codimension-one, $C^1$-foliation of a closed manifold $M$. 
If   $\cGF$ has    a   Markov minimal $K$, then $h(\cGF) >  0$. More precisely, we have $h(\cGF, K) > 0$.
\end{proposition}
Th\'{e}or\`{e}me~6.1 in \cite{GLW1988}  gives a proof of the converse for $C^2$-foliations of foliated bundles, using the Poincar\'{e}-Bendixson theory of levels (see also Theorem~13.5.3, \cite{CanCon2000}, and Theorem~3.6.1, \cite{Walczak2004}).
 The author showed that,  in fact, the full converse holds for $C^1$-foliations
\begin{theorem}[Theorem~1.1, \cite{Hurder2000b}; see also  Theorem~4.6.1, \cite{Walczak2004}]
   Let $\F$ be a codimension-one, $C^1$-foliation of a closed manifold $M$. 
Assume  $h_{loc}(\cGF,x) >  0$. Then for every open neighborhood, $x \in U \subset \cT$,  $\cGF$ has a chaotic, discrete Markov minimal set $K \subset \overline{U_{\cR}}$.
\end{theorem}

The proof of Th\'{e}or\`{e}me~6.1, \cite{GLW1988}  required the full theory of Poincar\'e-Bendixson Theory for $C^2$-foliations; they showed that if there is no resilient leaf, then the structure theory for codimension-one, $C^2$ foliations implies that $h(\F) = 0$. Thus, it used global methods to prove a global result. 

 The proof of Theorem~1.1, \cite{Hurder2000b} shows that  $h_{loc}(\cGF,x) > 0$ yields  normally expansive behavior along segments of orbits in every open neighborhood $x \in U \subset \cT$, which yields hyperbolic fixed-points which capture points in their   orbit, and thus generate homoclinic (resilient) orbit behavior. In contrast to the methods of \cite{GLW1988}, the techniques are $C^1$, local,  and mostly naive.

   Finally, we consider the relation between the Godbillon-Vey class $GV(\F)$ and  foliation dynamics in codimension-one. Recall that Theorem~\ref{thm-gv} showed that $GV(\F) \ne 0$ implies the set $\HF$ has positive Lebesgue measure. 

\begin{theorem} [Hurder \& Langevin, \cite{HLa2000}]  \label{thm-HL}
Let  $\F$ be a codimension-one, $C^1$-foliation  such that $\HR$ has positive Lebesgue measure.  
For every open set $U \subset \cT$ such that $U \cap \HR$ has positive Lebesgue measure,   
 $\cGF$ has a discrete Markov minimal set $K \subset \overline{U_{\cR}}$. 
\end{theorem}
The proof actually shows more:
\begin{theorem} [Hurder, \cite{Hurder2008b}]  
Let  $\F$ be a codimension-one, $C^1$-foliation   such that $\HR$ has positive Lebesgue measure.  
Then $h_{loc}(\cGF, x) > 0$ for almost every $x \in \HR$.
\end{theorem}
Hence we conclude:
\begin{corollary}   
Let  $\F$ be a $C^2$-foliation of codimension-one. Given $E \in \cBF$ such that  $GV(\F)|E \ne 0$, then there exists a set  $K \subset E$  of positive Lebesgue measure such that   $h_{loc}(\cGF, x) > 0$ for every $x \in K$.
\end{corollary}

Combining Theorems~\ref{thm-gv}, \ref{thm-excmeaszero} and \ref{thm-HL} we obtain:
\begin{theorem} [Hurder, \cite{Hurder2008b}]  \label{thm-grand}
Let  $\F$ be a $C^2$-foliation of codimension-one. Let $E \in \cBF$ be an exceptional local minimal set for $\cGF$. Then the Godbillon measure of $E$ vanishes. Hence, if $GV(\F) \ne 0$, then   there is an open saturated subset $U \subset M$ with
\begin{itemize}
\item $U$ contains the support of the cohomology class $GV(\F)$;
\item  $U$ contains a dense collection of discrete Markov minimal sets;
\item $\F | U$ is expansive.
\end{itemize}
\end{theorem}
 This yields a positive  solution to the Conjecture (page 239, \cite{CanCon2003}). The $C^2$-hypothesis is used in the 
  proof of Theorem~\ref{thm-grand} in three places:   to guarantee that the Godbillon-Vey class is defined;   in the proof of Theorem~\ref{thm-excmeaszero}; and finally   the existence of the open set $U$ requires the Poincar\'e-Bendixson Theory for $C^2$-foliations.

\section{Structure of minimal and transitive sets} \label{sec-minimal}

We conclude this survey with a discussion of some dynamical concepts and open problems related to the results discussed above, and  which suggest promising areas of investigation. One of the main obstacles to developing a full Pesin Theory for groupoid dynamics, is the general absence of invariant measures. For maps, the existence of an invariant measure has strong implications for the recurrence of its orbits. For a foliation, recurrence can be obtained by assuming it; that is, by studying the dynamics restricted to  minimal and transitive sets.  

Recall that a closed subset $K \in \cBF$ is minimal if for every $x \in K$ the orbit $\cO(x)$ is dense in $K$, while $K$ is transitive if for some $x \in K$, the orbit $\cO(x)$ is dense in $K$. The Poincar\'{e}-Bendixson Theory for codimension-one, $C^2$-foliations is focused on the properties of the (local open) minimal sets for the foliation. 

For dynamical systems on manifolds of dimension greater than $1$, the transitive sets are perhaps more important. This is best illustrated for the case where       $f \colon N^q \to N^q$ is an \emph{Axiom-A} diffeomorphism, so $q > 1$. By definition of Axiom-A, the restriction of $f$ to the non-wandering set $\Omega(f)$  is hyperbolic, and the periodic orbits of $f$ are dense in $\Omega(f)$. Thus, $\Omega(f)$ is not minimal, but $\Omega(f)$ is the finite union of its basic sets, which are closed, disjoint invariant   subsets on which $f$ is transitive.  The dynamics of $f$ are captured by its behavior near each basic set. Within $\Omega(f)$, the periodic points   play an important role, but their behavior is predicted by the Markov coding for the dynamics of $f$ on $K$  \cite{Bowen1978,Smale1967,Smale1980}. 
  For foliations with codimension $q > 1$, a similar conclusion is likely.

\begin{problem} 
Is there a theory of Axiom-A foliations? That is, can one impose sufficient hyperbolicity and other hypotheses on the non-wandering set $\Omega(\F)$ of a $C^2$-foliation such that the dynamics of $\F$ can be ``classified''?
 \end{problem}

The   Poincar\'{e}-Bendixson Theory for codimension-one foliations suggests that one should study the asymptotic properties of leaves in higher codimension.  Marzougui and Salhi introduced in  \cite{MarzouguiSalhi2003} a theory of levels for $C^1$-foliations, under the restriction that the foliations admit a transverse foliation. Their approach is based on the study of the open local minimal sets, and their main result is a structure theorem for the dynamics,  analogous to the conclusion for codimension-one foliations.

 In analogy, this suggests the definition:   a minimal (or possibly transitive) compact subset $E_0 \in \cBF$ has \emph{level zero}. Then define inductively, that a point $x \in \cT$ has \emph{level $k$} if the $\omega$-limit set $\omega(x)$ is a union of closed subsets with level less than $k$.  It is not   clear that this decomposition of the orbits of $\cGF$ carries similar import to the case of codimension-one.  
\begin{problem} 
How do the dynamical properties of $\F$, say as given by its decompositions into the sets at the start of section~\ref{sec-secclasses2}, and the theory of levels, either as in Marzougui and Salhi \cite{MarzouguiSalhi2003} or as above, determined one another?
For example,   how do the dynamical properties of the sets $E_0$ of level zero, influence the dynamics of an orbit $x \in \cT$ for which $E_0 \subset \omega(x)$.
\end{problem}

 In the case of codimension-one, the following problem remains open:
 \begin{problem} Let $\cM$ be a Markov system on $\mR$. Does there exists a $C^r$-foliation of a closed manifold $M$ such that $\cM$ is a Markov sub-pseudogroup of $\cGF$?
 \end{problem}
 In the case where $\cM$ is a discrete Markov system, Cantwell and Conlon prove this in Section~9,  \cite{CantwellConlon2002}. (Conlon gave a more detailed proof   in unpublished notes  \cite{Conlon2000}.) Their method does not seem to apply, however, to the case where $\cM$ does not satisfy the strong open set condition. 
 
Section~6.1 of \cite{Hurder2006b} discusses other open questions about the codimension-one case. They may be summarized by asking:
\begin{problem} 
Let $\F$ be a codimension-one, $C^r$-foliation  for $r \geq 1$. 
Characterize the compact minimal sets  $K \in \cBF$  for $\cGF$ which are not Markov.  
 \end{problem}

Next, assume that  $\F$ has codimension $q > 1$. If $\F$ is Riemannian, then the closure of each leaf $L$ is a minimal set, and moreover is a locally homogeneous space. 
\begin{problem} 
Let $K \in \cBF$ be a compact minimal set with   $K \subset \ER$.  Is the restriction $\cGF | K$ equicontinuous? Must $K_{\F}$ be a connected submanifold of $M$? 
 \end{problem}
In other words, does recurrence of the orbits in $K$ force   an invariant Borel Riemannian metric on the normal bundle, to be a continuous metric?

For parabolic foliations, we saw in Theorem~\ref{thm-solenoids} that there exists minimal sets for $\F$ which are generalized solenoids by construction. We know by Theorem~\ref{thm-hypfixpt2} 
that $K \cap \HR = \emptyset$ and that the restricted entropy $h(\cGF, K) = 0$.

\begin{problem} 
Let  $K \in \cBF$ be a compact minimal set such that $\cGF^K$ is distal. If $K$ is a Cantor set, is the saturation $K_{\F}$   a generalized solenoid?
 \end{problem}

  While it is totally speculative, one can ask if  there   exists a general classification   for parabolic minimal sets? 
  
\begin{problem} 
Suppose that $K \in \cBF$ is a compact minimal set such that  $\cGF^K$ is distal, or more generally just parabolic. 
Does there exists an analogue of the Furstenberg  structure theory for distal actions of a single transformation \cite{Ellis1978,Furtsenberg1963,Lindenstrauss1999,Zimmer1976}, which applies to $\cGF^K$? 
 \end{problem}

 For codimension $q > 1$, another phenomenon arises in the study of closed  invariant sets for $\cGF$. An exceptional  minimal set $K \in \cBF$ has no interior and is not discrete, hence $K \subset \mR$   implies that $K$ is totally disconnected and perfect. That is, $K$  is a Cantor set. However, for $K \subset \cR^q$ with $q > 1$, the set $K$ may have no interior, and yet not be totally disconnected. 
 \begin{definition}\label{def-exotic}
A minimal set   $K \in \cBF$ for $\cGF$ is an \emph{exotic minimal set} if $K$ has no interior, is perfect and is connected.
 \end{definition}
 For example, the invariant minimal sets for a Kleinian group of the second kind acting on the sphere $\mS^q$ has minimal set homeomorphic to a Sierpinski space, which is exotic. These examples are realized as minimal sets for smooth foliations using the standard suspension construction.   
 
Using the generalized suspension construction in \cite{BHS2006}, the authors constructed many families of  examples of smooth foliations for which the minimal sets are exotic subsets of $\mT^q$.  In fact, the torus $\mT^q$ can be replaced by any nil-manifold \cite{BisHurder2008}. 
 
    Exotic minimal sets also arise in the action of certain word-hyperbolic groups on their boundary at infinity \cite{KapovichBenakli2002}. It is not clear whether these examples can be realized as minimal sets for foliation pseudogroups.
 
There are few techniques developed for constructing foliations with prescribed minimal sets.  On the other hand, there is a massive literature for constructing continua defined by inverse limit constructions in the spaces $\mR^q$ 
with various remarkable properties \cite{KY1995}. The solenoid discussed previously is probably the simplest of these constructions, so one might expect other continua can be realized. We conclude with three open questions:

 \begin{problem} 
Which a compact continua in $ \mR^q$ are homeomorphic to a minimal (or transitive) invariant set for a pseudogroup $\cGF^r$, $r \geq 1$,  of a foliation $\F$ on a closed manifold $M$?
\end{problem}
 
  \begin{problem} 
Give an example of    $K \in \cBF$ which is an exotic minimal set for a $C^r$-pseudogroup $\cGF$, $r > 1$,    such that  $K$ has positive Lebesgue measure. If $K \subset \HR$, must $K$ have Lebesgue measure zero?
 \end{problem}
 
 \begin{problem} Let $\F$ be a $C^2$-foliation of codimension $q > 1$, and suppose that 
  $K \in \cBF$   is an exceptional minimal set. Show that $h_I(K) = 0$ for all Weil measures.
 \end{problem}


 \bibliographystyle{amsalpha}

\end{document}